%% file: main.tex
\documentclass{article}
\input{header}
\input{macro}
\begin{document}

\maketitle
\begin{abstract}
    In this paper, we introduce a \emph{Homogeneous Second-Order Descent Method} (HSODM) motivated from the homogenization trick in quadratic programming. The merit of homogenization is that only the leftmost eigenvector of a gradient-Hessian integrated matrix is computed at each iteration. Therefore, the algorithm is a single-loop method that does not need to switch to other sophisticated algorithms and is easy to implement. We show that HSODM has a global convergence rate of $O(\epsilon^{-3/2})$ to find an $\epsilon$-approximate second-order stationary point, and has a local quadratic convergence rate under the standard assumptions. The numerical results demonstrate the advantage of the proposed method over other second-order methods.
\end{abstract}

\section{Introduction}

In this paper, we consider the following unconstrained optimization problem:
\begin{equation}\label{Prob:main}
    \min_{x\in \mathbb{R}^n} f(x),
\end{equation}
where \(f:\mathbb{R}^n \mapsto \mathbb{R}\) is
a twice continuously differentiable function and \( f_{\inf}:=\inf f(x) > -\infty\).
Given a tolerance level $\epsilon > 0$, we aim to find an $\epsilon$-approximate second-order stationary point (SOSP) $x$ satisfying
\begin{subequations}
    \begin{align}
        \label{eq.approxfocp} & \|\nabla f(x)\|\leq O(\epsilon),                         \\
        \label{eq.approxsocp} & \lambda_1 (\nabla^2 f(x)) \geq \Omega(-\sqrt{\epsilon}),
    \end{align}
\end{subequations}
where $\lambda_1(A)$ denotes the smallest eigenvalue of matrix $A$.
When $f$ is nonconvex, it has been shown that the gradient descent (GD) method finds an $\epsilon$-approximate first-order stationary point satisfying \eqref{eq.approxfocp} in $O(\epsilon^{-2})$ iterations under the standard $L$-Lipschitz continuous gradient condition.
If the second-order condition \eqref{eq.approxsocp} is further required, first-order methods may fail, and a common practice is to consider second-order methods, that is, some variants of Newton's method \cite{connTrustRegionMethods2000}.

At each iteration, the Newton-type methods usually construct the second-order approximation at the current iterate $x_k$, and then compute the direction $d_{k}$ for the update. For example, Newton's method utilizes the following quadratic approximation:
\begin{equation}\label{eq.mk-Newton}
    d_{k} = \argmin_{d\in\real^n} \,\, m_k(d):=\gk^T d+\frac{1}{2}d^T \Hk d,
\end{equation}
where $\gk = \nabla f(\xk) $ and $ \Hk = \nabla^2 f(\xk)$.
In the nonconvex case, despite excellent performance in practice, \citet{cartis_complexity_2010} showed that Newton's method, perhaps surprisingly, has a worst-case complexity of $O(\epsilon^{-2})$ similar to that of GD. Therefore, some advanced techniques are needed to improve the convergence performance of Newton's method. \citet{nesterov_cubic_2006} introduced the cubic regularization (CR) and consider the following subproblem:
\begin{equation}\label{eq.arc}
    \dk^\mathsf{CR} = \arg\min_d \,\, m^{\mathsf{CR}}_k(d):=\gk^T d+\frac{1}{2}d^T \Hk d + \frac{\sigma_k}{3}\|d\|^3,
\end{equation}
where $\sigma_k > 0$. They showed that the cubic regularized Newton's method has an improved iteration complexity of $O(\epsilon^{-3/2})$.
\citet{cartis_adaptive_2011,cartis_adaptive_2011-1} later proposed an adaptive and inexact version of cubic regularization (ARC) with the same complexity.
Before the appearance of CR, a widely used classic algorithm is the trust-region (TR) method. It computes the update direction based on the same model function as Newton's method, but restrains it within the pre-specified trust-region radius $\Delta_k$, and accepts it if the corresponding acceptance ratio $\rho_k$ exceeds some threshold \cite{connTrustRegionMethods2000}:
\begin{subequations}\label{eq.tr}
    \begin{align}
        \label{eq:tr-subproblem}\dk^\mathsf{TR} ~ & =\arg \min_{\|d\|\le \Delta_k} m_k(d),              \\
        \label{eq:tr-tradi-rhok}\rho_k ~          & := \frac{f(x_k + \dk) - f(x_k)}{m_k(\dk) - m_k(0)}.
    \end{align}
\end{subequations}

However, it is more challenging to establish the improved $O(\epsilon^{-3/2})$ iteration complexity in this way. To our best knowledge, \citet{yeSecondOrderOptimization2005} provided the earliest $O(\epsilon^{-3/2})$ trust-region method by a fixed radius strategy. Recently, \citet{curtis_trust_2017} pointed out that the classical TR method (based on \eqref{eq.tr}) fails in satisfying the sufficient decrease property required to obtain the $O(\epsilon^{-3/2})$ complexity rate because it uses the classical $\rho_k$-based acceptance rule and linearly updated radius.

To overcome this issue, they developed an algorithm named TRACE \cite{curtis_trust_2017,curtisWorstCaseComplexityTRACE2023}, which achieves the desired complexity result but has a sophisticated rule of expanding and contracting $\Delta_k$ due to the nonlinearity between $\|\dk^\mathsf{TR}\|$ and the dual solution of the problem \cref{eq:tr-subproblem}. This complexity bound can also be achieved via a line search Newton CG framework proposed in
\citet{royer_complexity_2018}. Their algorithm alternates between Newton and regularized Newton steps based on the smallest eigenvalue of the Hessian $\Hk$, and the stepsize is chosen under a similar acceptance rule used in \cite{cartis_adaptive_2011,curtis_trust_2017}. Nevertheless, all the above methods solve Newton systems, where the cost of $O(n^3)$ is typical. It is possible to find inexact solutions with better complexity performance. In that sense, many classical algorithms are open for improvement with techniques such as negative curvature oracles and conjugate gradient method \cite{royer_complexity_2018,royer_newton-cg_2020,curtis2021newton-cg}.

\subsection{Our contribution}

Motivated by the homogenization technique to obtain semidefinite relaxations of quadratic programming \cite{sturm_cones_2003,ye_new_2003}, we propose a homogenized version of the local quadratic approximation $m_k(d)$.
We show that the resulting problem is essentially an eigenvalue problem and can be solved by the random-starting Lanczos algorithm \cite{kuczynski_estimating_1992}, which allows a dimension-independent complexity of $\tilde O(n(n+1)\epsilon^{-1/4})$ with high probability.
We demonstrate that the leftmost eigenvalue of the homogenized matrix is always negative; namely, the ``homogenized negative curvature" exists even when the original Hessian is near positive semidefinite. Similar to the gradient descent method, where a first-order stationary point is reached by moving along the negative gradient direction, we can attain a second-order stationary point by \emph{exclusively} moving along the direction corresponding to the homogenized negative curvature.

Secondly, we propose a new second-order method called the Homogeneous Second-Order Descent Method (HSODM) (\cref{alg.main alg})  with the homogenized quadratic model as subproblems. We offer two stepsize strategies to utilize the homogenized negative curvature, including the fixed-radius strategy and a simple backtracking line search method. Our method achieves a better iteration complexity of $O(\epsilon^{-3/2})$ to converge to a SOSP than the $O(\epsilon^{-2})$ complexity of the standard trust-region method \cite{curtis_concise_2018} and the negative-curvature based method \cite{curtis_exploiting_2019}. Accounting for the subproblems, it requires  $\tilde O((n+1)^2\epsilon^{-7/4})$ arithmetic operations.
In sharp comparison to \cite{carmon_accelerated_2018,agarwal_finding_2017,jin_how_2017,royer_complexity_2018}, HSODM only relies on the homogenized model and \emph{does not} alternate between different subroutines. The algorithm is elegant in a simple form and believed to be highly favorable to practitioners. To make a clear comparison, we provide the following \cref{tab:comparison} that includes the algorithms with the state-of-the-art complexity results.  Note that ARC \cite{cartis_adaptive_2011-1} and TRACE \cite{curtisWorstCaseComplexityTRACE2023} require Newton-type equations, from cubic regularized problems and trust-region subproblems, respectively. Both can be solved by applying matrix factorizations (in $O(n^3)$) with a suitable parameter search procedure in $O(n^2\log(1/\epsilon))$. To incorporate inexact subproblem solutions, the methods in \cite{royer_complexity_2018,royer_newton-cg_2020,curtis2021newton-cg} switch between the conjugate gradient method for linear equations and a randomized Lanczos method for extreme eigenvalue problems, so that the complexity rates can be improved to $\tilde O(n^2 \epsilon^{-1/4})$. For HSODM, only extreme eigenvalue problems are needed.

\begin{table}[h]
    \centering
    \footnotesize
    \caption{A brief comparison of several second-order algorithms. Here, $p \in (0,1)$ represents the failure probability of the randomized Lanczos method.  In the last column, we use \enquote{E} for the extreme eigenvalue problem and \enquote{N} for Newton-type equation.}
    \label{tab:comparison}
    \begin{tabular}{lccc}
        \toprule
        Algorithm                                                         & Iteration Complexity    & Subproblem Complexity                                 & Oracle(s) \\
        \midrule
        ARC \cite{cartis_adaptive_2011-1}                                 & $ O( \epsilon^{-3/2}) $ & $O(n^3 + n^2 \log(1/\epsilon))$                       & N         \\
        TRACE \cite{curtis_trust_2017,curtisWorstCaseComplexityTRACE2023} & $ O( \epsilon^{-3/2}) $ & $O(n^3 + n^2 \log(1/\epsilon))$                       & N         \\
        \cite[Algorithm 4.1]{curtis2021newton-cg}                         & $ O(\epsilon^{-3/2}) $  & $O(n^2 \epsilon^{-1/4} \log(n / p \epsilon))$         & N \& E    \\
        Newton-CG \cite{royer_complexity_2018,royer_newton-cg_2020}       & $O(\epsilon^{-3/2})$    & $O(n^2 \epsilon^{-1/4} \log(n / p\epsilon))$          & N \& E    \\
        HSODM                                                             & $ O(\epsilon^{-3/2})$   & $O((n+1)^2 \epsilon^{-1/4} \log(n(n+1) / p\epsilon))$ & E         \\
        \bottomrule
    \end{tabular}
\end{table}

Finally, the numerical results of the proposed method are also encouraging. In particular, two variants of HSODM outperform the standard second-order methods, including the classical trust-region method and the cubic regularized Newton method in the CUTEst dataset.

\subsection{Related works}

There is a recent trend in the study of
improved first-order algorithms \cite[]{carmon_accelerated_2018,agarwal_finding_2017,jin_how_2017} for $\epsilon$-approximate SOSP. Thus, this type of algorithm can serve as a scalable alternative to second-order ones.
Notably, some of these algorithms also enable faster first-order convergence to \eqref{eq.approxfocp} in $\tilde O(\epsilon^{-7/4})$ function and gradient evaluations. The basic idea is to extend Nesterov's accelerated gradient descent method (AGD) \cite[]{nesterovLecturesConvexOptimization2018} to the nonconvex case. This is achieved by properly embedding second-order information to make the AGD maintain its theoretical property in convex and semiconvex cases.
For example, \citet{carmon_accelerated_2018} applied Hessian-vector products and randomized Lanczos methods to explore the negative curvature (NC) {(we will define this formally in \eqref{eq.nc.suff})}, which is then used as a descent direction; otherwise, $f$ becomes locally semiconvex and AGD is invoked to solve the subproblem.
The later work in \cite[]{jin_how_2017,xu_first-order_2018} also requires NC but avoids Hessian-vector products, and the complexities remain the same. Beyond using NC, \citet[]{agarwal_finding_2017} achieved the same complexity bound by applying fast matrix inversion to cubic regularized steps. Recently, \citet{li_restarted_2022} introduced a restarted AGD that drops the logarithmic term $O(\log(\epsilon^{-1}))$ in the complexity bound if the solution is required to only satisfy the first-order condition, but it also losses second-order guarantees.
To make AGD work in a comfort zone, these algorithms create sophisticated nested loops that may be difficult to implement and tune.
Nevertheless, they are designed to be less ``dimension-dependent'' than pure second-order methods such as \cite{nesterov_cubic_2006,cartis_adaptive_2011} and are suitable for large-scale problems in theory.

Coming back to the second-order methods,  \citet{royer_complexity_2018} separated their method into two cases {if NC is absent}. In one case the smallest eigenvalue $\lambda_1(\Hk) > -\sqrt{\epsilon}$, regularized Newton step is used to provide the descent step. In the other case, when $\lambda_1(\Hk) > \sqrt \epsilon$ is certified, it turns to the ordinary Newton step. Therefore, in the worst case, this method must solve an eigenvalue problem and a Newton step at one iteration. It is unclear if one can unify these procedures as a whole. Recently, \citet{mishchenko2021regularized} proposed the Gradient Regularized Newton method for convex minimization problems. The subproblem at each iteration is simpler than that of the Cubic Regularized Newton method. Later, \citet{gratton2023yet} generalized the Gradient Regularized Newton method to exploit negative curvature for nonconvex optimization problems. However, their method alternates between regularized Newton and negative-curvature steps.

For trust-region methods \eqref{eq.tr}, $\lambda_1(\Hk) \le -\sqrt \epsilon$ implies that the Lagrangian dual variable is at least in the order of $\sqrt{\epsilon}$. When the curvature is properly utilized, it also implies a $\Omega(\epsilon^{3/2})$ progress as long as the stepsize is carefully selected. This fact can be easily recognized by using optimality conditions (for example, see \cite{yeSecondOrderOptimization2005,curtis_concise_2018}). Moreover, it remains true even when the subproblems are solved inexactly or suboptimally in some subspace \cite{cartis_adaptive_2011-1,curtisWorstCaseComplexityTRACE2023,zhang_drsom_2023}. \citet{curtis2021newton-cg} further proposed a trust-region method that does not alternate between steps but rather solves a slightly perturbed trust-region subproblem. {For fixed-radius strategies \cite{yeSecondOrderOptimization2005,zhang_drsom_2023},  the algorithm safely terminates if it is nearly convex.}

The situation is different for adaptive methods. Since the trust-region method uses an acceptance ratio $\rho_k$ in \eqref{eq.tr} and adjusts the radius linearly, a step may become too small with respect to the dual variable. A workaround can be found in \cite{curtis_trust_2017,curtisWorstCaseComplexityTRACE2023} with a delicate control over the progress of the function value and the gradient norm:
$$
    f_k-f_{k+1} \geq \Omega(\|\dk\|^3) \text { and }\left\|\dk\right\| \geq \Omega(\|g_{k+1}\|^{1 / 2}).
$$
Similar conditions are also needed in the analysis of cubic regularization methods \cite{cartis_adaptive_2011-1}. However, these adaptations can be less straightforward to understand, implement, and adjust.

In addition, our work is also related to solving trust-region subproblems by eigenvalue procedures \cite{rojas_new_2001}, which use the same $(n+1)$-dimensional symmetric matrix.  The idea is later extended to solve cubic regularized subproblems or generalized trust-region subproblems; see, for example, \cite{generalized-trs,felix-arc-generalized}. Both papers introduce matrix pencils that raise the dimension to $2(n+1)$ without providing the convergence analysis. While the aforementioned works mainly focus on solving subproblems, we use the homogenized matrix in a generic method that finds stationary points of a generic optimization problem. We also provide a complexity analysis of its global and local convergence. In addition, the matrices they construct have larger dimensions than ours, which brings more computational cost when solving the corresponding eigenvalue problem.

\subsection{Notations, assumptions, and organization of the paper}
In this subsection, we introduce the notations and assumptions used throughout the paper.

Let  \(\|\cdot\|\) be the standard Euclidean norm in space $\mathbb{R}^n$. Denote $B(x,R) = \{y: \|y-x\| \le R\}$ to be the closed ball with radius $R$ centered at $x$. For a matrix \(A \in \mathbb{R}^{n\times n}\), \(\|A\|\) represents the induced \(\ell_2\) norm, and $\lambda_1(A),\lambda_2(A),...,\lambda_{\max}(A)$ denotes its distinct eigenvalues in ascending order. For $n > 0$, $I_{n}$ denotes the $n$-dimensional identity matrix; we omit $n$ if it is clear from the context. At some iterate $\xk$, we denote \(\gk = \nabla f(\xk)\) and \(\Hk = \nabla^{2} f(\xk)\) for simplicity. We use order notation $O,\;\Omega,\; \Theta$ in the usual sense, while $\tilde O$ hides the logarithmic terms with respect to $O$. In particular, given two constants $A$ and $B$, we say $A = O(B)$ if there exists a constant $c > 0$ such that $A \le c\cdot B$, and $A = \Omega(B)$ if there exists a constant $c > 0$ such that $A \ge c\cdot B$. We say $A = \Theta(B)$ if $A = O(B)$ and $A = \Omega(B)$. We use $[a;b]$ (resp., $[a,b]$) to denote vertical (resp., horizontal) concatenation of arrays or numbers. For a vector $a\in \real^n$ and $0 \le j \le n$, we let $a_{[1:j]}$ be the first $j$ entries of $a$.

The rest of the paper is organized as follows.
In \cref{sec.method}, we briefly describe our approach based on the homogenized quadratic model. By solving the homogenized model as an eigenvalue problem, the corresponding HSODM is introduced in \cref{alg.main alg}.  In \cref{sec.global} and \cref{sec.local}, we give analyses of the global and local convergence of HSODM. Our results indicate that HSODM has a global complexity of $O(\epsilon^{-3/2})$ for an $\epsilon$-approximate second-order stationary point. If one does not early terminate the algorithm,  it converges at a local quadratic rate. We address the inexactness in HSODM in \cref{sec.inexact}, where a Lanczos method with skewed initialization is introduced to utilize the Ritz approximation to homogeneous curvature. In \cref{sec.experiment}, we demonstrate the effectiveness of our method by providing fruitful computational results in the CUTEst benchmark compared to other standard second-order methods.

\section{The Homogenized Quadratic Model and A Second-Order Descent Method}\label{sec.method}

\subsection{Motivation of homogenization}

Many optimization methods for nonconvex optimization use the Negative Curvature of the Hessian matrix. In particular, given an iterate $x_k$, it is often of interest to determine if there exists $\xi_k \in \real^n$ such that
\begin{equation}\label{eq.nc.suff}
    \mathcal R_k(\xi_k) := \frac{\xi_k^T\Hk \xi_k}{\|\xi_k\|^2} \le -\sqrt{\epsilon},
\end{equation}
for some tolerance $\epsilon > 0$, as it implies that $\lambda_{\min }\left(H_k\right) \leq-\sqrt{\epsilon}$.
    {Such a $\xi_k$ is referred to as the direction associated with negative curvature.}
Computationally, it is known that $\xi_k$ can be found at the cost of $\tilde O\left(n^2\cdot\epsilon^{-1/4}\right)$ arithmetic operations, using a randomized Lanczos method \cite{kuczynski_estimating_1992}.
When facilitating this direction with a proper stepsize $\eta$, the function value must decrease by $\Omega(\epsilon^{3/2})$ under second-order Lipschitz continuity. This property is widely used in the negative-curvature-based first-order methods \cite{carmon_accelerated_2018,jin_how_2017}.
However, if \eqref{eq.nc.suff} is invalid, one must switch to other subroutines \cite{carmon_accelerated_2018,agarwal_finding_2017,jin_how_2017,royer_complexity_2018}, which complicates the iteration procedure and thus is hard for efficient implementation and parameter tuning.

To alleviate this issue, we apply the homogenization trick (e.g., see \cite{ye_new_2003,sturm_cones_2003}) to the second-order Tayler expansion \eqref{eq.mk-Newton} at $\xk$:
\begin{align}
    \label{eq.homo.taylor}    {t^2} \left(m_k(d) - \frac{1}{2}\delta\right)
                        & = t^2\left(\gk ^T({v/t}) + \frac{1}{2} ({v/t})^T\Hk({v/t}) - \frac{1}{2} \delta\right) & (d:=v/t)~ \\
    \label{eq.motivate} & = t\cdot \gk^Tv + \frac{1}{2} v^T\Hk v - \frac{1}{2}{\delta} t^2
    = \frac{1}{2} \begin{bmatrix}
                      v \\ t
                  \end{bmatrix}^T
    \begin{bmatrix}
        \Hk   & \gk     \\
        \gk^T & -\delta
    \end{bmatrix}
    \begin{bmatrix}
        v \\ t
    \end{bmatrix},
                        & F_k := \begin{bmatrix}
                                     \Hk   & \gk     \\
                                     \gk^T & -\delta
                                 \end{bmatrix}.
\end{align}
The second equation is called \emph{homogenized quadratic model}.
One nice property of the \emph{homogenized matrix} $F_k$ is that: even if $\Hk$ is positive definite, $F_k$ is still indefinite, and thus the ``homogenized negative curvature'' can be computed from this $(n+1)$-dimensional lifted matrix.
To make a connection to the Rayleigh quotient given in \eqref{eq.nc.suff}, we impose a ball constraint $\|[v;t]\|\le 1$ and so \eqref{eq.motivate} is bounded. Furthermore, if we take $d=v/t$, the homogenized model and the second-order approximation \eqref{eq.mk-Newton} scaled by $t^2$ are equivalent up to some constant, i.e., $-\delta/2$.

% In this case, it is necessary to place a compact set for $[v;t]$ otherwise the problem is unbounded. We then recognize its correspondence to the eigenvalue problem if simply a ball constraint $\|[v;t]\|\le 1$ is concerned, since \eqref{eq.motivate} is simply a homogeneous indefinite quadratic problem. 

\subsection{Overview of the method}
We present the HSODM in \cref{alg.main alg}.
The rest of this paper discusses the method that uses the ``homogenized'' matrix in the iterates.  We formally define the homogenized quadratic model as follows. Given an iterate $\xk \in \mathbb{R}^n$, let $\psi_k(v, t; \delta)$ be the homogenized quadratic model,
\begin{equation}\label{eq.homoquadmodel}
    \psi_k(v, t; \delta) := ~ \begin{bmatrix}
        v \\ t
    \end{bmatrix}^T
    \begin{bmatrix}
        \Hk   & \gk     \\
        \gk^T & -\delta
    \end{bmatrix}
    \begin{bmatrix}
        v \\ t
    \end{bmatrix},~ v\in \mathbb{R}^n, t \in \mathbb{R},
\end{equation}
where $\delta \geq 0$ is a predefined constant. For each iteration, the HSODM minimizes the model at the current iterate $\xk$, i.e.,
\begin{equation}
    \label{eq.homo subproblem}
    \begin{aligned}
        \min_{\|[v; t]\| \le 1} \psi_k(v, t; \delta).
    \end{aligned}
\end{equation}
Denote the optimal solution of problem \eqref{eq.homo subproblem} as $[v_k; t_k]$. As the subproblem \eqref{eq.homo subproblem} is essentially an eigenvalue problem, and $[v_k; t_k]$ is the eigenvector corresponding to the smallest eigenvalue of $F_k$. Therefore, we can solve this subproblem using an eigenvector-finding procedure, see \cite{carmon_accelerated_2018, kuczynski_estimating_1992, royer_complexity_2018}.

After solving \eqref{eq.homo subproblem}, we construct a descent direction $d_k$ based on this optimal solution $[v_k;t_k]$ and carefully select the stepsize $\eta_k$ to ensure sufficient decrease.  According to \eqref{eq.homo.taylor}, $d_k = v_k/t_k$ would be the natural choice. However, the {extremal} case of $t_k = 0$ could make $d_k$ tend to infinity. Intuitively, if $|t_k|$ is sufficiently small,
it means that the Hessian matrix $H_k$ dominates the homogenized model, and thus we choose the truncated direction $v_k$ directly (\cref{line.exact.larga}). Otherwise, the predefined parameter $-\delta$ becomes significant, and we choose $v_k/t_k$ as the descent direction instead (\cref{line.exact.largb}).  We use $\sqrt{1/(1+\Delta^2)}$ and $\nu$ as the thresholds of $|t_k|$ to determine whether it is sufficiently small.
For the stepsize rule, we provide two strategies for selecting the stepsize: the first is to use line search to determine $\eta_k$, and the second is to adopt the idea of the fixed-radius trust-region method \cite{luenberger_linear_2021,zhang_drsom_2023} such that $\|\eta_kd_k\| = \Delta$, where $\Delta$ is some pre-determined constant. By iteratively performing this subroutine, our algorithm will converge to an $\epsilon$-approximate SOSP.

\begin{minipage}[t]{0.95\linewidth}
    \begin{algorithm}[H]
        \caption{Homogeneous Second-Order Descent Method (HSODM)}\label{alg.main alg}
        \KwData{initial point $x_1$, $\nu \in (0, 1/2)$, $\Delta  = \Theta(\sqrt{\epsilon})$}
        \For{$k = 1, 2, \cdots$}{
            Solve the subproblem \eqref{eq.homo subproblem}, and obtain the solution $[v_k; t_k]$\;
            \uIf(\tcp*[f]{small value case}){$|t_k| > \sqrt{1/(1+\Delta^2)}$}{
                $\dk \leftarrow v_k / t_k$ \;
                Update $x_{k+1} \leftarrow x_k + d_k$\;
                (Early) Terminate (or set $\delta = 0$ and proceed);
            }
            \eIf(\tcp*[f]{large value case (a)}){$|t_k| \geq \nu$}{
                $\dk \leftarrow v_k / t_k$ \label{line.exact.larga}
            }
            (\tcp*[f]{large value case (b)}){
                $\dk \leftarrow \text{sign}(-g_k^T v_k) \cdot v_k$
                \label{line.exact.largb}
            }
            Choose a stepsize $\eta_k$ by the fixed-radius strategy or the line search strategy (see Algorithm \ref{alg.backtracking})\;
            Update $x_{k+1} \leftarrow x_k + \eta_k \cdot d_k$\;
        }
    \end{algorithm}
\end{minipage}

\subsection{Preliminaries of the homogenized quadratic model}
In this subsection, we present some preliminary analysis of the homogenized quadratic model. First, we study the relationship between the smallest eigenvalues of the Hessian $\Hk$ and $F_k$, and the perturbation parameter $\delta$. Then we give the optimality conditions of problem \eqref{eq.homo subproblem} and provide some useful results based on those conditions.

\begin{lemma}[Relationship between $\lambda_1(F_k)$, $\lambda_1(H_k)$ and $\delta$]
    \label{lemma.relation of theta delta Hk}
    Let $\lambda_1(\Hk)$ and $\lambda_1(F_k)$ be the smallest eigenvalue of $H_k$ and $F_k$ respectively. Denote by $\mathcal{S}_{\lambda_1}$ the eigenspace corresponding to $\lambda_1(\Hk)$. If $g_k \neq 0$ and $H_k \neq 0$, then the following statements hold,
    \begin{enumerate}[(1)]
        \item $\lambda_1(F_k) < -\delta$ and $\lambda_1(F_k) \leq \lambda_1(\Hk)$;
        \item $\lambda_1(F_k) = \lambda_1(\Hk)$ only if $\lambda_1(\Hk)<0$ and $\gk \perp \mathcal{S}_{\lambda_1}$.
    \end{enumerate}
\end{lemma}
\begin{proof}
    We first prove the statement $(1)$. By the Cauchy interlace theorem \cite{parlett_symmetric_1998}, we immediately obtain $\lambda_1(F_k)  \leq \lambda_1(H_k) $. Now we need to prove that $\lambda_1(F_k) < -\delta $. It suffices to show that the matrix $F_k+\delta I$ has a negative eigenvalue.

    Let us consider the direction $\left[ -\eta \gk ; t \right]$, where $\eta,t>0$. Define the following function of $\left( \eta,t \right)$:
    \begin{equation*}
        \begin{aligned}
            f(\eta,t) & := \begin{bmatrix}
                               -\eta \gk \\
                               t
                           \end{bmatrix}^T (F_k+\delta I) \begin{bmatrix}
                                                              -\eta \gk \\
                                                              t
                                                          \end{bmatrix},     \\
                      & =\eta^2 \gk^T (\Hk+\delta I) \gk -2\eta t \|\gk\|^2.
        \end{aligned}
    \end{equation*}
    For any fixed $t>0$, we have
    \begin{equation*}
        f(0,t)=0 \quad  \text{and} \quad  \frac{\partial f(0,t)}{\partial \eta}=-2t \|\gk\|^2<0.
    \end{equation*}
    Therefore, for sufficiently small $\eta>0$, it holds that $f(\eta,t)<0$, which shows that $\left[ -\eta \gk ; t \right]$ is a negative curvature. Hence, $\lambda_1(F_k) < -\delta $.

    % Now we need to prove that $\lambda_1(F_k) < -\delta $. Denote $u_j = (u_j^{'}; t_j)$ with $u_j^{'} \in \mathbb{R}^{n}$ and $t_j \in \mathbb{R}$ for $j = 1, 2, \cdots, n+1$. Note that there exists at least one eigenpair of $F_k$ such that $t_j \neq 0$ and we take $ d_j = u_j^{'} / t_j$ for such $j$. By the definition of eigenvalue and eigenvector, we have 
    % \begin{equation*}
    %     F_k u_j = - \lambda_j u_j, 
    % \end{equation*}
    % implying that 
    % \begin{equation}
    %     \label{eq.equivoptc}
    %     (\Hk+\lambda_j I) d_j = -\gk \quad \text{and} \quad \gk^T d_j = \delta-\lambda_j.
    % \end{equation}

    % Let $\beta_i =  \gk^T v_i$ for all $i = 1, 2, \cdots, n$. We then have $\gk^T d_j = -\gk^T(\Hk+\lambda_j  I)^{\dagger}\gk = -\sum_{i=1}^n  \beta_i^2 / \left(  \sigma_i+\lambda_j \right)$. Therefore, $\lambda_j$ must satisfy
    % \begin{equation}
    %     \label{eq.propertytheta}
    %     \delta-\lambda_j =  -\sum_{i=1}^n\frac{\beta_i^2}{\sigma_i + \lambda_j}.
    % \end{equation}
    % Let us view both sides of the above equation as the functions of $\lambda_j$. The left-hand-side is a linear function in $\lambda_j$ while the right-hand-side is the sum of a sequence of inverse proportional functions. By the property of the two functions and the Cauchy interlace theorem \cite{parlett1998symmetric}, there exists a solution to \eqref{eq.propertytheta} such that $-\lambda_j<-\delta$, implying that $ \lambda_1 \left( F_k \right) = - \lambda_1 \leq \lambda_j <  - \delta$.  

    The proof of the statement $(2)$ is similar to the one of Theorem 3.1 in \cite{rojas_new_2001}, so we omit it here for the succinctness of the paper.
    % Then we prove the second part. We first start with $\lambda_1(\Hk)<0$, $\gk \perp \mathcal{S}_{\lambda_1}$, $\delta = \lambda_1(\Hk)-\gk^\top p_{\lambda_1}$.
    % To prove $\theta_k = -\delta_1$, it suffices to prove that under our choice of $\delta$, there is a eigenvector correspondings to $\delta_1$ with $t\neq0$,which is equivalent to the LHS and RHS intersects at a point with its first coordinate equals $-\delta_1$, this completes the proof. Hence we need to construct a eigenvector $[v;t]$ corresponding to $\delta_1$ with $t\neq0$. The construction is $[p_{\lambda_1}^T;1]$, note that by our construction of $p_{\lambda_1}$ and $\delta$, we have
    % \begin{equation}
    %     \label{eq.propertycons}
    %     (\Hk-\delta_1 I)p_{\lambda_1}= -(\Hk-\delta_1 I)(\Hk-\delta_1 I)^{\dagger}\gk=-\gk,
    % \end{equation}
    % where the second equality is due to that $(\Hk-\delta_1 I)(\Hk-\delta_1 I)^{\dagger}$ is a projection onto the image space $\mathcal{R}(\Hk-\delta_1 I)$ and that $\gk\perp \mathcal{S}_1$ thus belonging to the image space.
    % \begin{equation}
    %     \label{eq.checkeigenvec}
    %     \begin{bmatrix}
    %         \Hk   & \gk    \\
    %         \gk^T & \delta
    %     \end{bmatrix}\begin{bmatrix}
    %         p_{\lambda_1} \\1
    %     \end{bmatrix}=\begin{bmatrix}
    %         \Hk p_{\lambda_1}+\gk \\
    %         \delta+\gk^T p_{\lambda_1}
    %     \end{bmatrix}=\delta_1\begin{bmatrix}
    %         p_{\lambda_1} \\
    %         1
    %     \end{bmatrix}.
    % \end{equation}
    % The opposite part is similar so we omit it here. 
\end{proof}
\cref{lemma.relation of theta delta Hk} shows that we can control the smallest eigenvalue of the homogenized matrix $F_k$ by adjusting the perturbation parameter $\delta$. It helps us find a better direction to decrease the value of the objective function. We also note that the case $\gk\perp \mathcal{S}_{\lambda_1}$ is often regarded as a hard case in solving the trust-region subproblem. However, this challenge will not incapacitate HSODM in our convergence analysis. In the following, we will show the function value has a sufficient decrease under this scenario. Thus, the subproblem in HSODM is much easier to solve than the trust-region subproblem due to the non-existence of the hard case.

We remark that \cref{lemma.relation of theta delta Hk} is a simpler version of Lemma 3.3 in \cite{rojas_new_2001},
where the authors give a more detailed analysis of the relationship between the perturbation parameter $\delta$ and the eigenpair of the homogenized matrix $F_k$.
However, the difference is that they try to obtain a solution to the trust-region subproblem via the homogenization trick, while our goal is to seek a good direction to decrease the function value. Furthermore, if the homogenized model is used, then we can show that HSODM has the optimal $O(\epsilon^{-3/2})$ iteration complexity. However, if the homogenization trick is put on solving the trust-region subproblem as in \cite{rojas_new_2001}, one still needs a framework like the one in \citet{curtis_trust_2017} to guarantee the same convergence property. Moreover, a sequence of homogenized problems needs to be solved in each iteration of the framework.

In the following lemma, we characterize the optimal solution $[v_k; t_k]$ of problem \eqref{eq.homo subproblem} based on the optimality condition of the standard trust-region subproblem.

\begin{lemma}[Optimality condition]\label{lemma.optimal condition of subproblem}
    $[v_k; t_k]$ is the optimal solution of the subproblem \eqref{eq.homo subproblem} if and only if there exists a dual variable $\theta_k > \delta \geq 0$ such that
    \begin{align}
        \label{eq.homoeig.soc}
                                    & \begin{bmatrix}
                                          \Hk + \theta_k \cdot I & \gk              \\
                                          \gk^T                  & -\delta+\theta_k
                                      \end{bmatrix} \succeq 0, \\
        \label{eq.homoeig.foc}
                                    & \begin{bmatrix}
                                          \Hk + \theta_k \cdot I & \gk              \\
                                          \gk^T                  & -\delta+\theta_k
                                      \end{bmatrix}
        \begin{bmatrix}
            v_k \\ t_k
        \end{bmatrix} = 0,                                                      \\
        \label{eq.homoeig.norm one} & \|[v_k; t_k]\| = 1.
    \end{align}
    Moreover, $-\theta_k$ is the smallest eigenvalue of the perturbed homogenized matrix $F_k$, i.e., $-\theta_k = \lambda_1(F_k)$.
\end{lemma}
\begin{proof}
    By the optimality condition of the standard trust-region subproblem, $[v_k; t_k]$ is the optimal solution if and only if there exists a dual variable $\theta_k \geq 0$ such that
    \begin{equation*}
        \begin{bmatrix}
            \Hk + \theta_k \cdot I & \gk              \\
            \gk^T                  & -\delta+\theta_k
        \end{bmatrix} \succeq 0, \ \begin{bmatrix}
            \Hk + \theta_k \cdot I & \gk              \\
            \gk^T                  & -\delta+\theta_k
        \end{bmatrix} \begin{bmatrix}
            v_k \\ t_k
        \end{bmatrix} = 0, \ \text{and} \  \theta_k \cdot (\|[v_k; t_k]\| - 1) = 0.
    \end{equation*}
    With \cref{lemma.relation of theta delta Hk}, we have $\lambda_1(F_k) < -\delta \leq 0$.
    Therefore, $\theta_k \geq  - \lambda_1(F_k) >  \delta \geq 0$, and further $ \|[v_k; t_k]\| = 1 $. Moreover, by \eqref{eq.homoeig.foc}, we obtain
    \begin{equation*}
        \begin{bmatrix}
            \Hk   & \gk     \\
            \gk^T & -\delta
        \end{bmatrix} \begin{bmatrix}
            v_k \\ t_k
        \end{bmatrix} = -\theta_k \begin{bmatrix}
            v_k \\ t_k
        \end{bmatrix}.
    \end{equation*}
    Multiplying the equation above by $ \left[ v_k ; t_k \right]^T $, we have
    \begin{equation*}
        \min_{\|[v; t]\| \le 1} \psi_k(v, t; \delta) = -\theta_k
    \end{equation*}
    Note that with \eqref{eq.homoeig.norm one}, the optimal value of problem \eqref{eq.homo subproblem} is equivalent to the smallest eigenvalue of $F_k$, i.e., $\lambda_1(F_k)$. Thus, $-\theta_k = \lambda_1(F_k)$. The proof is then completed.
\end{proof}

With the above optimality condition, we can derive the following corollaries.
\begin{corollary}
    \label{corollary. foc}
    The equation \eqref{eq.homoeig.foc} in \cref{lemma.optimal condition of subproblem} can be rewritten as,
    \begin{equation}\label{eq.homoeig.foc cont}
        (\Hk + \theta_k I)  v_k =  - t_k g_k \quad \text{and} \quad \gk^T v_k  =  t_k (\delta - \theta_k).
    \end{equation}
    Furthermore,
    \begin{enumerate}
        \item[(1)] If $t_k = 0$, then we have
            \begin{equation}\label{eq.homoeig.foc t=0}
                (\Hk + \theta_k I)  v_k = 0 \quad \text{and} \quad \gk^T v_k  = 0,
            \end{equation}
            implying that $(-\theta_k, v_k)$ is the eigenpair of the Hessian matrix $H_k$.
        \item[(2)] If $t_k \neq 0$, then we have
            \begin{equation}
                \label{eq.homoeig.foc t neq 0}
                \gk^T d_k = \delta -\theta_k \quad \text{and} \quad (\Hk+\theta_k \cdot I)d_k =-\gk
            \end{equation}
            where $d_k =v_k / t_k$.
    \end{enumerate}
\end{corollary}
The corollary above is a direct application of \cref{lemma.optimal condition of subproblem}, so we omit its proof in the paper.

\begin{corollary}[Nontriviality of direction $v_k$] \label{coro.homo.exists.nc}
    If $g_k \neq 0$, then $v_k \neq 0$.
\end{corollary}
\begin{proof}
    We prove this by contradiction. Suppose that $v_k = 0$. Then, we have $t_k g_k = 0$ with equation \eqref{eq.homoeig.foc cont} in  \cref{corollary. foc}. It further implies that $t_k = 0$ due to $g_k \neq 0$. However, $[ v_k ; t_k]=0$ contradicts to the equation $\|[v_k; t_k]\| = 1$ in the optimality condition. Therefore, we have $v_k \neq 0$.
\end{proof}

This corollary shows that a nontrivial direction $v_k$ always exists, thus \cref{alg.main alg} will not get stuck.

\begin{corollary}\label{corollary. sign property}
    For the sign function value $\text{sign}(-g_k^T v_k)$, we always have $\text{sign}(-g_k^T v_k) \cdot t_k = |t_k|$.
\end{corollary}
\begin{proof}
    By the second equation of optimal condition \eqref{eq.homoeig.foc cont}, and $\delta < \theta_k$, we obtain that
    $$
        \text{sign}(-g_k^T v_k) = \text{sign}(t_k),
    $$
    and it implies
    $$
        \text{sign}(-g_k^T v_k) \cdot t_k = \text{sign}(t_k) \cdot t_k = |t_k|.
    $$
    This completes the proof.
\end{proof}

As a byproduct, we also have the following result.

\begin{corollary}[Trivial case, $g_k = 0$]
    Suppose that $g_k = 0$, then the following statements hold, \begin{enumerate}
        \item[(1)] If $\lambda_1(H_k) > -\delta$, then $t_k = 1$.
        \item[(2)] If $\lambda_1(H_k) < -\delta$, then $t_k = 0$.
    \end{enumerate}
\end{corollary}
\begin{proof}
    When $g_k = 0$, the homogenized matrix $F_k = [H_k, 0; 0, -\delta]$, and the subproblem \eqref{eq.homo subproblem} is
    \begin{equation*}
        \min_{\|[v; t]\| \le 1} \psi_k(v, t; \delta) = v^T H_k v - t^2 \cdot \delta.
    \end{equation*}
    We first prove the statement $(1)$ by contradiction. Suppose that $t_k \neq 1$, then we have $v_k \neq 0$ by the equation \eqref{eq.homoeig.norm one}. Thus,
    \begin{equation}
        \label{eq. gk=0 contra}
        \psi_k(v_k, t_k; \delta) = (v_k)^T H_k v_k - t_k^2 \cdot \delta > -\delta = \psi_k(0, 1; \delta),
    \end{equation}
    where the inequality holds due to $(v_k)^T H_k v_k \geq \lambda_1(H_k) \|v_k\|^2 > -\delta\|v_k\|^2$. The equation \eqref{eq. gk=0 contra} contradicts to the optimality of $(v_k, t_k)$, and thus $t_k = 1$. The second statement can be proved by the same argument, and we omit the proof here.
\end{proof}

\section{Global Convergence Rate}\label{sec.global}
In this section, we analyze the convergence rate of the proposed HSODM. To facilitate the analysis, we present two building blocks considering the large and small values of $\|d_k\|$, respectively.
For the large value case of $\|d_k\|$, we show that the function value decreases by at least $\Omega(\epsilon^{3/2})$ at every iteration after carefully selecting the perturbation parameter $\delta$. In the latter case, we prove that the next iterate $x_{k+1}$ is already an $\epsilon$-approximate SOSP, and thus the algorithm can terminate. Throughout the paper, we make the following standard assumptions.

\begin{assumption}
    \label{assm.lipschitz}
    Assume that $f$ has \(M\)-Lipschitz continuous Hessian on an open convex set X containing all the iterates $\xk$, i.e., for some $M > 0$, we have
    \begin{equation}\label{eq.assm.lipschitz}
        \|\nabla^2 f(x) - \nabla^2 f(y)\| \le M \|x-y\|, ~\forall x, y \in X,
    \end{equation}
    and that the Hessian matrix is bounded,
    \begin{equation}\label{eq.boundhess}
        \|\nabla^2 f(\xk) \| \le U_H, ~\forall k \ge 0,
    \end{equation}
    for some $U_H > 0$.
\end{assumption}

We also recall the next lemma for preparation.

\begin{lemma}[\citet{nesterovLecturesConvexOptimization2018}]\label{lem.lipschitz}
    If \(f:\mathbb{R}^n \mapsto \mathbb{R}\) satisfies \cref{assm.lipschitz}, then for all \(x,y\in \mathbb{R}^n\),
    \begin{subequations}
        \begin{align}
            % & |f(y)-f(x)-\nabla f(x)^T(y-x)|                                             \leq \frac{L}{2}\|y-x\|^{2}        \\
             & \left\|\nabla f(y)-\nabla f(x)-\nabla^{2} f(x)(y-x)\right\|                       \leq \frac{M}{2}\|y-x\|^{2}, \\
             & \left|f(y)-f(x)-\nabla f(x)^T(y-x)-\frac{1}{2}(y-x)^T\nabla^{2} f(x)(y-x)\right|  \leq \frac{M}{6}\|y-x\|^{3}.
        \end{align}
    \end{subequations}
\end{lemma}

\subsection{Analysis for the large value of $\|d_k\|$ }\label{subsection.decrease}

In HSODM, we define the large-value case of $\|d_k\|$ as the case that its norm is larger than the trust-region radius $\Delta$, i.e., $\|d_k\| > \Delta$. Note that in the case of $\nu \le |t_k| \le \sqrt{1/(1+\Delta^2)}$, we have $\|d_k\| = \|v_k\|/|t_k| = \sqrt{1-|t_k|^2}/|t_k| \ge \Delta$. Moreover, in the case of $|t_k| \le \nu$ with $\nu \in (0, 1/2)$, it holds that $\|d_k\| = \|v_k\| = \sqrt{1-|t_k|^2} \ge \sqrt{3}/2 \ge \Delta= \Theta(\sqrt{\epsilon})$. Therefore, we call these two cases the large value case (a) and (b) in \cref{alg.main alg}, respectively. In this situation, the homogenized direction can be either $d_k = \text{sign}(-g_k^T v_k) \cdot v_k$ or $d_k = v_k / t_k$. The following discussion shows that both stepsize selection strategies result in a sufficient decrease. The analysis for the fixed-radius strategy is more concise and clear, but it mainly serves as a theoretical result. On the contrary, the line search stepsize selection strategy is more practical in spite of a slightly more complicated analysis.

\subsubsection{Fixed-radius strategy}

For the fixed-radius strategy, the next iterate $x_{k+1}$ is constrained to satisfy $\|x_{k+1} - x_k\| = \Delta$, and hence the stepsize is selected as $\Delta/\|d_k\|$. Firstly, we will consider the scenario in which $|t_k| < \nu$ and $d_k = \text{sign}(-g_k^T v_k) \cdot v_k$.  We remark that this particular scenario encompasses the so-called ``hard case'' ($t_k = 0$) in trust-region methods \cite{rojas_new_2001}. When $t_k = 0$, \cref{corollary. foc} shows that $(-\theta_k, v_k)$ is an eigenpair of the Hessian $H_k$, and $v_k$ is a sufficiently negative curvature direction due to $-\theta_k < - \delta \leq 0$. Therefore, moving along the direction of $v_k$ with an appropriate stepsize will always decrease the function value~\cite{carmon_accelerated_2018}. We first present a lemma that applies to the case $|t_k| < \nu$, and it can be regarded as a generalized descent lemma.

\begin{lemma}\label{lemma.t < delta fix radius decrease lemma}
    Suppose that \cref{assm.lipschitz} holds and set $\nu \in (0, 1/2)$. If $|t_k| < \nu$, then let $d_k = \text{sign}(-g_k^T v_k) \cdot v_k$ and $\eta_k = \Delta/\|d_k\|$, we have
    \begin{equation}
        f(x_{k+1}) - f(x_k) \leq -\frac{\Delta^2}{2}\delta + \frac{M}{6}\Delta^3.
    \end{equation}
\end{lemma}
\begin{proof}
    When $d_k = \text{sign}(-g_k^T v_k) \cdot v_k$, with the optimality condition \eqref{eq.homoeig.foc cont} in \cref{corollary. foc} and \cref{corollary. sign property}, we obtain
    \begin{equation}
        \label{eq:bridge2}
        d_k^T H_k d_k = - \theta_k \|d_k\|^2 - t_k^2 \cdot (\delta - \theta_k) \quad \text{and} \quad
        g_k^T d_k = |t_k| \cdot(\delta - \theta_k).
    \end{equation}
    Since $\eta_k = \Delta / \|d_k\| \in (0, 1)$, then $  \eta_k - \eta_k^2 / 2 \geq 0 $, and further
    \begin{equation}
        \label{eq:bridge6}
        \left(\eta_k - \frac{\eta_k^2}{2}\right) \cdot (\delta - \theta_k) \leq 0.
    \end{equation}
    By the $M$-Lipschitz continuous property of $\nabla^2 f(x)$, we have
    \begin{subequations}
        \begin{align}
            f(x_{k+1}) - f(x_k) & = f(x_k + \eta_k d_k) - f(x_k)  \notag                                                                                                                                                                                            \\
                                & \leq \eta_k \cdot g_k^T d_k + \frac{\eta_k^2}{2} \cdot d_k^T H_k d_k + \frac{M}{6}\eta_k^3 \|d_k\|^3  \notag                                                                                                                      \\
                                & = \eta_k \cdot |t_k| \cdot(\delta - \theta_k) - \frac{\eta_k^2}{2} \cdot \theta_k\|d_k\|^2 - \frac{\eta_k^2}{2} \cdot t_k^2 \cdot (\delta - \theta_k) + \frac{M}{6}\eta_k^3 \|d_k\|^3 \label{subeq. substitute optimal condition} \\
                                & \leq \eta_k \cdot t_k^2 \cdot(\delta - \theta_k) - \frac{\eta_k^2}{2} \cdot \theta_k\|d_k\|^2 - \frac{\eta_k^2}{2} \cdot t_k^2 \cdot (\delta - \theta_k) + \frac{M}{6}\eta_k^3 \|d_k\|^3  \label{subeq. relax |tk|}               \\
                                & = \left(\eta_k - \frac{\eta_k^2}{2}\right) \cdot t_k^2 \cdot (\delta - \theta_k) - \frac{\eta_k^2}{2} \cdot \theta_k\|d_k\|^2 + \frac{M}{6}\eta_k^3 \|d_k\|^3  \notag                                                             \\
                                & \leq - \theta_k \cdot \frac{\Delta^2}{2} + \frac{M}{6}\Delta^3 \label{subeq. t=0 optimal condition}                                                                                                                               \\
                                & \leq -\frac{\Delta^2}{2}\delta + \frac{M}{6}\Delta^3, \label{subeq. t=0 delta}
        \end{align}
    \end{subequations}
    where \eqref{subeq. substitute optimal condition} follows from \eqref{eq:bridge2}, and \eqref{subeq. relax |tk|} holds due to $|t_k| < \nu < 1$ and $\delta - \theta_k < 0$. The inequality \eqref{subeq. t=0 optimal condition} holds by \eqref{eq:bridge6} and $\eta_k = \Delta/\|d_k\|$.
\end{proof}

Now we turn to the case $|t_k| \geq \nu$, and let the update direction $d_k = v_k / t_k$. When $\|d_k\|$ is large enough, i.e., $\|d_k\| > \Delta$, we can obtain the same decrease of function value in the next lemma.

\begin{lemma}\label{lemma.t > delta fix radius decrease lemma}
    Suppose that \cref{assm.lipschitz} holds and set $\nu \in (0, 1/2)$. If $|t_k| \geq \nu$ and $\|v_k / t_k\| > \Delta$, then let $d_k = v_k / t_k$ and $\eta_k = \Delta / \|d_k\|$, we have
    \begin{equation}
        f(x_{k+1}) - f(x_k) \leq -\frac{\Delta^2}{2}\delta + \frac{M}{6}\Delta^3.
    \end{equation}
\end{lemma}
\begin{proof}
    When $t_k \neq 0$, with equation \eqref{eq.homoeig.foc t neq 0} in \cref{corollary. foc}, we have
    \begin{equation}
        \label{eq:bridge3}
        d_k^T H_k d_k = -g_k^T d_k - \theta_k \|d_k\|^2 \quad \text{and} \quad g_k^T d_k = \delta - \theta_k \leq 0.
    \end{equation}
    Since $\eta_k = \Delta / \|d_k\| \in (0, 1)$, then $  \eta_k - \eta_k^2 / 2 \geq 0 $, and further
    \begin{equation}
        \label{eq:bridge4}
        \left(\eta_k - \frac{\eta_k^2}{2}\right) \cdot g_k^T d_k \leq 0.
    \end{equation}

    By the $M$-Lipschitz continuous property of $\nabla^2 f(x)$, we have
    \begin{subequations}
        \begin{align}
            f(x_{k+1}) - f(x_k) & = f(x_k + \eta_k d_k) - f(x_k)  \notag                                                                                                                                             \\
                                & \leq \eta_k \cdot g_k^T d_k + \frac{\eta_k^2}{2} \cdot d_k^T H_k d_k + \frac{M}{6}\eta_k^3 \|d_k\|^3  \notag                                                                       \\
                                & = \left(\eta_k - \frac{\eta_k^2}{2}\right) \cdot g_k^T d_k - \theta_k \cdot \frac{\eta_k^2}{2} \|d_k\|^2 + \frac{M}{6}\eta_k^3 \|d_k\|^3  \label{subeq. t neq 0 optimal condition} \\
                                & \leq -\theta_k \cdot \frac{\eta_k^2}{2} \|d_k\|^2 + \frac{M}{6}\eta_k^3 \|d_k\|^3 \label{subeq. t neq 0 drop gd}                                                                   \\
                                & \leq -\frac{\Delta^2}{2}\delta + \frac{M}{6}\Delta^3, \label{subeq. t neq 0 delta}
        \end{align}
    \end{subequations}
    where \eqref{subeq. t neq 0 optimal condition} holds due to equation \eqref{eq:bridge3}, \eqref{subeq. t neq 0 drop gd} follows from equation \eqref{eq:bridge4}, and in \eqref{subeq. t neq 0 delta} we substitute $\eta_k$ with $\Delta / \| d_k\| $ and use $\theta_k \geq \delta$.
\end{proof}

\subsubsection{Line search strategy}
For the line search strategy, we utilize a backtracking subroutine to determine the stepsize $\eta_k$, ensuring it produces a sufficient decrease. The details of the subroutine are provided below.

\begin{minipage}[t]{0.95\linewidth}
    \begin{algorithm}[H]
        \caption{Backtracking Line Search}\label{alg.backtracking}
        \KwData{Given current iterate $x_k$,  direction $d_k$, initial stepsize $\eta_k = 1$, $\gamma > 0$, $\beta \in (0, 1)$}
        \textbf{For} $j = 0, 1, 2, \cdots$ \textbf{do:} \\
        \quad Compute decrease quantity $D_k := f(x_k) - f(x_k + \eta_k d_k)$; \\
        \quad \textbf{If} $D_k \geq \gamma \eta_k^3\|d_k\|^3 / 6 $ \textbf{then:} \\
        \quad \quad \textbf{Break;} \\
        \quad \textbf{Else:} \\
        \quad \quad Update $\eta_k := \beta \cdot \eta_k$; \\
        \textbf{Output:} stepsize $\eta_k$.
    \end{algorithm}
    \vspace{1em}
\end{minipage}

Similarly, we derive the descent lemmas with the line search strategy and further upper bound the number of iterations required by the line search procedure. For the cases $|t_k| < \nu$ and $|t_k| \geq \nu$, we obtain the following two lemmas that characterize the sufficient decrease property.

\begin{lemma}\label{lemma.t < delta line search decrease lemma}
    Suppose that \cref{assm.lipschitz} holds and set $\nu \in (0, 1/2)$. If $|t_k| < \nu$, then let $d_k = \text{sign}(-g_k^T v_k) \cdot v_k$. The backtracking line search terminates with $\eta_k = \beta^{j_k}$,  and $j_k$ is upper bounded by
    $$
        j_N := \left \lceil \log_\beta\left(\frac{3\delta}{M+\gamma}\right) \right\rceil,
    $$and the function value associated with the stepsize $\eta_k$ satisfies,
    \begin{equation}\label{eq. function value decrease line search t < delta}
        f(x_{k+1}) - f(x_k) \leq -\min \left\{\frac{\sqrt{3}\gamma}{16}, \frac{9\gamma\beta^3\delta^3}{2(M+\gamma)}\right\}.
    \end{equation}
\end{lemma}
\begin{proof}
    Suppose that the backtracking line search terminate with $\eta_k = 1$, then we have
    $$
        \begin{aligned}
            f(x_k + \eta_kd_k) - f(x_k) & \leq -\frac{\gamma}{6}\eta_k^3\|d_k\|^3  = -\frac{\gamma}{6}\|v_k\|^3            \leq -\frac{\sqrt{3}\gamma}{16},
        \end{aligned}
    $$
    where the last inequality is due to $\|v_k\| = \sqrt{1-|t_k|^2} \geq \sqrt{1-\nu^2} \geq \sqrt{3}/2$. Suppose the algorithm does not stop at the iteration $j \geq 0$ and
    the condition in Line 4 is not met, i.e., $D_k < \frac{\gamma}{6}\beta^{3j}\|d_k\|^3 = \frac{\gamma}{6}\beta^{3j}\|v_k\|^3$. By using a similar argument in the proof of \cref{lemma.t < delta fix radius decrease lemma}, we have that
    \begin{equation}\label{ineq:line-search-iter-est}
        \begin{aligned}
            -\frac{\gamma}{6}\beta^{3j}\|v_k\|^3 & < f(x_k + \beta^j d_k) - f(x_k)                                                                                                                                                               \\
                                                 & \leq \beta^j \cdot g_k^T d_k + \frac{\beta^{2j}}{2} \cdot d_k^T H_k d_k + \frac{M}{6}\beta^{3j}\|d_k\|^3                                                                                      \\
                                                 & = \beta^j \cdot |t_k| \cdot (\delta-\theta_k) - \frac{\beta^{2j}}{2} \cdot \theta_k\|v_k\|^2 - \frac{\beta^{2j}}{2} \cdot t_k^2 \cdot (\delta - \theta_k) + \frac{M}{6}\beta^{3j}\|v_k\|^3    \\
                                                 & \leq \beta^j \cdot t_k^2 \cdot (\delta-\theta_k) - \frac{\beta^{2j}}{2} \cdot \theta_k\|v_k\|^2 - \frac{\beta^{2j}}{2} \cdot t_k^2 \cdot (\delta - \theta_k) + \frac{M}{6}\beta^{3j}\|v_k\|^3 \\
                                                 & = \left(\beta^j - \frac{\beta^{2j}}{2}\right) \cdot t_k^2 \cdot (\delta-\theta_k) -  \frac{\beta^{2j}}{2} \cdot \theta_k\|v_k\|^2 + \frac{M}{6}\beta^{3j}\|v_k\|^3                            \\
                                                 & \leq - \frac{\beta^{2j}}{2} \cdot \theta_k\|v_k\|^2 + \frac{M}{6}\beta^{3j}\|v_k\|^3                                                                                                          \\
                                                 & \leq - \frac{\beta^{2j}}{2} \cdot \delta\|v_k\|^2 + \frac{M}{6}\beta^{3j}\|v_k\|^3.
        \end{aligned}
    \end{equation}
    Therefore, $\beta^j > \frac{3\delta}{(M+\gamma)\|v_k\|}$ holds, which further implies that
    $$
        j < \log_\beta\left(\frac{3\delta}{(M+\gamma)\|v_k\|}\right).
    $$
    However, $j_N := \left \lceil \log_\beta\left(\frac{3\delta}{M+\gamma}\right) \right\rceil \geq \log_\beta\left(\frac{3\delta}{(M+\gamma)\|v_k\|}\right)$ due to $\|v_k\| \leq 1$.
    This means that the inequality \eqref{ineq:line-search-iter-est} does not hold when $j = {j_N}$, and thus
    the condition in Line 4 is satisfied in this case. Therefore, the iteration number of backtracking subroutine $j_k$ is upper bounded by $j_N$, and the function value decreases as
    $$
        \begin{aligned}
            f(x_k + \eta_k d_k) - f(x_k)
             & \leq -\frac{\gamma}{6}\beta^{3j_k}\|v_k\|^3             \\
             & = -\frac{\gamma \beta^3}{6} \beta^{3(j_k - 1)}\|v_k\|^3 \\
             & \leq - \frac{9\gamma\beta^3\delta^3}{2(M+\gamma)^3},
        \end{aligned}
    $$
    where the last inequality comes from $\beta^{j_k - 1} \geq \frac{3\delta}{(M+\gamma)\|v_k\|}$.
\end{proof}

% For the case of $|t_k| \geq \nu$, we have the following estimation on the iteration number of the line search and the value decrease on the objective function.

\begin{lemma}\label{lemma.t > delta line search decrease lemma}
    Suppose that \cref{assm.lipschitz} holds and set $\nu \in (0, 1/2)$. If $|t_k| \ge \nu$ and $\|v_k / t_k\| > \Delta$, then let $d_k = v_k/t_k$. The backtracking line search terminates with $\eta_k = \beta^{j_k}$, and $j_k$ is upper bounded by
    $$
        j_N := \left \lceil \log_\beta\left(\frac{3\delta\nu}{M+\gamma}\right) \right\rceil,
    $$and the function value associated with the stepsize $\eta_k$ satisfies,
    \begin{equation}\label{eq. function value decrease line search t > delta}
        f(x_{k+1}) - f(x_k) \leq -\min \left\{\frac{\gamma\Delta^3}{6}, \frac{9\gamma\beta^3\delta^3}{2(M+\gamma)^3}\right\}.
    \end{equation}
\end{lemma}
\begin{proof}
    Similarly, suppose that the backtracking line search terminates with $\eta_k = 1$, we have
    $$
        \begin{aligned}
            f(x_k + \eta_kd_k) - f(x_k) & \leq -\frac{\gamma}{6}\eta_k^3\|d_k\|^3 \\
                                        & \leq -\frac{\gamma}{6}\Delta^3,
        \end{aligned}
    $$
    where the last inequality comes from $\|d_k\| > \Delta$. If $\eta_k = 1$ does not lead to a sufficient decrease, then for any $j \geq 0$ where the condition in Line 4 is not met, we have
    \begin{equation}\label{ineq:line-search-iter-est2}
        \begin{aligned}
            -\frac{\gamma}{6}\beta^{3j}\|d_k\|^3 & < f(x_k + \beta^j d_k) - f(x_k)                                                                                                        \\
                                                 & \leq \beta^j \cdot g_k^T d_k + \frac{\beta^{2j}}{2} \cdot d_k^T H_k d_k + \frac{M}{6}\beta^{3j}\|d_k\|^3                               \\
                                                 & = (\beta^j - \frac{\beta^{2j}}{2}) \cdot (\delta_k - \theta_k) -\frac{\beta^{2j}}{2}\theta_k\|d_k\|^2 + \frac{M}{6}\beta^{3j}\|d_k\|^3 \\
                                                 & \leq -\frac{\beta^{2j}}{2}\delta\|d_k\|^2 + \frac{M}{6}\beta^{3j}\|d_k\|^3.
        \end{aligned}
    \end{equation}
    Therefore, $\beta^j \geq \frac{3\delta}{(M+\gamma)\|d_k\|}$ and it implies that
    $$
        j < \log_\beta\left(\frac{3\delta}{(M+\gamma)\|d_k\|}\right).
    $$
    Note that
    \begin{equation}
        \|d_k\| = \|v_k\|/|t_k| = \frac{\sqrt{1-|t_k|^2}}{|t_k|} \leq \frac{1}{\nu},
    \end{equation}
    and
    $j_N := \left \lceil \log_\beta\left(\frac{3\delta\nu}{M+\gamma}\right) \right\rceil \geq \log_\beta\left(\frac{3\delta}{(M+\gamma)\|d_k\|}\right)$,
    This means that the inequality \eqref{ineq:line-search-iter-est2} does not hold when $j = {j_N}$, and thus
    the condition in Line 4 is satisfied in this case. Therefore, the iteration number of backtracking subroutine $j_k$ is upper bounded by $j_N$,
    and the function value decreases as
    $$
        \begin{aligned}
            f(x_k + \eta_k d_k) - f(x_k) & \leq -\frac{\gamma}{6}\beta^{3j_k}\|d_k\|^3             \\
                                         & = -\frac{\gamma \beta^3}{6} \beta^{3(j_k - 1)}\|d_k\|^3 \\
                                         & \leq - \frac{9\gamma\beta^3\delta^3}{2(M+\gamma)^3},
        \end{aligned}
    $$
    where the last inequality is due to  $\beta^{j_k - 1} \geq \frac{3\delta}{(M+\gamma)\|d_k\|}$.
\end{proof}

Combining the above two lemmas, we now conclude a unified descent property for homogenized negative curvature equipped with a backtracking line search.

\begin{corollary}
    \label{corrollary. descend for line search}
    Suppose that \cref{assm.lipschitz} holds and set $\nu \in (0, 1/2)$. Let the backtracking line search parameters $\beta, \gamma$ satisfy $\beta \in (0, 1)$ and $\gamma > 0$. Then, after every outer iterate, the function value decreases as
    $$
        f(x_{k+1}) - f(x_k) \leq -\min\left\{\frac{\sqrt{3}\gamma}{16}, \frac{9\gamma\beta^3\delta^3}{2(M+\gamma)}, \frac{\gamma\Delta^3}{6}, \frac{9\gamma\beta^3\delta^3}{2(M+\gamma)^3}\right\}.
    $$
    and the inner iteration for backtracking line search is at most
    $$
        j_N \leq \max\left\{\left\lceil \log_\beta\left(\frac{3\delta}{M+\gamma}\right) \right\rceil, \left \lceil \log_\beta\left(\frac{3\delta\nu}{M+\gamma}\right) \right\rceil\right\} = \left \lceil \log_\beta\left(\frac{3\delta\nu}{M+\gamma}\right) \right\rceil.
    $$
\end{corollary}
\begin{remark}
    An interesting implication of \cref{corrollary. descend for line search} is that the amount of
    value decrease of the objective function
    is almost unaffected by the choice of $\nu$, the truncation parameter. The choice of $\nu$ only affects the number of iterations for backtracking line search, which is $O(\log_\beta(\delta\nu))$. Nevertheless, it is not suggested to choose small $\nu$,
    which will increase
    the complexity of line search as $\beta < 1$.
\end{remark}

\subsection{Analysis for the small value of $\|d_k\|$}
\label{subsection.convergence}

In this subsection, we consider the small value case where $\|d_k\| \le \Delta$. Note that in the case of $|t_k| \ge \sqrt{1/(1+\Delta)}$, we have $\|d_k\| = \|v_k\|/|t_k| = \sqrt{1-|t_k|^2}/|t_k| \le \Delta$, validating the name of small value case in \cref{alg.main alg}. Under this case, we prove that the next iterate $x_{k+1} = x_k + d_k$ is already an $\epsilon$-approximate SOSP. Therefore, we can terminate the algorithm after one iteration in the small value case. To prove this result, we provide an upper bound of $\|g_k\|$ for preparation.
\begin{lemma}
    \label{lemma.upper bound of g_k}
    Suppose that \cref{assm.lipschitz} holds. If $g_k \neq 0$, and $\|d_k\| \leq \Delta \leq \sqrt{2}/2$, then we have
    \begin{equation}\label{eq.upper bound of g_k}
        \|g_k\| \leq 2(U_H+\delta) \Delta.
    \end{equation}
\end{lemma}
\begin{proof}
    By \cref{lemma.relation of theta delta Hk}, we have $\theta_k - \delta > 0 $. Moreover, with equation \eqref{eq.homoeig.foc t neq 0} in \cref{corollary. foc}, we can give an upper bound of $\theta_k - \delta$, that is,
    \begin{equation}\label{eq.bound of theta - delta}
        \theta_k - \delta = -g_k^T d_k \leq \|g_k\| \|d_k\| \leq \Delta \|g_k\|.
    \end{equation}
    Define $h(t) = t^2 + \left(  g_k^T H_k g_k / \|g_k\|^2 + \delta  \right) t - \|g_k\|^2$. It is easy to see that the equation $h(t) = 0$ must have two real roots with opposite signs. Let its positive root be $t_2$. By $ \theta_k - \delta > 0$, we have $ \theta_k - \delta \geq t_2 $. Therefore, we must have
    \begin{equation*}
        h(\Delta \|g_k\|) =
        \ \Delta^2\|\gk\|^2+\left(\frac{g_k^T H_k g_k}{\|g_k\|^2} + \delta \right)\Delta\|\gk\|-\|\gk\|^2\geq 0.
    \end{equation*}
    After some algebra, we obtain
    \begin{align}
        \nonumber\|\gk\|   & \leq \frac{\left( g_k^T H_k g_k / \|g_k\|^2 + \delta \right) \Delta}{1-\Delta^2} \\
        \nonumber          & \leq \frac{(U_H+\delta)\Delta}{1-\Delta^2}                                       \\
        \label{eq.bound g} & \leq 2(U_H+\delta) \Delta.
    \end{align}
    The second inequality holds due to $H_k \preceq U_H I$, which implies $g_k^T H_k g_k / \|g_k\|^2 \leq U_H$. The last inequality follows from $\Delta \leq \sqrt{2}/2$.
\end{proof}

% When $\|d_k\| \leq \Delta$, we let the stepsize $\eta_k = 1$ and proceed the iteration by $x_{k+1} = x_k + d_k$. 
The following lemma shows that the norm of the gradient at $x_{k+1}$ has an upper bound, and the smallest eigenvalue of the Hessian at $x_{k+1}$ has a lower bound.

\begin{lemma}\label{lemma.norm of gk+1}
    Suppose that \cref{assm.lipschitz} holds. If $g_k \neq 0$, and $\| d_k \| \leq \Delta$, then let $\eta_k = 1$, we have
    \begin{align}
        \|g_{k+1}\| & \leq  2(U_H+\delta) \Delta^3 + \frac{M}{2} \Delta^2 + \delta \Delta, \label{eq.small norm of gk+1} \\
        H_{k+1}     & \succeq - \left(  2(U_H+\delta)\Delta^2 + M\Delta + \delta \right) I. \label{eq.psd of Hk+1}
    \end{align}
\end{lemma}
\begin{proof}
    We first prove \eqref{eq.small norm of gk+1}. By the optimality condition \eqref{eq.homoeig.foc t neq 0} in \cref{corollary. foc}, we have
    \begin{equation*}
        H_k d_k + g_k = -\theta_k d_k,
    \end{equation*}
    and with \eqref{eq.bound of theta - delta}, we have
    \begin{equation*}
        \theta_k \|d_k\| \leq \left( \delta + \Delta \|g_k\| \right) \|d_k\|.
    \end{equation*}
    Thus, it holds that
    \begin{equation}
        \label{eq:bridge1}
        \begin{split}
            \|H_k d_k + g_k \| = \theta_k \| d_k\| \leq \delta \Delta +  \|g_k\| \Delta^2.
        \end{split}
    \end{equation}
    Now we bound the norm of $\|g_{k+1}\|$ and obtain,
    \begin{subequations}
        \begin{align}
            \|g_{k+1}\| & \leq \|g_{k+1} - H_k d_k - g_k\| + \|H_k d_k + g_k\| \notag                                                              \\
                        & \leq \frac{M}{2}\|d_k\|^2 + \delta \Delta +  \|g_k\| \Delta ^2 \label{subeq. hessian lip}                                \\
                        & \leq \frac{M}{2} \Delta^2 + \delta \Delta + 2(U_H+\delta)\Delta \cdot \Delta^2 \label{subeq. convergence gk upper bound} \\
                        & =  2(U_H+\delta) \Delta^3 + \frac{M}{2} \Delta^2 + \delta \Delta, \notag
        \end{align}
    \end{subequations}
    where \eqref{subeq. hessian lip} holds due to the $M$-Lipschitz continuity of $\nabla^2 f(x)$ as well as equation \eqref{eq:bridge1}, and \eqref{subeq. convergence gk upper bound} follows from \cref{lemma.upper bound of g_k}.
    Now we prove \eqref{eq.psd of Hk+1}. Note that the optimality condition \eqref{eq.homoeig.soc} in \cref{lemma.optimal condition of subproblem} implies that
    \begin{equation*}
        H_k + \theta_k \cdot I \succeq 0.
    \end{equation*}
    With \eqref{eq.bound of theta - delta} and \eqref{eq.bound g}, we further obtain
    \begin{align}
        \nonumber   H_k      & \succeq - \theta_k I \succeq -(\Delta\|g_k\|+\delta)I \\
        \label{eq.psd pf Hk} & \succeq -2(U_H+\delta)\Delta^2I - \delta I.
    \end{align}
    To bound $H_{k+1}$, we have
    \begin{align}
        H_{k+1} & \succeq H_k - \|H_{k+1} - H_{k}\| I  \succeq H_k - M \|d_k\| I   \succeq H_k - M \Delta I,
    \end{align}
    where the second inequality holds by the $M$-Lipschitz continuity of $\nabla^2 f(x)$, and the last inequality follows from $\|d_k\| \leq \Delta$. Combining with \eqref{eq.psd pf Hk}, we arrive at
    \begin{equation}
        H_{k+1} \succeq -2(U_H+\delta)\Delta^2I - \delta I - M\Delta I.
    \end{equation}
    The proof is then complete.
\end{proof}

\subsection{The global convergence}
Putting the above pieces together, we present the formal global convergence results of HSODM in both the fixed-radius and line search strategies in \cref{thm.fix radius global convergence rate} and \cref{thm.line search global convergence rate}, respectively. It shows that our HSODM achieves $O(\epsilon^{-3/2})$ iteration complexity to find an $\epsilon$-approximate SOSP by properly choosing the perturbation parameter $\delta$ and the radius $\Delta$.
\begin{theorem}
    \label{thm.fix radius global convergence rate}
    Suppose that \cref{assm.lipschitz} holds.  Let $\delta = \sqrt{\epsilon}$, $\Delta = 2\sqrt{\epsilon} / M$ and $\nu \in (0, 1/2)$, then the homogeneous second-order descent method (HSODM) with the fixed-radius strategy terminates in at most $O\left(\epsilon^{-3/2}\right)$ steps, and the next iterate $x_{k+1}$ is a SOSP.
\end{theorem}
\begin{proof}
    Since we take $\delta = \sqrt{\epsilon}$ and $\Delta = 2\sqrt{\epsilon} / M $, by \cref{lemma.t < delta fix radius decrease lemma} and \cref{lemma.t > delta fix radius decrease lemma},  we immediately obtain that the function value decreases at least $\Omega(\epsilon^{3/2})$ for the large step case, i.e.,
    \begin{equation*}
        f(x_{k+1}) - f(x_k) \leq -\frac{2}{3M^2} \epsilon^{3/2}.
    \end{equation*}
    When the algorithm terminates, by \cref{lemma.norm of gk+1}, we have
    \begin{align}
        \nonumber   \|g_{k+1}\| & \leq 2(U_H+\delta) \Delta^3 +
        \frac{M}{2} \Delta^2 + \delta \Delta                                                                               \\
        \label{eq.gk+1.conv}    & \le \frac{16U_H \epsilon^{3/2} + 16\epsilon^2}{M^3} + \frac{4\epsilon}{M}\le O(\epsilon)
    \end{align}
    and
    \begin{align}
        \nonumber \lambda_1 (H_{k+1}) & \ge - \left(2(U_H+\delta)\Delta^2 + M\Delta + \delta\right)                                                     \\
        \label{eq.hk+1.conv}          & \geq - \left(\frac{8U_H\epsilon + 8\epsilon^{3/2}}{M^2} + 3\sqrt{\epsilon}\right) \ge \Omega(-\sqrt{\epsilon}).
    \end{align}
    Therefore, the next iterate $x_{k+1}$ is already a SOSP.
    Note that the total decreasing amount of the objective function value cannot exceed $f(x_1) - f_{\inf}$. Hence, the number of iterations for large step cases is upper bounded by
    \begin{equation*}
        O\left(\frac{3M^2}{2}\left(f(x_1) - f_{\inf}\right)\epsilon^{-3/2}\right),
    \end{equation*}
    which is also the iteration complexity of our algorithm.
\end{proof}

\begin{theorem}
    \label{thm.line search global convergence rate}
    Suppose that \cref{assm.lipschitz} holds. Let $\delta = \sqrt{\epsilon}$, $\Delta = 2\sqrt{\epsilon} / M$ and $\nu \in (0, 1/2)$, and the backtracking line search parameters $\beta, \gamma$ satisfy $\beta \in (0, 1)$ and $\gamma > 0$. Then the homogeneous second-order descent method (HSODM) with the backtracking line search terminates in at most $O\left(\epsilon^{-3/2}\log_\beta(\epsilon) \right)$ steps, and the next iterate $x_{k+1}$ is a SOSP. Specifically, the number of iterations is bounded by,
    \begin{equation*}
        O\left(\max\left\{\frac{2(M+\gamma)}{9\gamma\beta^3},\frac{3M^3}{4\gamma},\frac{2(M+\gamma)^3}{9\gamma\beta^3}\right\}\left \lceil \log_\beta\left(\frac{3\sqrt{\epsilon}\nu}{M+\gamma}\right) \right\rceil\left(f(x_1) - f_{\inf}\right)\epsilon^{-3/2}\right).
    \end{equation*}
\end{theorem}
\begin{proof}
    ~Since we take $\delta = \sqrt{\epsilon}$ and $\Delta = 2\sqrt{\epsilon} / M $, by \cref{corrollary. descend for line search}, we immediately obtain that the function value decreases at least $\Omega(\epsilon^{3/2})$ for the large step case, i.e.,
    \begin{equation*}
        \begin{aligned}
            f(x_{k+1}) - f(x_k) & \leq -\min\left\{\frac{\sqrt{3}\gamma}{16}, \frac{9\gamma\beta^3\delta^3}{2(M+\gamma)}, \frac{\gamma\Delta^3}{6}, \frac{9\gamma\beta^3\delta^3}{2(M+\gamma)^3}\right\} \\
                                & \leq -\min\left\{\frac{9\gamma\beta^3}{2(M+\gamma)},\frac{4\gamma}{3M^3},\frac{9\gamma\beta^3}{2(M+\gamma)^3}\right\}\epsilon^{3/2},
        \end{aligned}
    \end{equation*}
    and the inner iteration for backtracking line search is at most
    $$
        j_N \leq \left \lceil \log_\beta\left(\frac{3\delta\nu}{M+\gamma}\right) \right\rceil = \left \lceil \log_\beta\left(\frac{3\sqrt{\epsilon}\nu}{M+\gamma}\right) \right\rceil.
    $$
    When the algorithm terminates,  similar to \eqref{eq.gk+1.conv} and \eqref{eq.hk+1.conv}, we have
    \begin{equation*}
        \|g_{k+1}\|\leq O(\epsilon) \quad \text{and} \quad  \lambda_1 \left(H_{k+1}\right) \geq \Omega(-\sqrt{\epsilon}).
    \end{equation*}
    Therefore, the next iterate $x_{k+1}$ is already a SOSP.
    Note that the total decreasing amount of the objective function value cannot exceed $f(x_1) - f_{\inf}$. Hence, the number of iterations for large step case is upper bounded by
    \begin{equation*}
        O\left(\max\left\{\frac{2(M+\gamma)}{9\gamma\beta^3},\frac{3M^3}{4\gamma},\frac{2(M+\gamma)^3}{9\gamma\beta^3}\right\}\left \lceil \log_\beta\left(\frac{3\sqrt{\epsilon}\nu}{M+\gamma}\right) \right\rceil\left(f(x_1) - f_{\inf}\right)\epsilon^{-3/2}\right),
    \end{equation*}
    which is also the iteration complexity of our algorithm. Since $\beta < 1$, this completes the proof.
\end{proof}

Since $\delta = \sqrt{\epsilon}$, we see that the line-search version has an extra overhead of $O(\log_\beta\epsilon)$ compared to the fixed-radius strategy. In practice, the line-search version can choose steps that are much larger than $\Delta$, and thus has a fast rate of convergence. This benefit can be observed in the \cref{sec.experiment}.

\section{Local Convergence Rate}\label{sec.local}

In this section, we provide the local convergence analysis of HSODM. In particular, when $x_k$ is sufficiently close to a SOSP $x^*$,  we will show that the stepsize $\eta_k$ always equals $1$, and the line search procedure is not required.
Consequently, HSODM achieves a local quadratic convergence rate by setting the perturbation parameter $\delta = 0$ for the subsequent iterations.

% Instead of terminating \cref{alg.main alg} as soon as $\epsilon$-approximate second-order conditions \eqref{eq.approxfocp} and \eqref{eq.approxsocp} are met, we set the perturbation parameter $\delta=0$ and let HSODM continue. We show that the HSODM achieves a local quadratic convergence rate.

We first make the standard assumption \cite{connTrustRegionMethods2000,nocedal_numerical_2006,nesterovLecturesConvexOptimization2018} to facilitate the local convergence analysis.
\begin{assumption}
    \label{assm.local}
    Assume that HSODM converges to a strict local optimum $x^*$ satisfying that $\nabla f(x^*)=0$ and $\nabla^2 f(x^*) \succ 0$.
\end{assumption}

\begin{remark}
    From \cref{assm.local}, we immediately know that there exists a small neighborhood for some $R > 0$ and $\mu > 0$ such that
    \begin{equation} \label{eq.R.def}
        \forall x \in B(x^*, R) \quad\Rightarrow\quad \nabla^2 f(x) \succeq \mu \cdot I.
    \end{equation} In other words, $x_k$ arrives at the neighborhood of $x^*$ for some sufficiently large $k$, hence both $\Hk
    $ and $\Hk+\theta_k I$ are nonsingular.
\end{remark}

To prove the local convergence rate, we need the following auxiliary results for preparation.
\begin{corollary}
    \label{coro. large k t neq 0}
    Suppose that \cref{assm.local} holds, then $t_k \neq 0$ for sufficiently large $k$.
\end{corollary}
\begin{proof}
    We prove this by contradiction. Suppose that $t_k = 0$. Then by \cref{corollary. foc}, $(-\theta_k, v_k)$ is the eigenpair of $H_k$, implying that,
    $$
        \lambda_1(H_k) \leq -\theta_k.
    $$
    Recall that in \cref{lemma.optimal condition of subproblem}, we have $\theta_k > 0$, hence $\lambda_1(H_k) < 0$.
    % On the other hand, with the M-Lipschitz continuity of $\nabla^2 f$, we have
    % $$
    % -M  \|\xk - x^*\| \cdot I\preceq \Hk - \nabla^2 f (x^*) \preceq M \|\xk - x^*\|\cdot I,
    % $$
    % which leads to $\Hk \succeq \left(\mu - M \|\xk - x^*\|\right) I \succeq \frac{\mu}{2} \cdot I \succeq 0$ for sufficiently large $k$ by \cref{assm.local}. 
    % %Thus, we must have $\Hk \succeq 0$ as $\xk \to x^*$ for sufficiently large $k$.
    This contradicts $H_k \succ 0$. The proof is then completed.
\end{proof}

The following lemma demonstrates that the step $d_k$ generated by the HSODM eventually reduces to the small valued case for sufficiently large $k$. Consequently, we choose $\eta_k = 1$ and update the iteration by $x_{k+1} = x_k + d_k$ as shown in \cref{subsection.convergence}. We remark that it is similar to the case of the classical Newton trust-region method (see \cite[Theorem 4.9]{nocedal_numerical_2006}), where the updates become asymptotically similar to the pure Newton step.
\begin{lemma}\label{lemma. large k small step}
    For sufficiently large $k$, we have $\|\dk\|\leq \Delta$.
\end{lemma}
\begin{proof}
    Due to $t_k \neq 0$, by equation \eqref{eq.homoeig.foc t neq 0} in \cref{corollary. foc}, we have
    \begin{equation*}
        \dk = -(\Hk+\theta_k I)^{-1}\gk,
    \end{equation*}
    and further
    \begin{align}
        \nonumber\|\dk\|      & \leq \| (\Hk+\theta_k I)^{-1}\| \|\gk\|                       \\
        \label{eq.boundnormv} & \leq \frac{\|\gk\|}{\mu +\theta_k}  \leq \frac{\|\gk\|}{\mu}.
    \end{align}
    The above inequalities hold because of $\Hk \geq \mu I $ and $\theta_k>0$. Note that with \cref{assm.local}, $\|g_k\| \to 0$ as $k\to \infty$, then there exist a sufficiently large $K \ge 0$, such that
    \begin{equation}\label{eq.localboundg}
        \|\gk\|\leq \Delta \mu, \forall k \ge K. %\frac{2\mu\sqrt{\epsilon}}{M}.
    \end{equation}
    Combining \eqref{eq.boundnormv}, we conclude that $\|\dk\|\leq \Delta$ will be satisfied.
\end{proof}
In the local phase, we set the perturbation parameter $\delta=0$ and solve
\begin{equation}
    \label{eq.homo local subproblem}
    \begin{aligned}
        \min_{\|[v; t]\| \le 1} \psi_k(v, t; 0) := ~ &
        \begin{bmatrix}
            v \\ t
        \end{bmatrix}^T
        \begin{bmatrix}
            \Hk   & \gk \\
            \gk^T & 0
        \end{bmatrix}
        \begin{bmatrix}
            v \\ t
        \end{bmatrix}.
    \end{aligned}
\end{equation}
We also denote by $[v_k; t_k]$ the optimal solution to \eqref{eq.homo local subproblem}. Having gathered the above results, we are ready to prove the following theorem.
\begin{theorem}
    \label{theorem:local-quadratic}
    Suppose that \cref{assm.lipschitz} and \cref{assm.local} hold. For sufficiently large $k$, the HSODM converges to \(x^*\) quadratically, that is,
    $$
        \|x_{k+1} - x^*\| \le \left( \frac{M}{\mu} + \frac{\Delta (MR + \mu)^2}{\mu^2\left( 1- \Delta^2 \right)^2 } \right) \|x_k - x^*\|^2.
    $$
    % $$
    %     \|{x_{k+1}-x^*} \| \le O(\|\xk - x^*\|^{2}).
    % $$
    where $R$ is defined as in \eqref{eq.R.def}.
\end{theorem}
\begin{proof}
    By \cref{coro. large k t neq 0}, we have $t_k \neq 0$. Since we take $\delta = 0$, we have the equation \eqref{eq.homoeig.foc t neq 0} in \cref{corollary. foc}, we have
    \begin{equation*}
        g_k^T d_k = -\theta_k \quad \text{and} \quad (H_k + \theta_k I) d_k = - g_k,
    \end{equation*}
    implying that
    \begin{align}
        \nonumber\| H_k^{-1} g_k + d_k \| & =  \| - \theta_k H_k^{-1} d_k \|                   \\
        \nonumber                         & \leq  \| H_k^{-1}  \| \cdot | \theta_k | \| d_k \| \\
        \label{eq:bridge5}                & \leq \frac{1}{\mu} \| g_k\| \| d_k \|^2.
    \end{align}
    By \cref{lemma. large k small step}, we have $ x_{k+1} = x_k + \dk $. Therefore,
    \begin{subequations}
        \begin{align}
            \nonumber   \|x_{k+1} - x^*\| & = \|x_{k} + d_k + H_k^{-1} g_k - H_k^{-1} g_k - x^* \|              \\ \nonumber
                                          & \leq \|x_{k} - H_k^{-1}g_k - x^*\| + \|H_k^{-1}g_k + \dk\|          \\
            \label{eq.local.opt1}         & \leq \frac{M}{\mu}\|x_k - x^*\|^2 + \frac{1}{\mu} \|g_k\| \|\dk\|^2 \\
            \label{eq.local.opt2}         & \leq \frac{M}{\mu}\|x_k - x^*\|^2 + \Delta \|d_k\|^2,
            % \label{eq.local.opt3} & \leq \frac{M}{\mu}\|x_k - x^*\|^2 + \frac{L}{\mu}\|x_k - x^*\| \|d_k\|^2
        \end{align}
    \end{subequations}
    where \eqref{eq.local.opt1} holds due to the standard analysis of Newton's method \cite{nocedal_numerical_2006} and equation \eqref{eq:bridge5}, and \eqref{eq.local.opt2} follows from $ \|\gk\|\leq \Delta \mu $ as stated in \cref{lemma. large k small step}. Moreover, we have
    \begin{align}
        \nonumber        \|\dk\| & = \|  x_{k+1} - x^*  - \left( \xk - x^* \right) \|             \\
        \nonumber                & \leq \| x_{k+1} - x^* \|+\|\xk-x^*\|                           \\
        \nonumber                & \leq \frac{M}{\mu}\|\xk-x^*\|^2+\|\xk-x^*\|+\Delta \|d_k\|^2   \\ \nonumber
                                 & \leq  \frac{MR}{\mu}\|\xk-x^*\| +\|\xk-x^*\|+\Delta^2 \|d_k\|,
    \end{align}
    where the last inequality holds since $x_k \in B(x^*, R)$ and $\|d_k\| \le \Delta$. Rearranging the terms implies
    \begin{equation*}
        \|d_k\| \leq  \frac{MR + \mu}{\mu(1-\Delta^2)} \|\xk-x^*\| .
    \end{equation*}
    With \eqref{eq.local.opt2}, we conclude that
    \begin{equation*}
        \|x_{k+1} - x^*\| \le \frac{M}{\mu}\|x_k - x^*\|^2 + \Delta \|d_k\|^2 \le \left( \frac{M}{\mu} + \frac{\Delta (MR + \mu)^2}{\mu^2\left( 1- \Delta^2 \right)^2 } \right) \|x_k - x^*\|^2.
    \end{equation*}

    % By rearranging the terms and \cref{lemma. large k small step} we have
    % $$    (1-\Delta^2)\|\dk\| \leq O(\|\xk-x^*\|)+\|\xk-x^*\|,$$
    % which leads to that $\|\dk\|\leq O(\|\xk-x^*\|)$. With \eqref{eq.local.opt2}, we conclude that
    % \begin{equation}
    %     \begin{aligned}
    %         \|x_{k+1} - x^*\| & \le \frac{M}{\mu}\|x_k - x^*\|^2 + \Delta \|d_k\|^2 = O(\|x_k - x^*\|^2).
    %     \end{aligned}
    % \end{equation}
    This completes the proof.
\end{proof}

%% inexact

\section{An Inexact HSODM}\label{sec.inexact}
The above analysis relies on solving the subproblem \eqref{eq.homo subproblem} exactly, which requires matrix factorization with $O((n+1)^3)$ arithmetic operations. In this section, we propose an inexact HSODM (\cref{alg inexact HSODM}), which utilizes a Lanczos method (\cref{alg.randlanczos}) to approximately solve \eqref{eq.homo subproblem} in each iteration. After that, we construct the iterates based on the Ritz pair of $F_k$ instead of its exact leftmost eigenpair. We will prove later that this method provides a probabilistic worst-case arithmetic operation of $\tilde O((n+1)^2\epsilon^{-7/4})$, which has less dependence on $n$.

\subsection{A brief overview of the Lanczos method}
Before delving into the details, we briefly introduce the Lanczos method, which is utilized to compute the extremal eigenvalue of a symmetric matrix $A\in \real^{n\times n}$. At $j$-th iteration, the Lanczos method constructs an orthonormal basis $Q_j=[q_1,q_2,\ldots,q_j] \in \real^{n \times j}$ from $j$-th Krylov subspace $\cK(j;A,q_1) := \text{span} \{q_1, Aq_1, \ldots, A^{j-1} q_1\}$, keeping $T_j=Q_j^TAQ_j$ tridiagonal at the same time.  The next lemma provides some standard results of the Lanczos method.
\begin{lemma}[Basic properties of the Lanczos method \cite{golubMatrixComputations2013}]
    \label{lem.basic.lanczos}
    For any symmetric matrix $A\in \real^{n\times n}$, let $q_1\in \real^n$ and $\|q_1\|=1$. Suppose that the Lanczos method runs until iteration $J=\rank(\cK(n;A,q_1))$, then the following statements hold:
    \begin{enumerate}[(1)]
        \item For any $j=1,2,\ldots,J$, let $Q_j=[q_1,q_2,\ldots,q_j]$ be the orthonormal basis that spans $\cK(j;A,q_1)$, then
              \begin{equation*}
                  AQ_j=Q_j T_j+\xi_j (1_j)^T_{[1:j]} \quad \text{and} \quad Q_j\perp \xi_j,
              \end{equation*}
              where $T_j=Q_j^TAQ_j$ is a tridiagonal matrix, $1_j \in \real^n$ is the $j$-th column of $I_n$, and $\xi_j$ is the residual vector.
        \item Suppose $Y_j = Q_j S_j$ are computed from the $j$-th Krylov iteration of the Lanczos method and the real Schur decomposition $S_j^T T_jS_j = \Gamma_j$. Let $\gamma_i$ be the $i$-th entry on the diagonal of $\Gamma_j$, $y_i$ be the $i$-th column vector of $Y_j$, then we have the following error estimation:
              \begin{align*}
                  A y_i - \gamma_i y_i =  (1_j)^T_{[1:j]}S_j (1_i)_{[1:j]}\cdot \xi_j := s_{ji} \cdot  \xi_j\;\,\mbox{with}\;\, |s_{ji}| < 1 \;\, \mbox{such that}\;\,y_i \perp \xi_j,\;\, \forall i \le j.
              \end{align*}
              We call $(\gamma_i, y_i)$ the $i$-th Ritz pair.
    \end{enumerate}
\end{lemma}
For the rest of the paper, we sometimes omit the indexing ${[1:j]}$ for simplicity. It is understood that the matrix-vector operations are compatible in size.
With a slight abuse of notation, we let $[v_k;t_k]$ be the approximate solution. We still let $ -\theta_k =\lambda_1(F_k)$ be the smallest eigenvalue of $F_k$, and denote its eigenvector by $\chi_k$.
\begin{theorem}[Property of the approximate solution]
    \label{lemma.property of Lanczos}
    Suppose that the Lanczos method is used to approximately solve \eqref{eq.homo subproblem} and returns a Ritz pair $(-\gamma_k, [v_k;t_k])$. We have
    \begin{subequations}\label{eq.inexact optimal condition}
        \begin{align}
            \label{eq.inexact.opt}
            \begin{bmatrix} H_k   & g_k     \\g_k^T & -\delta\end{bmatrix}
            \begin{bmatrix}v_k \\t_k\end{bmatrix}+ \gamma_k \begin{bmatrix} v_k \\t_k\end{bmatrix}
             & = \begin{bmatrix} r_k\\ \sigma_k \end{bmatrix}, \\
            \label{eq.inexact.orth}      r_k^T v_k + \sigma_k \cdot t_k
             & = 0.
        \end{align}
    \end{subequations}
    where $[r_k; \sigma_k] \in \real^{n}\times\real$ is called the Ritz error.
\end{theorem}

The above theorem is a direct application of part (2) of \cref{lem.basic.lanczos}. Since $(-\gamma_k, [v_k;t_k])$ is only an approximate solution, we consider some error estimates $e_k > 0$ such that $|\theta_k - \gamma_k| \le e_k$. In the Lanczos method, $-\gamma_k$ is always an overestimate of $-\theta_k$ \cite{golubMatrixComputations2013}, thus we stop at $\theta_k - e_k \leq \gamma_k \le \theta_k$. We provide the following complexity estimates regarding a prescribed error $e_k$.

\begin{lemma}[Complexity of the Lanczos method]
    \label{lemma.complexity of lanczos method}
    Suppose that the Lanczos method is used to approximately solve \eqref{eq.homo subproblem}, and returns a Ritz pair $(-\gamma_k, [v_k;t_k])$ satisfying $\theta_k - e_k   \leq   \gamma_k \le \theta_k$ for some $e_k > 0$. Then, the number of required iterations can be upper-bounded by either of the following quantities.
    \begin{enumerate}[(1)]
        \item
              \begin{equation}\label{eq.gap-free}
                  1+\left\lceil 2\sqrt{\frac{\|F_k\|}{e_k}}\log\left(\frac{16\|F_k\|}{e_k(q_1^T \chi_k)^2}\right)\right\rceil,
              \end{equation}
              where $(-\theta_k, \chi_k)$ is the exact leftmost eigenpair of $F_k$ \cite{kuczynski_estimating_1992,royer_complexity_2018};
        \item
              \begin{equation}\label{eq.gap-dependent}
                  1 + \left\lceil\sqrt{\frac{2\|F_k\|}{\lambda_2(F_k) - \lambda_1(F_k)}}\log\left(\frac{8\|F_k\|}{e_k(q_1^T\chi_k)^2}\right)\right\rceil,
              \end{equation}
              where $\lambda_2(F_k)$ is the second-smallest eigenvalue of $F_k$ such that $\lambda_2(F_k) - \lambda_1(F_k) > 0$ \cite{kuczynski_estimating_1992}.
    \end{enumerate}
\end{lemma}
We also remark that the Lanczos method has finite convergence. Finally, we connect the Ritz error to the desired accuracy $e_k$.
\begin{lemma}
    \label{lem.convergence r sigma}
    Suppose that \cref{assm.lipschitz} holds, and $F_k$ is constructed {as in \eqref{eq.motivate}}, then
    \begin{align}
        \|F_k\|
        \label{eq.est.Fk.1} & \le \max\{U_H, \delta\} + \|\gk\|.
    \end{align}
    If we let  $\varsigma_k := \lambda_2(F_k)-\lambda_1(F_k)>0$, then for $[r_k;\sigma_k]$ in \eqref{eq.inexact optimal condition}, there exists $\tau_k \in [0,1]$ such that
    \begin{equation}
        \label{eq.bound r sigma in terms of e}
        \| [r_k;\sigma_k] \| \leq \tau_k e_k +2(\max\{U_H, \delta\} + \|\gk\|) \sqrt{\frac{e_k}{{\varsigma_k}}}.
    \end{equation}
\end{lemma}
We defer the proofs of \cref{lemma.complexity of lanczos method} and \cref{lem.convergence r sigma} to the Appendix as the results are mostly related to linear algebra.

\subsection{Overview of the inexact HSODM}

Now, we are ready to introduce the inexact HSODM in \cref{alg inexact HSODM}. It follows the basic idea of the exact HSODM but uses the Lanczos method to approximately solve \cref{eq.homo subproblem}. The inexactness brings several challenges to establishing the corresponding convergence result. First, since $\gamma_k$ in the Ritz pair is an inexact dual variable, we cannot guarantee that $\gamma_k$ exceeds $\delta$, which may result in an insufficient descent property. Second, the large Ritz error in the small value case (when $t_k>\sqrt{1/(1+\Delta^2)}$) may prevent the next iterate $x_{k+1}$ from being the SOSP when we update via $x_{k+1}=x_k+d_k$.

To overcome the first challenge, we propose a customized Lanczos method (\cref{alg.randlanczos}) with skewed randomization, which ensures that $\gamma_k$, in high probability, is always no smaller than $\delta$ (cf. \cref{thm.gamma ge delta}, \cref{thm.customized initialization}). For the second challenge, we discuss the magnitude of $\|r_k\|$. If $\|r_k\|$ is sufficiently small, we safely claim that $x_{k+1}=x_k+d_k$ is already a SOSP (\cref{lemma.inexact small step}). Otherwise, we increase the perturbation parameter $\delta$ and solve the subproblem \cref{eq.homo subproblem}.
By a delicate analysis of the spectrum, we show that the eigengap $\varsigma_k$ of the homogenized matrix $F_k$ is sufficiently large (e.g., in $\Omega(\sqrt \epsilon)$). This implies that it is possible to pursue a higher precision (\cref{line.increase.pert}) indicated by the gap-dependent complexity \eqref{eq.gap-dependent}.

\begin{minipage}[ht]{0.95\linewidth}
    \begin{algorithm}[H]
        \caption{Inexact Homogeneous Second-Order Descent Method}\label{alg inexact HSODM}
        \KwIn{Initial point $x_1$, $\nu \in (1/4, 1/2)$, $\Delta  = \sqrt{\epsilon}/{M}$, $\epsilon > 0$.}
        \For{$k = 1, 2, \cdots$}{
            Set  $\delta \leftarrow \sqrt{\epsilon}$, $e_k \leftarrow \sqrt{\epsilon}$, $J_{\max} \leftarrow n+1$\;
            Run \cref{alg.randlanczos} with $(\delta, e_k, J_{\max})$ to obtain the Ritz pair $(\gamma_k, [v_k; t_k])$ and the Ritz error $[r_k;\sigma_k]$\label{line.solve}\;
            \uIf(\tcp*[f]{small value case}\label{line.inexact.small}){$|t_k| > \sqrt{1/(1+\Delta^2)}$}{
                \eIf{$\|r_k\| \le 2\epsilon$}{
                    \label{line.largerde}
                    Set $\dk \leftarrow v_k / t_k$\;
                    Update $x_{k+1} \leftarrow x_k + d_k$\; \label{line.small.begin}
                    (Early) Terminate (or set $\delta = 0$ and proceed)\;\label{line.small.finish}
                }{
                    Set $\delta \leftarrow 3\sqrt\epsilon + 2\|\gk\|\Delta + (U_H + \gamma_k)\Delta^2$, $e_k = \min\left\{\epsilon,\frac{\epsilon^{\frac{5}{2}}}{4(U_H + U_g)^2}\right\}$\;\label{line.increase.pert}
                    Go to \cref{line.solve}\;
                }}
            \eIf(\tcp*[f]{large value case (a)}\label{line.inexact.largea}){$|t_k| \geq \nu$}{
                Set $\dk \leftarrow v_k / t_k$\;
            }
            (\tcp*[f]{large value case (b)}\label{line.inexact.largeb}){
                Set $\dk \leftarrow \text{sign}(-g_k^T v_k) \cdot v_k$\;
            }
            Choose a stepsize $\eta_k$ by fixed-radius strategy\;
            Update $x_{k+1}\leftarrow x_k + \eta_k \cdot d_k$\;
        }
    \end{algorithm}
\end{minipage}

In the following of this subsection, we analyze the descent properties under the large value cases (a) and (b) in the inexact HSODM (\cref{line.inexact.largea} and \cref{line.inexact.largeb} in \cref{alg inexact HSODM}). They follow in a similar manner to those in the exact HSODM, and our analysis shows that the inexactness indeed brings obstacles to the convergence analysis.

\begin{lemma}[Large value case (a)]
    \label{lemma.inexact large step v/t}
    Suppose that \cref{assm.lipschitz} holds and set $\nu \in (1/4, 1/2)$. If $|t_k| \ge \nu$ and $\|v_k/t_k\| \ge \Delta$, then let $ \dk = v_k/t_k$ and $\eta_k = \Delta/\|\dk\|$, we have
    \begin{equation*}
        f(x_{k+1}) - f(x_k) \le \left(\eta_k - \frac{1}{2}\eta_k^2\right)\left(\delta - \gamma_k\right) + 4\vert \sigma_k \vert - \frac{\gamma_k}{2} \Delta^2 + \frac{M}{6}\Delta^3.
    \end{equation*}
\end{lemma}
\begin{proof}
    By \eqref{eq.inexact.opt} and $d_k = v_k/t_k$, we have
    \begin{align*}
         & d_k^T H_k d_k + g_k^T d_k = -\gamma_k \|d_k\|^2 + \frac{r_k^T v_k}{t_k^2}, \\
         & g_k^T d_k = -\gamma_k + \delta + \frac{\sigma_k}{t_k}.
    \end{align*}
    Therefore, we obtain
    \begin{subequations}
        \begin{align}
            \nonumber      f(x_{k+1}) - f(x_k) = ~ & f(x_k + \eta_k \cdot d_k) - f(x_k)                                                                                                                                                                          \\
            \nonumber         \le                ~ & \eta_k \cdot g_k^T d_k + \frac{\eta_k^2}{2}\cdot d_k^T H_k d_k + \frac{M\eta_k^3}{6}\cdot \|d_k\|^3                                                                                                         \\
            \nonumber          =                 ~ & \eta_k \cdot \gk^T\dk + \frac{1}{2}\eta_k^2\left(\frac{r_k^Tv_k}{t_k^2} - \gk^T\dk - \gamma_k \|\dk\|^2\right) + \frac{M\eta_k^3}{6}\cdot \|d_k\|^3                                                         \\
            \nonumber        =                   ~ & \left(\eta_k - \frac{1}{2}\eta_k^2\right)\left(\frac{\sigma_k}{t_k} + \delta - \gamma_k\right) + \frac{\eta_k^2}{2}\left(\frac{r_k^Tv_k}{t_k^2}\right)  - \frac{\gamma_k}{2} \Delta^2 + \frac{M}{6}\Delta^3 \\
            \nonumber=                           ~ & \left(\eta_k - \frac{1}{2}\eta_k^2\right)\left(\delta - \gamma_k\right) -  \left(\eta_k^2 - \eta_k\right)\frac{\sigma_k}{t_k} - \frac{\gamma_k}{2} \Delta^2 + \frac{M}{6}\Delta^3.
        \end{align}
    \end{subequations}
    The last equality holds by \cref{eq.inexact.orth}. Since $\eta_k \in (0, 1)$, $|t_k| \ge \nu$ and $\nu \ge 1/4$, then it holds that
    \begin{equation*}
        - \left(\eta_k^2 - \eta_k\right)\frac{\sigma_k}{t_k} \le \left|\frac{\sigma_k}{\nu}\right| \le 4|\sigma_k|.
    \end{equation*}\
    Finally, we conclude
    \begin{equation*}
        f(x_{k+1}) - f(x_k) \le \left(\eta_k - \frac{1}{2}\eta_k^2\right)\left(\delta - \gamma_k\right) + 4\vert \sigma_k \vert - \frac{\gamma_k}{2} \Delta^2 + \frac{M}{6}\Delta^3.
    \end{equation*}
\end{proof}

\begin{lemma}[Large value case (b)]
    \label{lemma.inexact large step truncated}
    Suppose that \cref{assm.lipschitz} holds and set $\nu \in (1/4, 1/2)$. If $|t_k| \le \nu$, then let $ \dk = \text{sign}(-g_k^T v_k) \cdot v_k$ and $\eta_k = \Delta/\|\dk\|$, we have
    \begin{equation*}
        f(x_{k+1}) - f(x_k) \le \vert \sigma_k\vert -\frac{\gamma_k}{2}\Delta^2 +\frac{M}{6}\Delta^3.
    \end{equation*}
\end{lemma}

\begin{proof}

    From \eqref{eq.inexact.opt}, we obtain
    \begin{align*}
         & v_k^T H_k v_k = r_k^T v_k - \gamma_k \|v_k\|^2 - t_k g_k^T v_k, \\
         & g_k^T v_k = \sigma_k + t_k \cdot (\delta - \gamma_k).
    \end{align*}

    Consequently, it follows that
    \begin{subequations}
        \begin{align}
            \nonumber        f(x_{k+1}) - f(x_k) = & f(x_k + \eta_k \cdot d_k) - f(x_k)                                                                                                                                    \\
            \nonumber         \le                  & \eta_k \cdot g_k^T d_k + \frac{\eta_k^2}{2}\cdot d_k^T H_k d_k + \frac{M\eta_k^3}{6}\cdot \|d_k\|^3                                                                   \\
            \nonumber          =                   & \eta_k \cdot \operatorname{sign}(-\gk^T v_k) \gk^T v_k +\frac{1}{2}\eta_k^2(v_k)^T H_k v_k +\frac{M}{6}\eta_k^3\|v_k\|^3                                              \\
            \nonumber =                            & -\eta_k \cdot |g_k^T v_k| + \frac{1}{2}\eta_k^2 r_k^T v_k-\frac{1}{2}\eta_k^2t_k \gk^T v_k -\frac{1}{2}\eta_k^2 \gamma_k\|v_k\|^2 +\frac{M}{6}\eta_k^3\|v_k\|^3
            \\
            \nonumber \le                          & -\eta_k \cdot |g_k^T v_k| +\frac{1}{2}\eta_k^2 r_k^T v_k + \frac{1}{2}\eta_k^2 |t_k| |g_k^T v_k| -\frac{1}{2}\eta_k^2 \gamma_k\|v_k\|^2 +\frac{M}{6}\eta_k^3\|v_k\|^3
            \\
            \nonumber =                            & -\frac{1}{2}\eta_k^2 t_k \sigma_k  -\left(\eta_k - \frac{1}{2}\eta_k^2|t_k|\right)|g_k^T v_k| -\frac{\gamma_k}{2}\Delta^2 +\frac{M}{6}\Delta^3,
        \end{align}
    \end{subequations}
    where the last equality holds due to \eqref{eq.inexact.orth} and $ \eta_k \|v_k\| =  \eta_k \|d_k\| = \Delta$. Since $\eta_k < 1$ and $|t_k| \le \nu < 1$, we have $\eta_k^2 |t_k| \le \eta_k < 1$, and thus
    \begin{equation*}
        f(x_{k+1}) - f(x_k) \le \vert \sigma_k\vert -\frac{\gamma_k}{2}\Delta^2 +\frac{M}{6}\Delta^3.
    \end{equation*}
\end{proof}
The above two lemmas illustrate how the Ritz error $[r_k;\sigma_k]$ and the inexact dual variable $\gamma_k$ obstruct the descent property. To ensure the convergence of the inexact HSODM, the Lanczos method should guarantee $\gamma_k \ge \delta$ and provide a sufficiently small Ritz error. However, the classical Lanczos method with random start \cite{kuczynski_estimating_1992} cannot satisfy the need. In the next subsection, we propose a customized Lanczos method with skewed randomization to overcome this challenge, which may be of independent interest.

We close this subsection by introducing the following assumption, which is widely adopted in the analysis of second-order algorithms \cite{cartis_adaptive_2011,royer_complexity_2018}.
\begin{assumption}\label{assum.bounded gradient}
    Assume that there exists a constant $U_g > 0$ independent of $k$, such that $\|\nabla f(x_k)\| \le U_g, ~ \forall k \ge 1$.
\end{assumption}
Since the inexact HSODM is monotone (as established later in \cref{thm.inexact complexity}), the above assumption can be easily satisfied whenever the sublevel set $\{x: f(x) \le f(x_1)\}$ is compact. According to \cref{lem.convergence r sigma}; this assumption implies that $U_H + U_g$ serves as an upper bound of $\|F_k\|$, which is necessary to establish the properties of the customized Lanczos method in \cref{thm.gamma ge delta}.

\subsection{A customized Lanczos method with skewed randomization}
In this subsection, we develop a Lanczos method with skewed randomization, which allows us to attain a convergence behavior akin to that of the exact HSODM. The crux of our Lanczos method lies in the skewed randomization of the initial vector $q_1$  (\cref{line.init} in \cref{alg.randlanczos}). The basic idea is to assign a greater weight to the last entry of $q_1$. Namely, we first sample $b_i$ i.i.d. from a standard normal distribution $\mathcal{N}(0, 1), ~i = 1, \ldots, n+1,$ and multiply the last entry $b_{n+1}$ with a large constant $\Psi_k$. Let $b = [b_1, \cdots, b_n, \Psi_k \cdot b_{n+1}]^T$, then we choose the normalized vector $q_1 := b/\|b\|$ as the initial vector for the Lanczos method.

\begin{minipage}[h]{0.95\textwidth}
    \begin{algorithm}[H]
        \caption{A Lanczos Method with Skewed Randomization}\label{alg.randlanczos}
        \KwIn{Iterate $x_k$, $\gk$, $\Hk$; $\delta > 0$, $p \in (\exp(-n),1), e_k > 0, J_{\max} \ge 0$}
        \textbf{Initialization:} {sample $b_1, b_2, \ldots, b_{n+1}$ i.i.d. from a standard normal distribution $\mathcal{N}(0,1)$\;
        Set $\Psi_k$ by \eqref{eq.psi.choice}, $b: = [b_1, \cdots, b_n, \Psi_k \cdot b_{n+1}]^T$} and $q_1 = b / \|b\|$\label{line.init}\;
        Construct $F_k$ with $\chi_k $ being its exact leftmost eigenvector and let $J_m = \min\left\{J_{\max}, 1 + \sqrt{\frac{2\|F_k\|}{e_k}}\log\left(\frac{8}{e_k(q_1^T\chi_k)^2}\right)\right\}$\label{line.setmaxiter}\;
        \While{$j =1,..., J_m$}{
        Compute $F_k Q_j = Q_j T_j + \xi_j (1_j)^T_{[1:j]}$\;
        \uIf{$\|\xi_j\| \le \epsilon$}{Break\;}
        $j \leftarrow j+1$\;
        }
        Compute Schur decomposition of $T_j$ such that $S_j^T T_jS_j = \Gamma_j$\;
        Compute Ritz approximation $(-\gamma_k, [v_k; t_k])$\;
        \Return{$(-\gamma_k, [v_k; t_k])$ and the corresponding Ritz error $[r_k;\sigma_k]$}
    \end{algorithm}
    \vspace{1em}
\end{minipage}
For ease of theoretical analysis, since $\|q_1\| = 1$, one can rewritte it as $q_1 :=  \sqrt{1-\alpha^2} \cdot [u;0] + \alpha \cdot [0;1] \in \mathbb{R}^{n+1}$, where $u \in \mathbb{R}^{n}$ and $\|u\| = 1$. The following theorem shows that when $|\alpha|$ exceeds a certain threshold, the inequality $\gamma_k \ge \delta$ is guaranteed. Surprisingly, the magnitude of the last entry in the Ritz error can also be bounded by $|\alpha|$.
\begin{theorem}
    \label{thm.gamma ge delta}
    Suppose that \cref{assm.lipschitz} and \cref{assum.bounded gradient} hold. For the homogenized matrix $F_k$, suppose that the Lanczos method is run with the initial vector $q_1:=  \sqrt{1-\alpha^2} \cdot [u;0] + \alpha \cdot [0;1] \in \mathbb{R}^{n+1}$, where $u \in \mathbb{R}^{n}$ and $\|u\| = 1$, then for any $|\alpha| \ge 1/2$, the following statements holds:
    \begin{enumerate}[(1)]
        \item After the $j$-th iteration ($j \ge 2$), the last entry of the Lanczos vector $q_j = [\ell_j; \beta_j] \in \real^n \times \real$ is bounded, i.e. $|\beta_j| \le 2\sqrt{1-\alpha^2}$.
        \item After the $j$-th iteration ($j \ge 4$), the last entry of the Ritz error $[r_k; \sigma_k]$ is bounded, i.e.
              \begin{equation}\label{eq.safeguard.sigma}
                  |\sigma_k| \le U_\sigma \sqrt{1-\alpha^2},
              \end{equation}
              where $U_\sigma$ is a constant independent of $k$:
              \begin{equation}\label{eq.usigma}
                  U_\sigma := \sqrt{(U_H + U_g)^2 + (\delta+U_g)^2}\sqrt{U_g^2+\delta^2} + 4\sqrt{n}(U_g + \max\{U_H, \delta\}).
              \end{equation}
        \item Suppose that
              \begin{equation}\label{eq.safeguard.gammadelta}
                  \alpha \cdot g_k^T u \le 0 \quad \mbox{and} \quad |\alpha| \ge \frac{U_H+\delta}{\sqrt{(U_H+\delta)^2 + 4(\gk^Tu)^2}},
              \end{equation}
              then the inexact dual variable $\gamma_k$ is sufficiently large, i.e., $\gamma_k \ge \delta$.
    \end{enumerate}
\end{theorem}

Based on the above theorem, we next show \cref{alg.randlanczos} fits the purpose by selecting $\Psi_k$ properly.

\begin{theorem}\label{thm.customized initialization}
    Suppose that \cref{assm.lipschitz} and \cref{assum.bounded gradient} hold. Consider the skewed initialization (\cref{line.init}) in \cref{alg.randlanczos},
    and choose $\Psi_k$ such that
    \begin{equation}\label{eq.psi.choice}
        \Psi_k = \frac{\sqrt{10n}}{\sqrt{\pi}p}\cdot\max\left\{\frac{16M^2U_\sigma}{\epsilon^2}, \sqrt{1+\frac{(U_H+\delta)^2 }{2 p^2 \pi \|\gk\|^2}}, \frac{2}{\sqrt{3}}\right\},
    \end{equation}
    where $U_\sigma$ is defined in \eqref{eq.usigma}.
    Recalling that $\chi_k = [\chi_{k,1},...,\chi_{k,n+1}]$ is the exact leftmost eigenvector of $F_k$,
    then for any constant $p \in (\exp(-n),1)$ and $\epsilon > 0$, with a probability of at least {$1-4p$}, it holds that
    \begin{align}
        \label{eq.psi.bound.ip}    & (q_1^T \chi_k)^2   \ge \min\left\{\frac{\epsilon^4}{256M^4U^2_\sigma}, \left({1+\frac{(U_H+\delta)^2 }{2 p^2 \pi \|\gk\|^2}}\right)^{-1}, \frac{3}{4}\right\} \cdot \frac{\pi^2 p^4 \sum_{i=1}^n\chi_{k,i}^2}{100n(n + 1)} + \frac{p^2\pi \chi_{k, n+1}^2}{10(n + 1)}, \\
        \label{eq.psi.bound.sigma} & |\sigma_k|    \le \frac{\epsilon^2}{16M^2} \quad \text{and} \quad|\alpha| \ge \frac{U_H+\delta}{\sqrt{(U_H+\delta)^2 + 4(\gk^Tb_{[1:n]})^2}}.
    \end{align}
\end{theorem}

We delay the proofs of the above two theorems to the Appendix, {as they are quite technical}.
The above two theorems show that skewed randomization can guarantee sufficiently small $\sigma_k$ with high probability. Due to the symmetry of normal distribution, one can always ensure $\alpha \cdot \gk^Tb_{[1:n]} \le 0$ by flipping the sign of $b_{[1:n]}$, guaranteeing that the inexact dual variable $\gamma_k$ satisfies $\gamma_k \ge \delta$. Furthermore, we show that $(q_1^T\chi_k)^2$ is bounded away from $0$, which generally attains the first term in \eqref{eq.psi.bound.ip} (i.e., in $\Omega(\epsilon^4/n(n+1)$); this enables a later complexity analysis of our method.

\begin{remark}
    Note that \cref{alg.randlanczos} may rely on a priori $\|F_k\|$. Technically, one can slightly refine \cref{alg.randlanczos} with the bound estimation \cite[Algorithm 5]{royer_newton-cg_2020}, in which case $\|F_k\|$ can be estimated by some $\hat F_k$ such that
    \begin{equation}
        \|F_k\| \in [\hat F_k / 2, \hat F_k],
    \end{equation}
    in the first $O(\log(n))$ iterations with high probability (\cite[Lemma 10]{royer_newton-cg_2020}). Then the dependency on a priori $\|F_k\|$ can be removed (\cref{line.setmaxiter} in \cref{alg.randlanczos}) at the cost of one trial run.
\end{remark}

In the following corollary, we show that a sufficient decrease can be achieved in the large value cases by the customized Lanczos method with skewed randomization.
\begin{corollary}
    \label{lemma.inexact function value decrease}
    Suppose that \cref{assm.lipschitz} and \cref{assum.bounded gradient} hold. If we run \cref{alg.randlanczos} and set the parameters $e_k = \delta = \sqrt{\epsilon}$ and $\Delta = \frac{\sqrt{\epsilon}}{M}$. Then for any $\epsilon > 0$ and $p \in (\exp(-n),1)$, under the two large value cases, it holds that
    \begin{equation*}
        f(x_{k+1}) - f(x_k) \le - \frac{\delta}{4} \Delta^2 + \frac{M}{6}\Delta^3
    \end{equation*}
    with a probability of at least ${1-4p}$.
\end{corollary}
\begin{proof}
    From \cref{thm.gamma ge delta} and \cref{thm.customized initialization}, with a probability of at least ${1-4p}$, it holds that
    \begin{equation*}
        |\sigma_k| \le \frac{\epsilon^2}{16M^2} \le \frac{\epsilon^{3/2}}{16M^2} = \frac{\delta}{16}\Delta^2 \quad \text{and} \quad\gamma_k \ge \delta.
    \end{equation*}
    Therefore, the large step case (a) (\cref{lemma.inexact large step v/t}) implies that
    $$f(x_{k+1}) - f(x_k) \le 4|\sigma_k| -\frac{\gamma_k}{2}\Delta^2 +\frac{M}{6}\Delta^3.$$
    Combining the large case (b) (\cref{lemma.inexact large step truncated}):
    \begin{align*}
        f(x_{k+1}) - f(x_k) & \le |\sigma_k| -\frac{\gamma_k}{2}\Delta^2 +\frac{M}{6}\Delta^3,
    \end{align*}
    we have
    \begin{align*}
        f(x_{k+1}) - f(x_k) & \le 4|\sigma_k| -\frac{\gamma_k}{2}\Delta^2 +\frac{M}{6}\Delta^3                                                                    \\
                            & \le  \frac{\delta}{4}\Delta^2 - \frac{\delta}{2} \Delta^2 + \frac{M}{6}\Delta^3 = -\frac{\delta}{4} \Delta^2 + \frac{M}{6}\Delta^3.
    \end{align*}
    This completes the proof. \end{proof} \endproof

\subsection{Small value case in the inexact HSODM}
For the small value case, as before, it occurs when
$|t_k| \geq \nu$ and $d_k = v_k/t_k$. Under this scenario, we show that the Hessian matrix at iterate $x_k$ is nearly positive semidefinite. In this view, \cref{alg inexact HSODM} tests whether the Ritz error $r_k$ is sufficiently small. If not, it increases the perturbation parameter $\delta$ and recalculates the Ritz pair by \cref{alg.randlanczos}. In this case, we show that the eigengap of the homogenized matrix $F_k$ now exceeds $\Omega(\sqrt \epsilon)$. These results are summarized in the following lemma.

\begin{lemma}[Small value case]
    \label{lemma.inexact small step}
    Suppose that \cref{assm.lipschitz} and \cref{assum.bounded gradient} hold. If $|t_k| > \sqrt{1/(1+\Delta^2)}$ and \cref{alg.randlanczos} is run with $e_k = \delta = \sqrt{\epsilon}$ and $\Delta = \frac{\sqrt{\epsilon}}{M}$, where $\epsilon \le \min \{(2M U_g / (2U_H + U_g) )^2, 3M^2, 1\}$, then the following statements hold:
    \begin{itemize}
        \item[(1)] for any $p \in (\exp(-n),1)$, it holds that
            \begin{align*}
                \lambda_1(H_k) \ge -2 \delta - 2\|\gk\|\Delta - (U_H + \gamma_k)\Delta^2 \ge -2\left(1 + \frac{2U_g}{M}\right)\sqrt{\epsilon}
            \end{align*}
            with a probability of at least ${1-4p}$.
        \item[(2)] If the Ritz error $r_k$ satisfies $\|r_k\| \le 2\epsilon$ (\cref{line.largerde}), then the next iterate $x_{k+1} = x_k + d_k$ is already an $\epsilon$-approximate SOSP.
        \item[(3)] Otherwise, when resetting $\delta = 3\sqrt\epsilon + 2\|\gk\|\Delta + (U_H + \gamma_k)\Delta^2$ (\cref{line.increase.pert}),   the eigengap of the resulting homogenized matrix holds that: $\varsigma_k = \lambda_2(F_k) - \lambda_1(F_k) \ge \sqrt{\epsilon}$.
    \end{itemize}
\end{lemma}
\begin{proof}
    From \eqref{eq.inexact optimal condition}, we have
    \begin{align*}
        -\gamma_k = -\delta t_k^2 + 2t_k \gk^Tv_k + v_k^T H_kv_k.
    \end{align*}
    Rearranging the terms gives
    \begin{align*}
        (\gamma_k - \delta) t_k^2 & =- 2t_k \gk^Tv_k - \left(\gamma_k + \frac{v_k^T H_kv_k}{\|v_k\|^2}\right) \|v_k\|^2               \\
                                  & \le 2t_k\sqrt{1-t_k^2} \|\gk\| - \left(\gamma_k + \frac{v_k^T H_kv_k}{\|v_k\|^2}\right) \|v_k\|^2 \\
                                  & \le 2t_k\sqrt{1-t_k^2} \|\gk\| - \left(\gamma_k + \lambda_1(\Hk)\right) (1-t_k^2),
    \end{align*}
    where the first equality holds since $\|v_k\|^2 + t_k^2=1$. This further implies that
    \begin{align*}
        \begin{split}
            \gamma_k - \delta & \le 2 \Delta \|\gk\| + |\lambda_1(\Hk) + \gamma_k|\Delta^2 \\
            & \le 2\Delta\|\gk\| + (U_H + \gamma_k)\Delta^2,
        \end{split}
    \end{align*}
    where the first inequality follows from $ \Delta \geq \sqrt{1-t_k^2}/t_k$.
    Recall that $H_k + \theta_k I \succeq 0$ and $\theta_k \le \gamma_k + e_k = \gamma_k + \delta$, {we further have
            \begin{equation}\label{eq.smalleigenlb}
                \lambda_1(H_k) + \theta_k \geq 0 \quad
                \Rightarrow \quad \lambda_1(H_k) +  2\delta + 2\|\gk\|\Delta + (U_H + \gamma_k)\Delta^2 \ge 0.
            \end{equation}
        }
    Since $\delta = \sqrt{\epsilon}$, $\Delta = \sqrt{\epsilon} / M$, $\|g_k\| \le U_g$ and $\gamma_k \le \|F_k\| \le U_H + U_g$, we conclude
    \begin{align*}
        \lambda_1(H_k) & \ge -2\sqrt{\epsilon} - \frac{2\|g_k\|}{M}\sqrt{\epsilon} - \frac{(U_H +\gamma_k)\epsilon}{M^2} \\
                       & \ge -2\sqrt{\epsilon} - \frac{2 U_g}{M}\sqrt{\epsilon} - \frac{(2U_H + U_g)\epsilon}{M^2}       \\
                       & \ge -2\left(1 + \frac{2U_g}{M}\right)\sqrt{\epsilon},
    \end{align*}
    where the last inequality holds since $\epsilon \le (2M U_g / (2U_H + U_g) )^2$. For the case of $\|r_k\| \le 2\epsilon$, since $\|d_k\| = \|v_k/t_k\| \le \Delta$, using the similar argument in \cref{lemma.norm of gk+1} gives $\lambda_1(H_{k+1}) \ge \Omega(-\sqrt{\epsilon})$. Now we inspect the value of $\|g_{k+1}\|$. By the second-order Lipschitz continuity, we have
    \begin{equation}
        \begin{aligned}
            \|\gkn\| & \leq \|\gkn-\gk -H_k d_k\|+\| \gk+ H_k d_k\|             \\
                     & \leq \frac{M}{2}\|\dk\|^2 +\|\gk+H_k d_k\|               \\
                     & =\frac{M}{2}\|\dk\|^2+\|r_k / t_k -\gamma_k d_k\|        \\
                     & \le \frac{M}{2}\Delta^2 + \nu\|r_k\| + |\gamma_k|\Delta.
        \end{aligned}
        \notag
    \end{equation}
    where the equality holds due to \cref{eq.inexact.opt}. Recall that $|\sigma_k| \le \epsilon^2 / 16M^2$ holds with a probability of at least ${1-4p}$ and $\gamma_k = \delta + \sigma_k/t_k - g_k^T d_k$, we have
    \begin{align*}
        |\gamma_k| & \le  |\delta| + \left|\frac{\sigma_k}{t_k}\right| + |g_k^T d_k| \\
                   & \le \sqrt{\epsilon} + \frac{\epsilon^2}{8M^2}+U_g \Delta ,
    \end{align*}
    where the second inequality holds since $|t_k| > \sqrt{1/(1+\Delta^2)} = M /\sqrt{\epsilon + M^2} \ge 1/2$ for any $\epsilon \le 3M^2$.
    Combining the above results with $\nu \le 1/2$ and $\|r_k\| \le 2\epsilon$, for any $0<\epsilon < 1$, we have that
    \begin{align*}
        \|\gkn\| \le \frac{\epsilon}{2M} + \epsilon + \frac{\epsilon}{M} + \frac{ \epsilon^{\frac{5}{2}}}{8M^3} + \frac{U_g \epsilon}{M} \le \left(\frac{5}{2M} + \frac{U_g}{M} + \frac{1}{8M^3} + 1\right)\epsilon,
    \end{align*}
    which means that $x_{k+1}$ is already an $\epsilon$-approximate SOSP. For the last statement, note that the homogenized matrix admits the form
    \begin{equation*}
        F_k =  \begin{bmatrix}
            H_k   & g_k     \\
            g_k^T & -\delta
        \end{bmatrix}.
    \end{equation*}
    {
    In view of \eqref{eq.smalleigenlb}, as initially $\delta := \sqrt\epsilon$, we have a low bound on $\lambda_1(H_k)$,
    \begin{equation}\label{eq.eigenlb1}
        \lambda_1(H_k) +  2\sqrt\epsilon + 2\|\gk\|\Delta + (U_H + \gamma_k)\Delta^2 \ge 0.
    \end{equation}
    Since we reset $\delta := 3\sqrt\epsilon + 2\|\gk\|\Delta + (U_H + \gamma_k)\Delta^2$, the Cauchy interlace theorem gives that
    \begin{align*}
        \varsigma_k & = \lambda_2(F_k) - \lambda_1(F_k) \ge \lambda_1(H_k) + \delta \stackrel{\eqref{eq.eigenlb1}}{\ge} \sqrt{\epsilon},
    \end{align*}
    which completes the proof.
    }
\end{proof}

It remains to characterize the scenario in which the increased perturbation is used (\cref{line.increase.pert}). For the newly calculated Ritz pair $[v_k;t_k]$, if it falls into the large value case, the function value decreases, and we proceed to the next iteration. The key aspect is that if $[v_k;t_k]$ falls again into the small value case, then $\|r_k\| \le 2\epsilon$ must hold, indicating that $x_{k+1} = x_k + d_k$ is an $\epsilon$-approximate SOSP. This argument is formalized as follows.
\begin{lemma}\label{lemma.inexact small step perturb}
    Suppose that \cref{assm.lipschitz} and \cref{assum.bounded gradient} hold, and we reset
    \begin{equation*}
        \delta = 3\sqrt\epsilon + 2\|\gk\|\Delta + (U_H + \gamma_k)\Delta^2,\quad e_k = \min\left\{\epsilon,\frac{\epsilon^{\frac{5}{2}}}{4(U_H + U_g)^2}\right\}, \quad \mbox{and}\quad \Delta = \frac{\sqrt{\epsilon}}{M}
    \end{equation*}
    in \cref{line.increase.pert} of \cref{alg inexact HSODM}. For any $0 < \epsilon < 1$ and $p \in (\exp(-n),1)$, if $|t_k| > \sqrt{1/(1+\Delta^2)}$, then $\|r_k\| \le 2\epsilon$ holds with a probability of at least ${1-4p}$.
\end{lemma}
\begin{proof}
    Note that for the increased $\delta$, by \cref{lemma.inexact small step} it holds that $\varsigma_k = \lambda_2(F_k) - \lambda_1(F_k) \ge \sqrt{\epsilon}$. From \cref{lem.convergence r sigma}, we have
    \begin{align*}
        \|r_k\| & \leq \tau_k e_k +2(\max\{U_H, \delta\} + \|\gk\|) \sqrt{\frac{e_k}{{\varsigma_k}}} \\
                & \leq \tau_k e_k + 2(U_H + U_g) \sqrt{\frac{e_k}{{\sqrt{\epsilon}}}}                \\
                & \leq 2\epsilon,
    \end{align*}
    where the last inequality holds because of $\tau_k < 1$ (cf. \Cref{lem.convergence r sigma}, \eqref{eq.bound r sigma in terms of e}).
\end{proof} \endproof

\subsection{Global convergence analysis of the inexact HSODM}
\label{sec.inexact.complexity}
Finally, we are ready to analyze the complexity of the Lanczos method.

\begin{corollary}[Complexity of \cref{alg.randlanczos}]\label{cor.complexity lanczos}
    Suppose that \cref{assm.lipschitz} and \cref{assum.bounded gradient} hold. When \cref{alg.randlanczos} is called in \cref{line.solve} in the inexact HSODM, for any constant $p \in (\exp(-n),1)$, with a probability of at least ${1-4p}$, its number of iterations to complete one call is upper bounded by
    \begin{equation*}
        O\left(\sqrt{\|F_k\|}\epsilon^{-1/4}\log\left(\frac{n(n+1)}{p\epsilon}\right)\right).
    \end{equation*}
\end{corollary}
\begin{proof}
    Recall that \cref{thm.customized initialization} shows that the inner product $q_1^T \chi_k > 0$ with a probability of at least ${1-4p}$, which facilitates the application of the complexity result in \cref{lemma.complexity of lanczos method}. Specifically, we know $(q_1^T\chi_k)^2$ is bounded away from $0$, and it generally attains the first term in \eqref{eq.psi.bound.ip}, which is in the order of $\Omega(\epsilon^4/n(n+1))$, as the second term in \eqref{eq.psi.bound.ip} is almost constant (like the last term) as $\epsilon < 1$ is small.

    Note that only two cases may occur when \cref{alg.randlanczos} is called in the inexact HSODM at some iteration $k$. In the first case, we set $e_k = \sqrt\epsilon$. By \eqref{eq.gap-free}, the worst-case complexity is thus
    $
        O\left(\sqrt{\|F_k\|}\epsilon^{-1/4}\log(n(n+1)/(p\epsilon))\right).
    $
    In the second case, we set $\delta$ to a larger value (\cref{line.increase.pert}), and by \cref{lemma.inexact small step}, we know $\varsigma_k =
        \lambda_2(F_k) - \lambda_1(F_k) \ge \sqrt{\epsilon}.
    $
    Hence, we are safe to use a higher accuracy while keeping the complexity in the same order by the gap-dependent estimate \eqref{eq.gap-dependent}.  \end{proof} \endproof

In summary, we show that in any case, the Lanczos method in \cref{alg inexact HSODM} is guaranteed to terminate in $\tilde O(\epsilon^{-1/4})$ iterations. However, contrasting with the complexity result presented in \cite{royer_complexity_2018,royer_newton-cg_2020}, which depends on $\|H_k\|$ and can be capped by $U_H$, our approach necessitates the magnitude of $\|F_k\|$, which is upper bounded by $U_H + U_g$.
In the following theorem, we prove the arithmetic complexity of inexact HSODM.
\begin{theorem}[Complexity of the inexact HSODM]
    \label{thm.inexact complexity}
    Suppose that \cref{assm.lipschitz} and \cref{assum.bounded gradient} hold.  For any constant $p \in (\exp(-n),1)$, the inexact HSODM (\cref{alg inexact HSODM}) terminates in
    $$ K= 12(f(x_1) - f_{\inf})M^2\epsilon^{-3/2}$$
    iterations and returns an iterate $\xkn$ such that
    \begin{equation*}
        \|g_{k+1}\|\leq O(\epsilon) \quad \text{and} \quad  \lambda_1 (H_{k+1}) \geq \Omega(-\sqrt{\epsilon})
    \end{equation*}
    with a probability of at least $({1-4p})^{2K}$. Furthermore, the arithmetic operations required by \cref{alg inexact HSODM} are bounded from above by
    $$O\left((n+1)^2 \epsilon^{-7/4} (f(x_1) - f_{\inf})M^2\sqrt{U_H + U_g} \log(n(n+1)/(p\epsilon))\right).$$
\end{theorem}

\begin{proof}

    For the two large cases in \cref{alg inexact HSODM}, \cref{lemma.inexact function value decrease} implies that the function value decreases at least
    \begin{equation*}
        f(x_{k+1}) - f(x_k) \le -\frac{\delta}{4}\Delta^2 + \frac{M}{6}\Delta^3 = -\frac{\epsilon^{3/2}}{12M^2}
    \end{equation*}
    by selecting $\delta=\sqrt{\epsilon}$ and $\Delta = \frac{\sqrt{\epsilon}}{M}$.
    While according to \cref{lemma.inexact small step} and \cref{lemma.inexact small step perturb}, in the small value case, the algorithm will terminate at an $\epsilon$-approximate SOSP or come back to the large value case.
    Consequently, we obtain that the number of iterations is bounded above by $K= 12(f(x_1) - f_{\inf})M^2\epsilon^{-3/2}$ before reaching an $\epsilon$-approximate SOSP. At each iteration, one inquiry of \cref{alg.randlanczos} is needed if we have the large value case. Otherwise, we have to reset the parameters (cf. \cref{line.increase.pert}). In that case, we either fall into the large value case and proceed, or again into the small value case. The latter implies that $\|r_k\| \le 2\epsilon$ as shown in \cref{lemma.inexact small step perturb} and will terminate the algorithm. To sum up, each iteration needs at most 2 inquiries of \cref{alg.randlanczos} in high probability.
    Since there is a probability of at least ${1-4p}$ that the Lanczos method will succeed, we have no incorrect termination of \cref{alg.randlanczos} occurs in the $K$ iterations with a probability of at least $({1-4p})^{2K}$. Combining these results with \cref{cor.complexity lanczos}, the complexity of the arithmetic operations can be established.

\end{proof} \endproof

{We remark that $(1-4p)^{2K} \ge 1 - 8Kp$ holds for some $p$ satisfying $p < 1/2K$. Recall $p \in (\exp(-n),1)$, this condition can be easily met when $n \ge \Omega(-\log\epsilon)$. For example, setting $\epsilon = 10^{-8}$ yields $n \approx 20$. Therefore, ``with a probability of at least $(1-4p)^{2K}$" in the theorem can be replaced by ``with a probability of at least $1 - 8Kp$" while remaining informative.} Since our algorithm requires arithmetic operations on a homogenized matrix of dimension $(n+1)$, its dependency on dimension and the complexity associated with eigenvalue procedure (\cref{cor.complexity lanczos}) are slightly worse compared to prior second-order algorithms, such as \cite{royer_complexity_2018,royer_newton-cg_2020,curtis2021newton-cg,carmon_accelerated_2018,agarwal_finding_2017}. Regarding the Lipschitz constants, the dependency on Hessian Lipschitz constant $M$ in our bound is comparatively inferior to those in \cite{agarwal_finding_2017,carmon_accelerated_2018}, as our algorithm does not explicitly incorporate this constant; rather, it is only invoked in establishing the overall computational complexity. Nevertheless, our algorithm, HSODM, is characterized by its conciseness and unity, requiring only the eigenvalue procedure at each iteration. Specifically, it achieves computational efficiency superior to Newton-type methods when encountering degeneracy in the Hessian matrix \cite{he_homogeneous_2023}. Furthermore, the subsequent section also demonstrates the promising practical performance of HSODM.

\section{Numerical Experiments}
\label{sec.experiment}
In this section, we provide the computational results of HSODM on a few classes of nonconvex optimization problems.
% To be specific, we include a set of nonconvex $L_2-L_p$ minimization problem that arises from compressed sensing and has long been one of great interests in the community. In addition,
We include the CUTEst problems \cite{gould_cutest_2015} since they serve as a standard dataset to test the performance of algorithms for nonlinear problems. Because the HSODM belongs to the family of second-order methods, we focus on comparisons with Newton trust-region method and adaptive cubic regularized Newton method \cite{cartis_adaptive_2011}. Our implementation in Julia \cite{Julia-2017} is provided at \url{https://github.com/bzhangcw/DRSOM.jl}. All experiments are conducted in Julia, and the development is handled by a desktop of MacOS with a 3.2 GHz 6-Core Intel Core i7 processor.

\subsection{Implementation details}
Apart from the original form of HSODM (see \cref{alg.main alg}), we add a few techniques for practical implementations.
We first note that a practical HSODM may not explicitly use the Hessian matrix $\Hk$. In the computation of $F_k \cdot [v; t]$ where $v\in\real^n, t\in \real$, we have
$$
    F_k \cdot \begin{bmatrix}v \\ t\end{bmatrix} = \begin{bmatrix}
        \Hk \cdot v + t\cdot \gk \\
        \gk^T v - t\cdot \delta
    \end{bmatrix}.
$$
From the above fact, a matrix-free option by utilizing the Hessian-vector product $\Hk v$ is provided as in other inexact Newton-type methods \cite{cartis_adaptive_2011,curtis_inexact_2018}.

Not limited to the backtrack line-search algorithm for theoretical analysis, in practice, the homogeneous direction should work with any well-defined line-search method. In our implementation, we apply the Hager-Zhang line-search method with default parameter settings \cite{hager_algorithm_2006}. For eigenvalue problems, we use the Lanczos method to solve homogenized subproblems with a given tolerance, $10^{-6}$. Since these methods are readily provided by a few efficient Julia packages, we directly use the line-search algorithms from LineSearches.jl \cite{kmogensenOptimMathematicalOptimization2018}, and the Lanczos method from KrylovKit.jl \cite{haegeman_krylovkit_2024}. For hyperparameters, we set $\delta = -\sqrt\epsilon$, $\nu = 0.01$, and $\Delta = 10^{-4}$.

\paragraph{The benchmark algorithms}

\citet{orban_juliasmoothoptimizers_2019} provided highly efficient Julia packages in the JuliaSmoothOptimizers organization that include the Newton trust-region method utilizing the Steihaug-Toint conjugate-gradient method (\newtontrst) and an adaptive cubic regularization (\arc{}) with necessary subroutines and techniques including subproblem solutions and Krylov methods.  The numerical results are recently reported in \cite{dussaultScalableAdaptiveCubic2023}.  We use the original implementation in \cite{orban_juliasmoothoptimizers_2019} and the default settings therein.
% \input{l2lp.tex}
%%%%%%%%%%%%%%%%%
% cutest
\subsection{Unconstrained problems in CUTEst}

We next present the results on a selected subset of the CUTEst dataset. To set a comprehensive comparison, we provide the results of HSODM with readily Hessian matrices, named after \hsodm{}, and a version facilitated by Hessian-vector products (\hsodmhvp{}).
We set an iteration limit of $20,000$ and termination criterion as $\|\nabla f(\xk)\| \le 10^{-5}$ for all the tested algorithms; we check if this criterion is ensured else marked as failed.
We focus on the unconstrained problems with the number of variables $n \in [4, 5000]$.
% ; we regard a problem as small if the problem dimension $n \le 200$
For each problem in the CUTEst, if it has different parameters, we select all instances that fit the criterion.
% we choose the smallest instance that fits the criterion. 
Then we have 200 instances in total where a few instances cannot be solved by any method. The complete result can be found in \cref{tab.cutest.kt} and \cref{tab.cutest.fx}.

\paragraph{Overall comparison of the algorithms.}
The following \cref{tab.perf.geocutest} presents a summary of tested algorithms. In this table, we let $\mathcal K$ be the number of successful instances. Besides, we compute performance statistics based on scaled geometric means (\textsf{SGM}), including $\overline t_G,\overline k_G,\overline k_G^f,\overline k_G^g,\overline k_G^H$ as (geometric) mean running time, mean iteration number, mean function evaluations, mean gradient evaluations, and mean Hessian evaluations, respectively. The running time is scaled by 1 second, and other metrics are scaled by 50 evaluations or iterations accordingly.
% We let $\overline t,\overline k,\overline k^f,\overline k^g,\overline k^H$ be the mean running time, mean iteration number, mean function evaluations, mean gradient evaluations, and mean Hessian evaluations, respectively. These average numbers are taken among successful instances. 
Note that the cubic regularization \arc{}, \newtontrst{}, and \hsodmhvp{} use Hessian-vector products, so that $\overline k_G^H = 0$ and the gradient evaluations in $\overline k_G^g$ actually include the number of Hessian-vector products.

% Moreover, we also calculate the scaled geometric means that are marked with a subscript $G$, 
\begin{table}[h]
    \centering
    \caption{Performance in \textsf{SGM} of different algorithms on the CUTEst dataset. Note $\overline t_{G}, \overline k_{G}$ are scaled geometric means (scaled by $1$ second and $50$ iterations, respectively). If an instance is failed, its iteration number and solving time are set to $20,000$.}
    \label{tab.perf.geocutest}
    \begin{tabular}{lcccccc}
        \toprule
        Method      & $\mathcal K$ & $\overline t_G$ & $\overline k_G$ & $\overline k_G^f$ & $\overline k_G^g$ & $\overline k_G^H$ \\
        \midrule
        % \newtontrst         & 155.00       & 15.41           & 216.59          & 211.99            & 219.58            & 203.82            \\
        \newtontrst & 165.00       & 6.14            & 170.44          & 170.44            & 639.64            & 0.00              \\
        \arc        & 167.00       & 5.32            & 185.03          & 185.03            & 888.35            & 0.00              \\
        \hsodmhvp   & 173.00       & 4.79            & 111.24          & 200.60            & 787.32            & 0.00              \\
        \hsodm      & 174.00       & 4.86            & 113.30          & 197.46            & 256.20            & 111.28            \\
        \bottomrule
    \end{tabular}

\end{table}

\begin{figure}[h]
    \centering
    \subfloat[Performance of iteration number]{
        \includegraphics[width=0.45\linewidth]{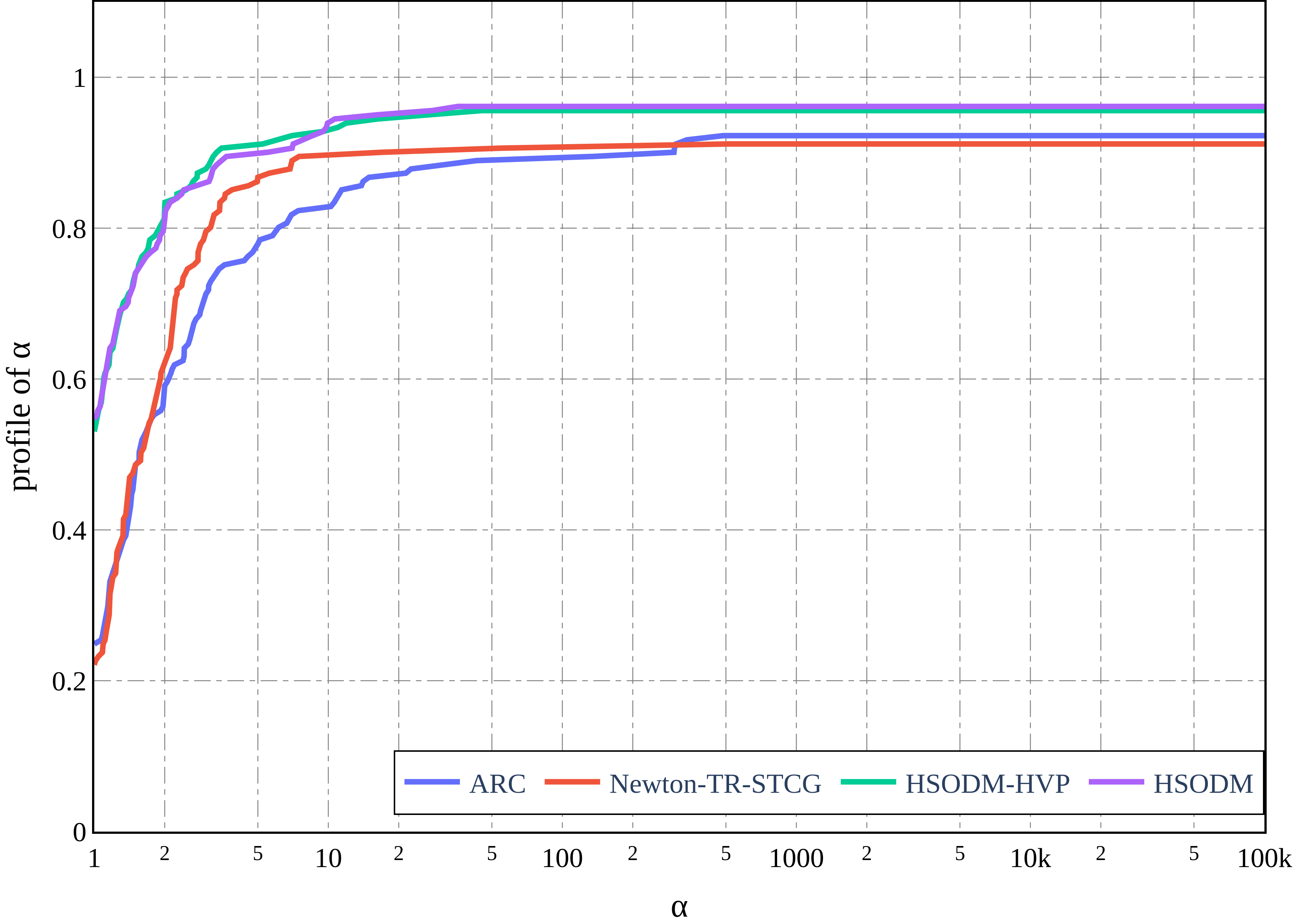}
    }
    ~
    \subfloat[Performance of gradient evaluations]{
        \includegraphics[width=0.45\linewidth]{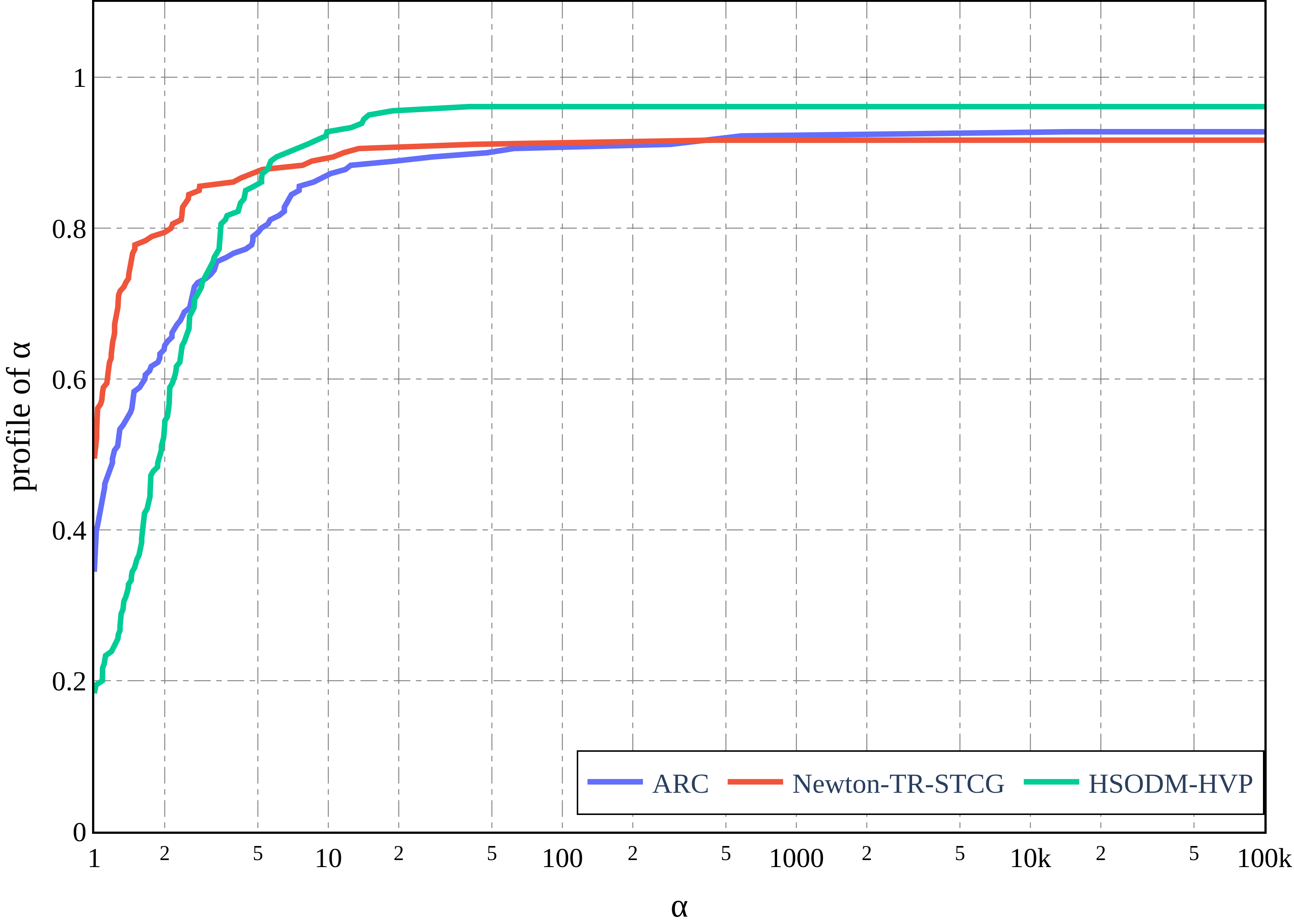}
    }\label{fig.perf.ig}
    \caption{Performance profiles of the second-order methods for CUTEst problems. In (a), we report the iteration number. Figure (b) includes the results of gradient evaluations; we only include methods using Krylov subspaces.}
    \label{fig.perfprof}
\end{figure}

Apart from metrics measured by \textsf{SGM}, we use the performance profile on iteration number as defined in \cite{dolan_benchmarking_2002}.
In essence, the performance profile at point $\alpha$ in \cref{fig.perfprof} of an algorithm indicates the probability of successfully solved instances within $2^\alpha$ times the best iteration number amongst competitors.

\normalsize

The results from these preliminary implementations show that \hsodm{} and \hsodmhvp{} outperformed the standard second-order methods, including \newtontrst{} and \arc{} on average. \hsodmhvp{} and \hsodm{} had better iteration complexity and running time in terms of $\overline k_G, \overline t_G$ among competing algorithms. The HVP variant \hsodmhvp{} used comparable gradient evaluations with \arc{}.
Since HSODM needs fewer iterations, more gradient evaluations seem necessary. More function evaluations are needed by extra overhead from the line searches.
It is also interesting to see \hsodm{} and also \hsodmhvp{} (see \texttt{EXTROSNB}), \newtontrst{} (see \texttt{ARGLINC}) and \arc{}  (see \texttt{OSCIGRAD}) all had instances on which they performed best.

In terms of performance profile, we see both \hsodm{} and \hsodmhvp{} had an advantage in iteration numbers. \newtontrst{} has the best performance on gradient evaluations in its succeeded instances. \hsodmhvp{} needs more gradient evaluations since it uses a slightly larger $n+1$ dimensional system. Nevertheless, this disadvantage seems to be mild in practice.

% However, \hsodm{} still permits an advantage over \lbfgs{} if one seeks to find solutions with high precision (for example, $\|\gk\| \le 10^{-8}$).
%Whether this is a consequence of the global rate is not clear. 
%Since the current implementations for \hsodm{} only apply the standard trust-region update rules, i.e., the ratio-based acceptance and linear radius adjustments, we
%are confident that using the provably better frameworks like those in \cite{cartis_adaptive_2011,curtis_trust_2017} would have a promising improvement to \hsodm{}.

\section{Conclusion}\label{sec.conclusion}

In this paper, we introduce a homogenized second-order descent method (HSODM) whose global rate of complexity is optimal among a certain broad class of second-order methods (see \cite{cartis_evaluation_2022}). The HSODM utilizes the homogenization trick to the quadratic model, which comes from the standard second-order Taylor expansion, such that the resulting homogenized quadratic form can be solved as an eigenvalue problem. We have shown that the homogenized idea is well-defined in both convex and nonconvex cases, where a negative curvature direction always exists. Using the model all along, one can safely stop at a small step to obtain an $\epsilon$-approximate SOSP without switching to other methods.

We provide comprehensive experiments of HSODM on nonlinear optimization problems in the CUTEst benchmark. Two variants of HSODM show promising results in these experiments. One future direction is to utilize the method for constrained optimization problems.

\clearpage
\addcontentsline{toc}{section}{References}
\bibliography{homo}
\bibliographystyle{plainnat}
\clearpage

\appendix
\section{Appendix}
\input{appendix_proof}

\input{appendix_cutest}

\end{document}

%% file: header.tex
\usepackage{amssymb,amsmath,amsthm,enumerate}
\usepackage{array}
\usepackage{mathtools}
\usepackage[parfill]{parskip}
\usepackage{graphicx}
\usepackage{geometry}
\usepackage{mathrsfs}
\usepackage{caption}
\usepackage{subcaption}
\usepackage{bm,amsfonts,amscd}
\usepackage[]{units}
\usepackage{listings}
\usepackage{multicol,multirow}
\usepackage{longtable}
\usepackage{lscape}
\usepackage{tcolorbox,csquotes}
\usepackage{float}
\usepackage{tikz}
\usepackage[ruled, linesnumbered]{algorithm2e}
\usepackage{authblk}
\usepackage{booktabs}
\usepackage[colorlinks]{hyperref}
\usepackage[numbers]{natbib}

\newtheorem{theorem}{Theorem}

\newtheorem{corollary}{Corollary}
\newtheorem{lemma}{Lemma}
\newtheorem{assumption}{Assumption}

\newtheorem{remark}{Remark}

\numberwithin{assumption}{section}
\numberwithin{definition}{section}
\numberwithin{equation}{section}
\numberwithin{lemma}{section}
\numberwithin{corollary}{section}
\numberwithin{theorem}{section}
\numberwithin{figure}{section}
\numberwithin{table}{section}

\usepackage[nameinlink,capitalize]{cleveref}
\crefname{assumption}{Assumption}{Assumptions}
\crefname{Assumption}{Assumption}{Assumptions}
\crefname{Lemma}{Lemma}{Lemmata}
\crefname{lemma}{Lemma}{Lemmata}
\crefname{Theorem}{Theorem}{Theorems}
\crefname{theorem}{Theorem}{Theorems}
\crefname{corollary}{Corollary}{Corollaries}
\crefname{proposition}{Proposition}{Propositions}
\crefname{claim}{Claim}{Claims}
\crefname{procedure}{Procedure}{Procedures}
\crefname{algorithm}{Algorithm}{Algorithms}
\crefname{figure}{Figure}{Figures}
\crefname{remark}{Remark}{Remarks}
\crefname{section}{Section}{Sections}
\crefname{procedure}{Procedure}{Procedures}
\crefname{definition}{Definition}{Definitions}
\crefname{example}{Example}{Examples}
\crefname{table}{Table}{Tables}
\crefname{line}{Line}{Lines}
\crefname{equation}{}{}
\crefname{enumi}{}{}

\author[1]{\small Chuwen Zhang\thanks{This research is partially supported by the National Natural Science Foundation of China (NSFC) [Grant NSFC-72150001, 72225009, 11831002] and the Natural Science Foundation of Shanghai [23ZR1445900].}}
\author[2]{\small Dongdong Ge}
\author[1]{\small Chang He}
\author[1]{\small Yuntian Jiang}
\author[1]{\small Chenyu Xue}
\author[1]{\small Bo Jiang}
\author[3]{\small Yinyu Ye}

\affil[1]{\footnotesize School of Information Management and Engineering\\ Shanghai University of Finance and Economics}
\affil[2]{\footnotesize Antai College of Economics and Management\\ Shanghai Jiao Tong University}
\affil[3]{\footnotesize Department of Management Science and Engineering, Stanford University}

\title{A Homogeneous Second-Order Descent Method for Nonconvex Optimization}

%% file: macro.tex
\newcommand{\real}{\mathbb{R}}

\newcommand{\argmin}{\arg\min}

\newcommand{\xk}{x_k}
\newcommand{\Hk}{H_k}

\newcommand{\gk}{g_k}
\newcommand{\dk}{d_k}
\newcommand{\gkn}{g_{k+1}}
\newcommand{\xkn}{x_{k+1}}

\newcommand{\hsodm}{\textsf{HSODM}}
\newcommand{\hsodmhvp}{\textsf{HSODM-HVP}}

\newcommand{\newtontrst}{\textsf{Newton-TR-STCG}}

\newcommand{\arc}{\textsf{ARC}}

\newcommand{\rank}{\text{rank}}
\newcommand{\cK}{\mathcal{K}}

\definecolor{rosemary}{RGB}{190,88,69}
\definecolor{stone}{RGB}{42,121,10}
\definecolor{deepblack}{RGB}{90,90,90}

\hypersetup{
    colorlinks,
    linkcolor={deepblack},
    citecolor={rosemary},
    urlcolor={rosemary}
}

%% file: appendix_proof.tex
\section{Additional Proofs}\label{sec.addproofs}
\subsection{Proof of \cref{lemma.complexity of lanczos method}} We provide a sketch here as the results are the combination of the complexity estimates in \cite{jin_lecture_2021} and Lemma 9 in \cite{royer_complexity_2018}. Consider the positive semidefinite matrix $F_k^\prime: = \|F_k\|I - F_k$, and substituting $\epsilon:= \frac{e_k}{2\|F_k\|}$ into the complexity results in \cite{jin_lecture_2021}, the Lanczos method returns an estimate $\gamma_{\max}(F_k^\prime)$ satisfies
\begin{align*}
    \gamma_{\max}(F_k^\prime) \ge \left(1 - \frac{e_k}{2\|F_k\|}\right)\lambda_{\max}(F_k^\prime)
\end{align*}
if it starts with the vector $q_1$ and runs at most
\begin{equation*}
    1+2\sqrt{\frac{\|F_k\|}{e_k}}\log\left(\frac{16\|F_k\|}{e_k(q_1^T \chi_k)^2}\right)
\end{equation*}
iterations (gap-free version). Since $\gamma_{\max}(F_k^\prime) = \|F_k\| - \gamma_k$ and $\lambda_{\max}(F_k^\prime) = \|F_k\| - \lambda_1(F_k)$, following the same argument of Lemma 9 in \cite{royer_complexity_2018}, we obtain
\begin{equation*}
    \gamma_k \le \lambda_1(F_k) + e_k.
\end{equation*}
The result of the gap-dependent version can be established similarly, and thus we omit it here.
$\qed$

\subsection{Proof of \cref{lem.convergence r sigma}}
For the first statement, note that
\begin{align*}
    \|F_k\|   & = \max_{\|[v;t]\| = 1} \begin{bmatrix}
                                           v \\ t
                                       \end{bmatrix}^T \begin{bmatrix}
                                                           H_k   & g_k     \\
                                                           g_k^T & -\delta
                                                       \end{bmatrix} \begin{bmatrix}
                                                                         v \\ t
                                                                     \end{bmatrix}                  \\
    \nonumber & \le \max_{\|[v;t]\| = 1} \begin{bmatrix}
                                             v \\ t
                                         \end{bmatrix}^T \begin{bmatrix}
                                                             H_k & 0       \\
                                                             0   & -\delta
                                                         \end{bmatrix} \begin{bmatrix}
                                                                           v \\ t
                                                                       \end{bmatrix} + \max_{\|[v;t]
    \| = 1} \begin{bmatrix}
                v \\ t
            \end{bmatrix}^T \begin{bmatrix}
                                0     & g_k \\
                                g_k^T & 0
                            \end{bmatrix} \begin{bmatrix}
                                              v \\ t
                                          \end{bmatrix}                                             \\
              & \le \max\{U_H, \delta\} + \|\gk\|,
\end{align*}
which completes the proof. For the second argument, multiplying $[v_k;t_k]$ on both sides of \eqref{eq.inexact optimal condition} yields
\begin{equation}\label{eq.regularity d}
    [v_k;t_k]^T F_k [v_k;t_k] +\gamma_k =0.
\end{equation}
Since $[v_k;t_k]$ is a unit vector, we can rewrite $[v_k;t_k] = \tau_k \cdot \chi_k+s$ for some $\tau_k \in [0,1]$and $s \perp \chi_k$ satisfying $\tau_k^2+\|s\|^2=1$. Substituting into \eqref{eq.regularity d} gives
\begin{align}
    \nonumber  -\theta_k+e_k & \geq -\gamma_k  =-\theta_k \tau_k^2+s^T F_k s        \notag      \\
                             & \geq -\theta_k \tau_k^2+(-\theta_k+{\varsigma_k})\|s\|^2, \notag
\end{align}
where the equality is obtained by the fact $s\perp \chi_k$. It implies
\begin{equation}
    \label{eq.bound s final}
    \|s\|^2\leq \frac{e_k}{{\varsigma_k}}.
\end{equation}
Thus from \eqref{eq.inexact optimal condition} we have
\begin{equation}
    \label{eq.bound r}
    \begin{aligned}
        [r_k;\sigma_k] & = F_k [v_k;t_k]+\gamma_k [v_k;t_k]                     \\
                       & = (F_k+\gamma_k I)(\tau_k \chi_k +s)                   \\
                       & = \tau_k (\gamma_k-\theta_k)\chi_k +(F_k+\gamma_k I)s.
    \end{aligned}
\end{equation}
Hence, the norm of the residual follows
\begin{equation}
    \label{eq.bound r final}
    \begin{aligned}
        \|r_k\| & \le \|[r_k;\sigma_k]\|                                                              \\ & \leq \tau_k(\theta_k-\gamma_k)+\|(F_k+\gamma_k I)s\|                             \\
                & \leq \tau_k e_k + \|(F_k+\gamma_k I)\|\sqrt{\frac{e_k}{{\varsigma_k}}}              \\
                & \leq \tau_k e_k +2(\max\{U_H, \delta\} + \|\gk\|) \sqrt{\frac{e_k}{{\varsigma_k}}}.
    \end{aligned}
\end{equation}
This completes the proof.
$\qed$

\subsection{Proof of \cref{thm.gamma ge delta}}\label{sec.proof.gammagedelta}
For part (1), due to the mechanism of the Lanczos method, for any orthonormal basis $q_j = [\ell_j; \beta_j]$ with $j\ge 2$, we have $q_j \perp q_1$. Therefore, it holds that
\begin{align*}
    \beta_j \alpha & = - \ell_j^Tu \sqrt{1-\alpha^2},
\end{align*}
and it implies
\begin{equation}\label{eq.bound.lastentry}
    |\beta_j| \le \frac{\sqrt{1-\alpha^2}\|\ell_j\|\|u\|}{|\alpha|} \le 2\sqrt{1-\alpha^2}
\end{equation}
for any $|\alpha| \ge 1/2$.

For part (2), denote $\zeta_{n+1} = F_k 1_{n+1}$ and $y = \zeta_{n+1} - (\zeta_{n+1}^Tq_1) \cdot q_1 - (\zeta_{n+1}^Tq_2) \cdot q_2$, then $q_1,q_2,y$ are mutually orthogonal. Therefore, let $\Pi$ be the projection matrix onto the subspace spanned by $q_1, q_2$, then $y$ is the residual of $\zeta_{n+1}$ after projecting on this subspace, and it follows
\begin{align*}
    \|y\| = \|(I_{n+1} - \Pi) \zeta_{n+1}\|.
\end{align*}
By denoting $\varphi := \zeta_{n+1} - F_kq_1$, we conclude,
\begin{align*}
    \|y\| = \|(I_{n+1} - \Pi) \zeta_{n+1}\| = \|(I_{n+1} - \Pi) (F_kq_1 + \varphi) \| = \|(I_{n+1} - \Pi)\varphi \| \le \|\varphi\| \\= \left\|\begin{bmatrix}
        (1-\alpha) \cdot g_k -\sqrt{1-\alpha^2} \cdot H_ku \\
        (1-\alpha) \cdot (-\delta) - \sqrt{1-\alpha^2} \cdot g_k^Tu
    \end{bmatrix}\right\|
\end{align*}
since $F_kq_1 \in \mathcal K(2;F_k,q_1)$ and $\|I_{n+1} - \Pi\| = 1$.
    {In this view, we have,
        \begin{align}\label{eq.residual.by2projection}
            \|y\|^2 \le ((1-\alpha^2)(\|\Hk\| + \|\gk\|)^2 + (1-\alpha)^2(\delta+\|\gk\|)^2) \le ((U_H + U_g)^2 + (\delta+U_g)^2) \cdot (1-\alpha^2)
        \end{align}
    }
as $1-\alpha \le \sqrt{1-\alpha^2}$ holds for $\alpha \in (0,1)$. Recall that for the Lanczos method, it holds that
\begin{equation*}
    F_k Q_j - Q_jT_j = \xi_j 1_j^T,~Q_j = [q_1, \ldots, q_j] \in \mathbb{R}^{(n+1)\times j},~T_j \in \mathbb{R}^{j \times j}.
\end{equation*}
Consider the last term of the residual $\xi_j$, by $\xi_{j,n+1}$, for $j \ge 3$, it follows
\begin{align*}
    \xi_{j,n+1} & = 1_{n+1}^T\xi_j = 1_{n+1}^T\xi_j1_j^T1_j                                      \\
                & = 1_{n+1}^TF_kQ_j1_j - 1_{n+1}^TQ_jT_j1_j                                      \\
                & = \zeta_{n+1}^Tq_j - [\beta_1, ..., \beta_{j}] [0,...,0, T_{j-1,j}, T_{j,j}]^T \\
                & = \zeta_{n+1}^Tq_j - \beta_{j-1}T_{j-1,j} - \beta_{j} T_{j,j}
\end{align*}
Since $q_j$ is perpendicular to $q_1$ and $q_2$, we have $q_j^T\zeta_{n+1} = q_j^Ty$, and thus
\begin{subequations}
    \begin{align}
        \nonumber |\xi_{j,n+1}| & = |\zeta_{n+1}^Ty - \beta_{j-1}T_{j-1,j} - \beta_{j} T_{j,j}|                                                \\
        \nonumber               & \le \|\zeta_{n+1}\| \cdot \|y\| +  |T_{j-1,j}| \cdot |\beta_{j-1}| +  |T_{j,j}| \cdot |\beta_j|              \\
        \label{eq.bd.1}         & \le \sqrt{1-\alpha^2}\left(\sqrt{(U_H + U_g)^2 + (\delta+U_g)^2}\sqrt{U_g^2+\delta^2} + 4\|T\|_\infty\right) \\
        \label{eq.bd.2}         & \le \sqrt{1-\alpha^2}U_\sigma
        = O(\sqrt{1-\alpha^2})
    \end{align}
\end{subequations}
where \eqref{eq.bd.1} follows from \eqref{eq.bound.lastentry} and \eqref{eq.residual.by2projection}.
The last inequality \eqref{eq.bd.2} follows from the fact that $\|T\|_\infty \le \sqrt{n}\|T\|\le \|F_k\|$ since the spectra of $T$ is bounded by that of $F_k$ (see, e.g., \cite[Theorem 10.1.2]{golubMatrixComputations2013}). By taking $U_\sigma := \sqrt{(U_H + U_g)^2 + (\delta+U_g)^2}\sqrt{U_g^2+\delta^2} + 4\sqrt{n}(U_g + \max\{U_H, \delta\})$, and by the fact of Ritz approximation (Section 10.1.4 in \cite{golubMatrixComputations2013}), we conclude
$$|\sigma_k| \le |\xi_{j,n+1}| \le \sqrt{1-\alpha^2}U_\sigma.$$

For part (3), from the shift-invariant property of the Krylov subspace, we have
\begin{equation*}
    \mathcal{K}(j;U_F I_{n+1} - F_k) = \mathcal{K}(j;F_k) := \left\{q_1, F_kq_1, \ldots, F_k^j q_1\right\}.
\end{equation*}
Since the Ritz value $-\gamma_k$ is generated in a larger Krylov subspace, $-\gamma_k \le q_1^T F_k q_1$ holds. Therefore, it is sufficient to establish that $q_1^T F_k q_1 \le -\delta$. The selection \eqref{eq.safeguard.gammadelta} implies
\begin{equation*}
    (U_H + \delta)^2 \cdot (1-\alpha^2) \le 4(\gk^Tu)^2 \cdot \alpha^2,
\end{equation*}
and thus, it follows
\begin{align*}
    \begin{bmatrix}
        \sqrt{1-\alpha^2} u \\
        \alpha
    \end{bmatrix}^T
    \begin{bmatrix}
        H_k   & g_k     \\
        g_k^T & -\delta
    \end{bmatrix}
    \begin{bmatrix}
        \sqrt{1-\alpha^2} u \\
        \alpha
    \end{bmatrix}
     & = -\delta \cdot \alpha^2 + 2\alpha\sqrt{1-\alpha^2} \cdot \gk^Tu + (1-\alpha^2) \cdot u^THu \\
     & \le -\delta \cdot \alpha^2 - (1 - \alpha^2) \cdot (U_H + \delta) + (1 - \alpha^2) \cdot U_H \\
     & \le -\delta.
\end{align*}

$\qed$

\subsection{Proof of \cref{thm.customized initialization}}
\label{sec.proof.customized}
We first provide a useful inequality. For any given constant $x \ge 0$, it holds that
$$
    1-\exp \left(-4 x^2 / \pi\right) \geq \operatorname{erf}(x)^2,
$$
where $\operatorname{erf}(x)=\frac{2}{\sqrt{\pi}} \int_0^x e^{-t^2} d t$. It implies that for any random variable $X \sim \mathcal{N}(0,1)$, we have
\begin{equation*}
    \operatorname{Prob}(|X| \le x) = \operatorname{erf}(x/\sqrt{2}) \le \sqrt{1-\exp \left(-2 x^2 / \pi\right)}.
\end{equation*}
Now we begin with establishing a lower bound of $q_1^T \chi_k$. Since
\begin{equation*}
    (q_1^T \chi_k)^2 = \frac{(\chi_k^T b)^2}{\|b\|^2} = \frac{(\sum_{i=1}^n \chi_{k, i} b_i + \Psi_k \cdot \chi_{k, n+1} b_{n+1})^2}{\sum_{i=1}^n b_i^2 + \Psi_k^2 b_{n+1}^2},
\end{equation*}
it is sufficient to provide a lower bound of $(\chi_k^T b)^2$ and an upper bound of $\|b\|^2$, respectively. For the term $(\chi_k^T b)^2$, recall $b_1, \ldots, b_{n+1} \overset{\textrm{i.i.d.}}{\sim} \mathcal{N}(0,1)$, then it holds that $\sum_{i=1}^n \chi_{k, i} b_i +  \Psi_k \cdot \chi_{k, n+1} b_{n+1} \sim \mathcal{N}(0, \sum_{i=1}^n \chi_{k, i}^2 + \Psi_k^2 \chi_{k, n+1}^2)$. Hence, we have
\begin{equation}\label{eq.bound gaussian}
    \begin{aligned}
              & \operatorname{Prob}\left(\left|\chi_k^T b\right| \le p\sqrt{\frac{\pi\left(\sum_{i=1}^n \chi_{k,i}^2 + \Psi_k^2 \chi_{k, n+1}^2\right)}{2}}\right)    \\
        = ~   & \operatorname{Prob}\left(\left|\frac{\chi_k^T b}{\sqrt{\sum_{i=1}^n \chi_{k,i}^2 + \Psi_k^2 \chi_{k, n+1}^2}}\right| \le p\sqrt{\frac{\pi}{2}}\right) \\
        \le ~ & \sqrt{1-\operatorname{exp}\left(-p^2\right)} \le p
    \end{aligned}
\end{equation}
for any constant $0 < p < 1$. Consequently, with a probability of at least $1 - p$, we conclude
\begin{equation*}
    (\chi_k^T b)^2 \ge \frac{p^2 \pi\left(\sum_{i=1}^n\chi_{k,i}^2 + \Psi_k^2 \chi_{k, n+1}^2\right)}{2}.
\end{equation*}
Then we consider the upper bound of $\|b\|^2$. Note that $\|b\|^2 \le \Psi_k^2 \cdot \sum_{i=1}^{n+1} b_{n+1}^2$, and $\sum_{i=1}^{n+1} b_{n+1}^2$ follows the chi-square distribution with $n+1$ degrees of freedom. Applying the tail bound (Lemma 1 in \cite{laurent2000adaptive}) gives that
\begin{equation*}
    \operatorname{Prob}\left(\sum_{i=1}^{n+1} b_{n+1}^2 \ge 5(n +1) \right) \le \exp(-(n+1)).
\end{equation*}
Hence, it holds that $\|b\|^2 \le 5\Psi_k^2(n + 1)$ with a probability of at least $1 - \exp(-(n+1))$. Therefore, {by applying the union bound, we conclude
        \begin{equation}\label{eq.lower of inner w.r.t. Psi}
            (q_1^T \chi_k)^2 = \frac{(\chi_k^T b)^2}{\|b\|^2} \ge \frac{\pi p^2 \sum_{i=1}^n\chi_{k,i}^2}{10\Psi_k^2(n + 1)} + \frac{\pi p^2 \chi_{k, n+1}^2}{10(n + 1)}
        \end{equation}
        with a probability of at least $1 - p - \exp(-(n+1))$.}

Now we justify the relationship between $\Psi_k$ and accuracy $\epsilon$. Motivated by \cref{thm.gamma ge delta}, we consider the following condition
\begin{equation}\label{eq.u1-a}
    1 - \alpha^2 = \frac{\sum_{i=1}^n b_i^2}{\sum_{i=1}^n b_i^2 + \Psi_k^2 b_{n+1}^2} \le \min\left\{\frac{\epsilon^4}{256M^4U^2_\sigma}, \left({1+\frac{(U_H+\delta)^2 }{2 p^2 \pi \|\gk\|^2}}\right)^{-1}, \frac{3}{4}\right\}.
\end{equation}
To guarantee \eqref{eq.u1-a}, it suffices to choose $\Psi_k$ such that
\begin{equation*}
    \frac{ \Psi_k^2 b_{n+1}^2}{\sum_{i=1}^n b_i^2} \ge \max\left\{\frac{256M^4U^2_\sigma}{\epsilon^4},1+\frac{(U_H+\delta)^2 }{2 p^2 \pi \|\gk\|^2}, \frac{4}{3}\right\}.
\end{equation*}
{Since $\sum_{i=1}^n b_i^2$ follows the chi-square distribution with $n$ degrees of freedom and $b_{n+1} \sim \mathcal{N}(0, 1)$, from a similar argument of \eqref{eq.lower of inner w.r.t. Psi}, we see that
\begin{equation}\label{eq.larger.p}
    \Psi_k^2 b_{n+1}^2 \ge \frac{\Psi_k^2p^2\pi}{2} \quad\textrm{and}\quad \sum_{i=1}^n b_i^2 \le 5n
\end{equation}
with a probability at least of $1 - \exp(-n)$. Therefore, choosing
\begin{equation*}
    \Psi_k = \frac{\sqrt{10n}}{\sqrt{\pi}p}\cdot\max\left\{\frac{16M^2U_\sigma}{\epsilon^2}, \sqrt{1+\frac{(U_H+\delta)^2 }{2 p^2 \pi \|\gk\|^2}}, \frac{2}{\sqrt{3}}\right\}
\end{equation*}
guarantees that \eqref{eq.u1-a} holds with a probability at least of $1 - \exp(-n)$.} Substituting the choice of $\Psi_k$ into \eqref{eq.lower of inner w.r.t. Psi} gives
\begin{equation}\label{eq. lower of inner final}
    (q_1^T \chi_k)^2
    \ge \min\left\{\frac{\epsilon^4}{256M^4U^2_\sigma}, \left({1+\frac{(U_H+\delta)^2 }{2 p^2 \pi \|\gk\|^2}}\right)^{-1}, \frac{3}{4}\right\} \cdot \frac{\pi^2 p^4 \sum_{i=1}^n\chi_{k,i}^2}{100n(n + 1)} + \frac{p^2\pi \chi_{k, n+1}^2}{10(n + 1)}.
\end{equation}
Combining \eqref{eq.safeguard.sigma} of Theorem \ref{thm.gamma ge delta}, we have
\begin{equation}\label{eq.upper bound of sigma}
    |\sigma_k| \le U_\sigma\sqrt{1-\alpha^2} \le \frac{\epsilon^2}{16M^2}.
\end{equation}
Finally, by the middle term in \eqref{eq.u1-a}, we have $
    \alpha^2 \ge \frac{(U_H+\delta)^2}{(U_H+\delta)^2 + 2 p^2 \pi \|\gk\|^2}.
$
Since $\gk^Tb_{[1:n]} \sim \mathcal{N}(0, \|g_k\|^2)$, from a similar argument of \eqref{eq.bound gaussian}, it holds that
\begin{equation}\label{eq.bound gb}
    (\gk^Tb_{[1:n]})^2 \ge \frac{p^2\pi\|g_k\|^2}{2}
\end{equation}
with a probability at least of $1-\exp(-n)$, implying $\alpha^2 \ge \frac{(U_H+\delta)^2}{(U_H+\delta)^2 + 4(\gk^Tb_{[1:n]})^2}$. {Recall that \eqref{eq. lower of inner final} holds with a probability at least of $1 - p - \exp(-n) - \exp(-(n+1))$ due to the union bound. Choosing $p \in (\exp(-n), 1)$ guarantees the inequalities \eqref{eq. lower of inner final} and \eqref{eq.bound gb} hold with a probability at least of $1-4p$.}
This completes the proof.

$\qed$

%% file: appendix_cutest.tex
\section{Detailed Computational Results of CUTEst Dataset}
For brevity, we use the abbreviations in \cref{tab.abbr}.
\begin{table}[hbt!]
        \centering
        \caption{Abbreviations of the Methods}\label{tab.abbr}
\begin{tabular}{ll}
                \toprule
                name          & abbreviation \\
                \arc{}        & A            \\
                \hsodm{}      & H            \\
                \hsodmhvp{}   & Hv           \\
                \newtontrst{} & N            \\
                \bottomrule
        \end{tabular}
\end{table}

\setlength{\tabcolsep}{2pt}

{
        \scriptsize
        \begin{longtable}{ll|rrrr|rrrr|rrrr}
                \caption{Complete Results on CUTEst Dataset, iteration \& time}
                \small
                \label{tab.cutest.kt}                                                                                                                                                                                                                                                                                 \\
                \toprule
                \multirow{2}{*}{name} & \multirow{2}{*}{$n$} & \multicolumn{4}{c}{$k$} & \multicolumn{4}{c}{$k^g$} & \multicolumn{4}{c}{$t$}                                                                                                                                                                          \\
                                      & {}                   & A                       & H                         & Hv                      & N              & A                 & H                & Hv               & N               & A                & H                & Hv               & N                \\
                \midrule
                \endfirsthead
                \caption[]{Complete Results on CUTEst Dataset, iteration \& time}                                                                                                                                                                                                                                     \\
                \toprule
                \multirow{2}{*}{name} & \multirow{2}{*}{$n$} & \multicolumn{4}{c}{$k$} & \multicolumn{4}{c}{$k^g$} & \multicolumn{4}{c}{$t$}                                                                                                                                                                          \\
                                      & {}                   & A                       & H                         & Hv                      & N              & A                 & H                & Hv               & N               & A                & H                & Hv               & N                \\
                \midrule
                \endhead
                \midrule
                \multicolumn{12}{r}{{Continued on next page}}                                                                                                                                                                                                                                                         \\
                \midrule
                \endfoot

                \bottomrule
                \endlastfoot
                \texttt{ARGLINA}      & \texttt{200}         & \texttt{5}              & \texttt{3}                & \texttt{3}              & \texttt{3}     & \texttt{14}       & \texttt{12}      & \texttt{16}      & \texttt{8}      & \texttt{3.5e-03} & \texttt{9.4e-01} & \texttt{3.5e+00} & \texttt{1.8e-03} \\
                \texttt{ARGLINB}      & \texttt{200}         & \texttt{3}              & \texttt{4963}             & -                       & \texttt{3}     & \texttt{8}        & \texttt{28756}   & -                & \texttt{8}      & \texttt{1.6e-03} & \texttt{2.0e+02} & \texttt{0.0e+00} & \texttt{1.6e-03} \\
                \texttt{ARGLINC}      & \texttt{200}         & \texttt{32}             & \texttt{5049}             & \texttt{5}              & \texttt{3}     & \texttt{897590}   & \texttt{55465}   & \texttt{52}      & \texttt{8}      & \texttt{2.0e+02} & \texttt{2.0e+02} & \texttt{1.2e+02} & \texttt{1.6e-03} \\
                \texttt{ARGTRIGLS}    & \texttt{200}         & \texttt{7}              & \texttt{13}               & \texttt{13}             & \texttt{7}     & \texttt{964}      & \texttt{57}      & \texttt{4006}    & \texttt{963}    & \texttt{7.4e-01} & \texttt{7.8e-01} & \texttt{3.6e+00} & \texttt{6.5e-01} \\
                \texttt{ARWHEAD}      & \texttt{1000}        & \texttt{7}              & \texttt{6}                & \texttt{6}              & \texttt{7}     & \texttt{25}       & \texttt{24}      & \texttt{41}      & \texttt{25}     & \texttt{1.2e-03} & \texttt{4.0e-03} & \texttt{3.7e-02} & \texttt{1.9e-03} \\
                                      & \texttt{100}         & \texttt{7}              & \texttt{6}                & \texttt{6}              & \texttt{7}     & \texttt{22}       & \texttt{24}      & \texttt{41}      & \texttt{22}     & \texttt{3.5e-03} & \texttt{6.1e-02} & \texttt{3.1e-02} & \texttt{4.1e-03} \\
\end{longtable}

        \begin{longtable}{rr|llll|lllllll}
                \caption{Complete Results on CUTEst Dataset, function value \& norm of the gradient}
                \scriptsize
                \label{tab.cutest.fx}                                                                                                                                                                                                    \\
                \toprule
                \multirow{2}{*}{name} & \multirow{2}{*}{$n$} & \multicolumn{4}{c}{$f$} & \multicolumn{4}{c}{$\|g\|$}                                                                                                                     \\
                                      & {}                   & A                       & H                           & Hv                & N                 & A                & H                & Hv               & N                \\
                \midrule
                \endfirsthead
                \caption[]{Complete Results on CUTEst Dataset, function value \& norm of the gradient}                                                                                                                                   \\
                \toprule
                \multirow{2}{*}{name} & \multirow{2}{*}{$n$} & \multicolumn{4}{c}{$f$} & \multicolumn{4}{c}{$\|g\|$}                                                                                                                     \\
                                      & {}                   & A                       & H                           & Hv                & N                 & A                & H                & Hv               & N                \\
                \midrule
                \endhead
                \midrule
                \multicolumn{12}{r}{{Continued on next page}}                                                                                                                                                                            \\
                \midrule
                \endfoot

                \bottomrule
                \endlastfoot
                \texttt{ARGLINA}      & \texttt{200}         & \texttt{+1.2e-22}       & \texttt{+7.0e-28}           & \texttt{+6.7e-29} & \texttt{+2.8e-26} & \texttt{2.2e-11} & \texttt{3.7e-13} & \texttt{3.1e-14} & \texttt{3.4e-13} \\
                \texttt{ARGLINB}      & \texttt{200}         & \texttt{+5.0e+01}       & \texttt{+5.0e+01}           & \texttt{+0.0e+00} & \texttt{+5.0e+01} & \texttt{1.7e-03} & \texttt{1.4e+01} & \texttt{0.0e+00} & \texttt{2.0e-03} \\
                \texttt{ARGLINC}      & \texttt{200}         & \texttt{+5.1e+01}       & \texttt{+5.1e+01}           & \texttt{+5.1e+01} & \texttt{+5.1e+01} & \texttt{1.2e+01} & \texttt{5.0e+01} & \texttt{1.8e-01} & \texttt{5.0e-04} \\
                \texttt{ARGTRIGLS}    & \texttt{200}         & \texttt{+2.1e-19}       & \texttt{+7.1e-23}           & \texttt{+6.9e-17} & \texttt{+2.2e-19} & \texttt{1.4e-08} & \texttt{2.4e-08} & \texttt{3.3e-06} & \texttt{1.5e-08} \\
                \texttt{ARWHEAD}      & \texttt{1000}        & \texttt{+0.0e+00}       & \texttt{+0.0e+00}           & \texttt{+0.0e+00} & \texttt{+0.0e+00} & \texttt{1.2e-13} & \texttt{5.1e-09} & \texttt{5.1e-09} & \texttt{1.2e-13} \\
                                      & \texttt{100}         & \texttt{+0.0e+00}       & \texttt{+0.0e+00}           & \texttt{+0.0e+00} & \texttt{+0.0e+00} & \texttt{1.1e-11} & \texttt{1.4e-07} & \texttt{1.4e-07} & \texttt{1.1e-11} \\
                \texttt{BDQRTIC}      & \texttt{1000}        & \texttt{+3.8e+02}       & \texttt{+3.8e+02}           & \texttt{+3.8e+02} & \texttt{+3.8e+02} & \texttt{1.1e-08} & \texttt{1.2e-10} & \texttt{4.7e-07} & \texttt{1.1e-08} \\
                                      & \texttt{100}         & \texttt{+4.0e+03}       & \texttt{+4.0e+03}           & \texttt{+4.0e+03} & \texttt{+4.0e+03} & \texttt{1.4e-08} & \texttt{2.9e-06} & \texttt{3.0e-06} & \texttt{1.4e-08} \\
\end{longtable}

}

%% file: homo.bib
@misc{he_homogeneous_2023,
  title      = {Homogeneous {Second}-{Order} {Descent} {Framework}: {A} {Fast} {Alternative} to {Newton}-{Type} {Methods}},
  shorttitle = {Homogeneous {Second}-{Order} {Descent} {Framework}},
  url        = {http://arxiv.org/abs/2306.17516},
  doi        = {10.48550/arXiv.2306.17516},
  urldate    = {2023-07-15},
  publisher  = {arXiv},
  author     = {He, Chang and Jiang, Yuntian and Zhang, Chuwen and Ge, Dongdong and Jiang, Bo and Ye, Yinyu},
  month      = jun,
  year       = {2023},
  note       = {arXiv:2306.17516 [math]}
}

@article{laurent2000adaptive,
	title = {Adaptive estimation of a quadratic functional by model selection},
	volume = {28},
	issn = {0090-5364},
	url = {https://projecteuclid.org/journals/annals-of-statistics/volume-28/issue-5/Adaptive-estimation-of-a-quadratic-functional-by-model-selection/10.1214/aos/1015957395.full},
	doi = {10.1214/aos/1015957395},
	number = {5},
	urldate = {2024-10-28},
	journal = {The Annals of Statistics},
	author = {Laurent, B. and Massart, P.},
	month = oct,
	year = {2000},
}

@article{curtisWorstCaseComplexityTRACE2023,
	title = {Worst-{Case} {Complexity} of {TRACE} with {Inexact} {Subproblem} {Solutions} for {Nonconvex} {Smooth} {Optimization}},
	volume = {33},
	issn = {1052-6234, 1095-7189},
	url = {https://epubs.siam.org/doi/10.1137/22M1492428},
	doi = {10.1137/22M1492428},
	language = {en},
	number = {3},
	urldate = {2023-08-22},
	journal = {SIAM Journal on Optimization},
	author = {Curtis, Frank E. and Wang, Qi},
	month = sep,
	year = {2023},
	pages = {2191--2221},
}

@misc{jin_lecture_2021,
  title  = {Lecture 14: {Lanczos} {Algorithm} ({March} 22, 2021), {ELE539}/{COS512}: Optimization for Machine Learning},
  author = {Jin, Chi},
  year   = {2021},
  institution = {Princeton University},
  url    = {https://sites.google.com/view/cjin/teaching/ece539cos512-2021-ver?authuser=0}
}

@book{nesterovLecturesConvexOptimization2018,
	address = {Cham},
	series = {Springer {Optimization} and {Its} {Applications}},
	title = {Lectures on {Convex} {Optimization}},
	volume = {137},
	isbn = {978-3-319-91577-7 978-3-319-91578-4},
	url = {http://link.springer.com/10.1007/978-3-319-91578-4},
	urldate = {2023-10-26},
	publisher = {Springer International Publishing},
	author = {Nesterov, Yurii},
	year = {2018},
	doi = {10.1007/978-3-319-91578-4},
}

@book{golubMatrixComputations2013,
	address = {Baltimore},
	edition = {Fourth edition},
	series = {Johns {Hopkins} studies in the mathematical sciences},
	title = {Matrix {Computations}},
	isbn = {978-1-4214-0794-4},
	publisher = {The Johns Hopkins University Press},
	author = {Golub, Gene H. and Van Loan, Charles F.},
	year = {2013},
}

@article{royer_newton-cg_2020,
  title    = {A {Newton}-{CG} algorithm with complexity guarantees for smooth unconstrained optimization},
  volume   = {180},
  issn     = {1436-4646},
  url      = {https://doi.org/10.1007/s10107-019-01362-7},
  doi      = {10.1007/s10107-019-01362-7},
  language = {en},
  number   = {1},
  urldate  = {2024-01-09},
  journal  = {Mathematical Programming},
  author   = {Royer, Clément W. and O’Neill, Michael and Wright, Stephen J.},
  month    = mar,
  year     = {2020},
  pages    = {451--488}
}

@article{dussaultScalableAdaptiveCubic2023,
  title    = {Scalable adaptive cubic regularization methods},
  issn     = {1436-4646},
  url      = {https://doi.org/10.1007/s10107-023-02007-6},
  doi      = {10.1007/s10107-023-02007-6},
  language = {en},
  urldate  = {2023-12-29},
  journal  = {Mathematical Programming},
  author   = {Dussault, Jean-Pierre and Migot, Tangi and Orban, Dominique},
  month    = oct,
  year     = {2023}
}

@misc{haegeman_krylovkit_2024,
  title     = {{KrylovKit}},
  url       = {https://zenodo.org/records/10884302},
  urldate   = {2024-04-21},
  publisher = {Zenodo},
  author    = {Haegeman, Jutho},
  month     = mar,
  year      = {2024},
  doi       = {10.5281/zenodo.10884302}
}

@article{kmogensenOptimMathematicalOptimization2018,
  title      = {Optim: {A} mathematical optimization package for {Julia}},
  volume     = {3},
  issn       = {2475-9066},
  shorttitle = {Optim},
  url        = {http://joss.theoj.org/papers/10.21105/joss.00615},
  doi        = {10.21105/joss.00615},
  number     = {24},
  urldate    = {2023-12-29},
  journal    = {Journal of Open Source Software},
  author     = {K Mogensen, Patrick and N Riseth, Asbjørn},
  month      = apr,
  year       = {2018},
  pages      = {615}
}

@misc{orban_juliasmoothoptimizers_2019,
  title     = {{JuliaSmoothOptimizers}},
  url       = {https://zenodo.org/record/2655082},
  urldate   = {2023-04-02},
  publisher = {Zenodo},
  author    = {Orban, Dominique and Siqueira, Abel},
  month     = apr,
  year      = {2019},
  doi       = {10.5281/zenodo.2655082},
  note      = {Language: eng}
}

@book{cartis_evaluation_2022,
  address    = {Philadelphia, PA},
  title      = {Evaluation {Complexity} of {Algorithms} for {Nonconvex} {Optimization}: {Theory}, {Computation} and {Perspectives}},
  isbn       = {978-1-61197-698-4 978-1-61197-699-1},
  shorttitle = {Evaluation {Complexity} of {Algorithms} for {Nonconvex} {Optimization}},
  url        = {https://epubs.siam.org/doi/book/10.1137/1.9781611976991},
  language   = {en},
  urldate    = {2022-12-29},
  publisher  = {Society for Industrial and Applied Mathematics},
  author     = {Cartis, Coralia and Gould, Nicholas I. M. and Toint, Philippe L.},
  month      = jan,
  year       = {2022},
  doi        = {10.1137/1.9781611976991}
}

@article{dolan_benchmarking_2002,
  title    = {Benchmarking optimization software with performance profiles},
  volume   = {91},
  issn     = {1436-4646},
  url      = {https://doi.org/10.1007/s101070100263},
  doi      = {10.1007/s101070100263},
  language = {en},
  number   = {2},
  urldate  = {2023-04-30},
  journal  = {Mathematical Programming},
  author   = {Dolan, Elizabeth D. and Moré, Jorge J.},
  month    = jan,
  year     = {2002},
  pages    = {201--213}
}

@article{gratton2023yet,
  title   = {Yet another fast variant of Newton's method for nonconvex optimization},
  author  = {Gratton, Serge and Jerad, Sadok and Toint, Philippe L},
  journal = {arXiv preprint arXiv:2302.10065},
  year    = {2023}
}

@article{mishchenko2021regularized,
  title   = {Regularized Newton Method with Global  {$O(1/k^2)$} Convergence},
  author  = {Mishchenko, Konstantin},
  journal = {arXiv preprint arXiv:2112.02089},
  year    = {2021}
}

@inproceedings{agarwal_finding_2017,
  title     = {Finding approximate local minima faster than gradient descent},
  booktitle = {Proceedings of the 49th {Annual} {ACM} {SIGACT} {Symposium} on {Theory} of {Computing}},
  author    = {Agarwal, Naman and Allen-Zhu, Zeyuan and Bullins, Brian and Hazan, Elad and Ma, Tengyu},
  year      = {2017},
  pages     = {1195--1199}
}

@article{carmon_accelerated_2018,
  title    = {Accelerated {Methods} for {NonConvex} {Optimization}},
  volume   = {28},
  issn     = {1052-6234, 1095-7189},
  url      = {https://epubs.siam.org/doi/10.1137/17M1114296},
  doi      = {10.1137/17M1114296},
  language = {en},
  number   = {2},
  urldate  = {2022-07-05},
  journal  = {SIAM Journal on Optimization},
  author   = {Carmon, Yair and Duchi, John C. and Hinder, Oliver and Sidford, Aaron},
  month    = jan,
  year     = {2018},
  pages    = {1751--1772}
}

@article{cartis_adaptive_2011,
  title      = {Adaptive cubic regularisation methods for unconstrained optimization. {Part} {I}: motivation, convergence and numerical results},
  volume     = {127},
  issn       = {0025-5610, 1436-4646},
  shorttitle = {Adaptive cubic regularisation methods for unconstrained optimization. {Part} {I}},
  url        = {http://link.springer.com/10.1007/s10107-009-0286-5},
  doi        = {10.1007/s10107-009-0286-5},
  language   = {en},
  number     = {2},
  urldate    = {2022-05-25},
  journal    = {Mathematical Programming},
  author     = {Cartis, Coralia and Gould, Nicholas I. M. and Toint, Philippe L.},
  month      = apr,
  year       = {2011},
  pages      = {245--295}
}

@article{cartis_adaptive_2011-1,
  title      = {Adaptive cubic regularisation methods for unconstrained optimization. {Part} {II}: worst-case function- and derivative-evaluation complexity},
  volume     = {130},
  issn       = {0025-5610, 1436-4646},
  shorttitle = {Adaptive cubic regularisation methods for unconstrained optimization. {Part} {II}},
  url        = {http://link.springer.com/10.1007/s10107-009-0337-y},
  doi        = {10.1007/s10107-009-0337-y},
  language   = {en},
  number     = {2},
  urldate    = {2022-05-25},
  journal    = {Mathematical Programming},
  author     = {Cartis, Coralia and Gould, Nicholas I. M. and Toint, Philippe L.},
  month      = dec,
  year       = {2011},
  pages      = {295--319}
}

@article{cartis_complexity_2010,
  title    = {On the {Complexity} of {Steepest} {Descent}, {Newton}'s and {Regularized} {Newton}'s {Methods} for {Nonconvex} {Unconstrained} {Optimization} {Problems}},
  volume   = {20},
  issn     = {1052-6234},
  url      = {https://epubs.siam.org/doi/abs/10.1137/090774100},
  doi      = {10.1137/090774100},
  number   = {6},
  urldate  = {2022-06-28},
  journal  = {SIAM Journal on Optimization},
  author   = {Cartis, Coralia and Gould, Nicholas I. M. and Toint, Philippe L.},
  month    = jan,
  year     = {2010},
  note     = {Publisher: Society for Industrial and Applied Mathematics},
  pages    = {2833--2852}
}

@book{connTrustRegionMethods2000,
	address = {Philadelphia, Pa},
	series = {{MPS}-{SIAM} series on optimization},
	title = {Trust-{Region} {Methods}},
	isbn = {978-0-89871-460-9},
	language = {eng},
	publisher = {Society for Industrial and Applied Mathematics},
	author = {Conn, Andrew R. and Gould, Nicholas I. M. and Toint, Philippe L.},
	collaborator = {{Society for Industrial and Applied Mathematics}},
	year = {2000},
	doi = {10.1137/1.9780898719857},
}

@article{curtis_concise_2018,
  title    = {Concise complexity analyses for trust region methods},
  volume   = {12},
  issn     = {1862-4480},
  url      = {https://doi.org/10.1007/s11590-018-1286-2},
  doi      = {10.1007/s11590-018-1286-2},
  language = {en},
  number   = {8},
  urldate  = {2022-06-30},
  journal  = {Optimization Letters},
  author   = {Curtis, Frank E. and Lubberts, Zachary and Robinson, Daniel P.},
  month    = dec,
  year     = {2018},
  pages    = {1713--1724}
}

@article{curtis_exploiting_2019,
  title    = {Exploiting negative curvature in deterministic and stochastic optimization},
  volume   = {176},
  issn     = {1436-4646},
  url      = {https://doi.org/10.1007/s10107-018-1335-8},
  doi      = {10.1007/s10107-018-1335-8},
  language = {en},
  number   = {1},
  urldate  = {2022-08-28},
  journal  = {Mathematical Programming},
  author   = {Curtis, Frank E. and Robinson, Daniel P.},
  month    = jul,
  year     = {2019},
  pages    = {69--94}
}

@article{curtis_trust_2017,
  title    = {A trust region algorithm with a worst-case iteration complexity of {$O(\epsilon^{-3/2})$} for nonconvex optimization},
  volume   = {162},
  issn     = {1436-4646},
  url      = {https://doi.org/10.1007/s10107-016-1026-2},
  doi      = {10.1007/s10107-016-1026-2},
  language = {en},
  number   = {1},
  urldate  = {2022-06-30},
  journal  = {Mathematical Programming},
  author   = {Curtis, Frank E. and Robinson, Daniel P. and Samadi, Mohammadreza},
  month    = mar,
  year     = {2017},
  pages    = {1--32}
}

@article{gould_cutest_2015,
  title      = {{CUTEst}: a {Constrained} and {Unconstrained} {Testing} {Environment} with safe threads for mathematical optimization},
  volume     = {60},
  issn       = {1573-2894},
  shorttitle = {{CUTEst}},
  url        = {https://doi.org/10.1007/s10589-014-9687-3},
  doi        = {10.1007/s10589-014-9687-3},
  language   = {en},
  number     = {3},
  urldate    = {2022-08-31},
  journal    = {Computational Optimization and Applications},
  author     = {Gould, Nicholas I. M. and Orban, Dominique and Toint, Philippe L.},
  month      = apr,
  year       = {2015},
  pages      = {545--557}
}

@article{hager_algorithm_2006,
  title      = {Algorithm 851: {CG}\_DESCENT, a conjugate gradient method with guaranteed descent},
  volume     = {32},
  shorttitle = {Algorithm 851},
  number     = {1},
  journal    = {ACM Transactions on Mathematical Software (TOMS)},
  author     = {Hager, William W. and Zhang, Hongchao},
  year       = {2006},
  note       = {Publisher: ACM New York, NY, USA},
  pages      = {113--137}
}

@inproceedings{jin_how_2017,
  title     = {How to escape saddle points efficiently},
  booktitle = {International {Conference} on {Machine} {Learning}},
  publisher = {PMLR},
  author    = {Jin, Chi and Ge, Rong and Netrapalli, Praneeth and Kakade, Sham M. and Jordan, Michael I.},
  year      = {2017},
  pages     = {1724--1732}
}

@article{kuczynski_estimating_1992,
  title    = {Estimating the {Largest} {Eigenvalue} by the {Power} and {Lanczos} {Algorithms} with a {Random} {Start}},
  volume   = {13},
  issn     = {0895-4798, 1095-7162},
  url      = {http://epubs.siam.org/doi/10.1137/0613066},
  doi      = {10.1137/0613066},
  language = {en},
  number   = {4},
  urldate  = {2022-10-13},
  journal  = {SIAM Journal on Matrix Analysis and Applications},
  author   = {Kuczyński, J. and Woźniakowski, H.},
  month    = oct,
  year     = {1992},
  pages    = {1094--1122}
}

@misc{li_restarted_2022,
  title      = {Restarted {Nonconvex} {Accelerated} {Gradient} {Descent}: {No} {More} {Polylogarithmic} {Factor} in the {$O(\epsilon^{-7/4})$} {Complexity}},
  shorttitle = {Restarted {Nonconvex} {Accelerated} {Gradient} {Descent}},
  url        = {http://arxiv.org/abs/2201.11411},
  doi        = {10.48550/arXiv.2201.11411},
  urldate    = {2022-11-14},
  publisher  = {arXiv},
  author     = {Li, Huan and Lin, Zhouchen},
  month      = may,
  year       = {2022},
  note       = {arXiv:2201.11411 [cs, math]}
}

@book{luenberger_linear_2021,
  address   = {Cham},
  series    = {International {Series} in {Operations} {Research} \& {Management} {Science}},
  title     = {Linear and {Nonlinear} {Programming}},
  volume    = {228},
  isbn      = {978-3-030-85449-2 978-3-030-85450-8},
  url       = {https://link.springer.com/10.1007/978-3-030-85450-8},
  language  = {en},
  urldate   = {2022-02-15},
  publisher = {Springer International Publishing},
  author    = {Luenberger, David G. and Ye, Yinyu},
  year      = {2021},
  doi       = {10.1007/978-3-030-85450-8}
}

@article{nesterov_cubic_2006,
  title    = {Cubic regularization of {Newton} method and its global performance},
  volume   = {108},
  issn     = {1436-4646},
  url      = {https://doi.org/10.1007/s10107-006-0706-8},
  doi      = {10.1007/s10107-006-0706-8},
  language = {en},
  number   = {1},
  urldate  = {2022-06-29},
  journal  = {Mathematical Programming},
  author   = {Nesterov, Yurii and Polyak, B.T.},
  month    = aug,
  year     = {2006},
  pages    = {177--205}
}

@book{nocedal_numerical_2006,
	address = {New York, NY},
	edition = {Second},
	series = {Springer series in operations research and financial engineering},
	title = {{Numerical} {Optimization}},
	isbn = {978-0-387-30303-1},
	language = {eng},
	publisher = {Springer},
	author = {Nocedal, Jorge and Wright, Stephen J.},
	year = {2006},
}

@book{parlett_symmetric_1998,
  title     = {The {Symmetric} {Eigenvalue} {Problem}},
  url       = {http://epubs.siam.org/doi/book/10.1137/1.9781611971163},
  language  = {en},
  urldate   = {2024-04-21},
  publisher = {Society for Industrial and Applied Mathematics},
  author    = {Parlett, Beresford N.},
  month     = jan,
  year      = {1998},
  doi       = {10.1137/1.9781611971163}
}

@article{rojas_new_2001,
  title    = {A {New} {Matrix}-{Free} {Algorithm} for the {Large}-{Scale} {Trust}-{Region} {Subproblem}},
  volume   = {11},
  issn     = {1052-6234},
  url      = {https://epubs.siam.org/doi/abs/10.1137/S105262349928887X},
  doi      = {10.1137/S105262349928887X},
  number   = {3},
  urldate  = {2022-10-23},
  journal  = {SIAM Journal on Optimization},
  author   = {Rojas, Marielba and Santos, Sandra A. and Sorensen, Danny C.},
  month    = jan,
  year     = {2001},
  note     = {Publisher: Society for Industrial and Applied Mathematics},
  pages    = {611--646}
}

@article{royer_complexity_2018,
  title   = {Complexity analysis of second-order line-search algorithms for smooth nonconvex optimization},
  volume  = {28},
  number  = {2},
  journal = {SIAM Journal on Optimization},
  author  = {Royer, Clément W. and Wright, Stephen J.},
  year    = {2018},
  note    = {Publisher: SIAM},
  pages   = {1448--1477}
}

@article{sturm_cones_2003,
	title = {On {Cones} of {Nonnegative} {Quadratic} {Functions}},
	volume = {28},
	issn = {0364-765X},
	url = {https://pubsonline.informs.org/doi/abs/10.1287/moor.28.2.246.14485},
	doi = {10.1287/moor.28.2.246.14485},
	number = {2},
	urldate = {2024-10-23},
	journal = {Mathematics of Operations Research},
	author = {Sturm, Jos F. and Zhang, Shuzhong},
	month = may,
	year = {2003},
	note = {Publisher: INFORMS},
	pages = {246--267},
}

@article{xu_first-order_2018,
  title   = {First-order stochastic algorithms for escaping from saddle points in almost linear time},
  volume  = {31},
  journal = {Advances in neural information processing systems},
  author  = {Xu, Yi and Jin, Rong and Yang, Tianbao},
  year    = {2018}
}

@article{ye_new_2003,
  title   = {New results on quadratic minimization},
  volume  = {14},
  number  = {1},
  journal = {SIAM Journal on Optimization},
  author  = {Ye, Yinyu and Zhang, Shuzhong},
  year    = {2003},
  note    = {Publisher: SIAM},
  pages   = {245--267}
}

@misc{yeSecondOrderOptimization2005,
	title = {Second {Order} {Optimization} {Algorithms} {I}},
	url = {https://web.stanford.edu/class/msande311/lecture12.pdf},
	language = {en},
	author = {Ye, Yinyu},
	year = {2005},
}

@misc{zhang_drsom_2023,
  title      = {{DRSOM}: {A} {Dimension} {Reduced} {Second}-{Order} {Method}},
  shorttitle = {{DRSOM}},
  url        = {http://arxiv.org/abs/2208.00208},
  doi        = {10.48550/arXiv.2208.00208},
  urldate    = {2023-01-03},
  publisher  = {arXiv},
  author     = {Zhang, Chuwen and Ge, Dongdong and He, Chang and Jiang, Bo and Jiang, Yuntian and Ye, Yinyu},
  month      = jan,
  year       = {2022},
  note       = {arXiv:2208.00208 [cs, math]}
}

@article{Julia-2017,
  title     = {Julia: A fresh approach to numerical computing},
  author    = {Bezanson, Jeff and Edelman, Alan and Karpinski, Stefan and Shah, Viral B},
  journal   = {SIAM {R}eview},
  volume    = {59},
  number    = {1},
  pages     = {65--98},
  year      = {2017},
  publisher = {SIAM},
  doi       = {10.1137/141000671},
  url       = {https://epubs.siam.org/doi/10.1137/141000671}
}

@article{curtis_inexact_2018,
  author  = {Curtis, Frank E and Robinson, Daniel P and Samadi, Mohammadreza},
  title   = {An inexact regularized Newton framework with a worst-case iteration complexity of $\mathcal{O}(\varepsilon^{-3/2})$ for nonconvex optimization},
  journal = {IMA Journal of Numerical Analysis},
  volume  = {39},
  number  = {3},
  pages   = {1296-1327},
  year    = {2018},
  month   = {05}
}

@article{curtis2021newton-cg,
	title = {Trust-{Region} {Newton}-{CG} with {Strong} {Second}-{Order} {Complexity} {Guarantees} for {Nonconvex} {Optimization}},
	volume = {31},
	issn = {1052-6234},
	url = {https://epubs.siam.org/doi/10.1137/19M130563X},
	doi = {10.1137/19M130563X},
	number = {1},
	urldate = {2023-05-25},
	journal = {SIAM Journal on Optimization},
	author = {Curtis, Frank E. and Robinson, Daniel P. and Royer, Clément W. and Wright, Stephen J.},
	month = jan,
	year = {2021},
	note = {Publisher: Society for Industrial and Applied Mathematics},
	pages = {518--544},
}

@article{generalized-trs,
  author  = {Adachi, Satoru and Iwata, Satoru and Nakatsukasa, Yuji and Takeda, Akiko},
  title   = {Solving the Trust-Region Subproblem By a Generalized Eigenvalue Problem},
  journal = {SIAM Journal on Optimization},
  volume  = {27},
  number  = {1},
  pages   = {269-291},
  year    = {2017}
}

@article{felix-arc-generalized,
  author  = {Lieder, Felix},
  title   = {Solving Large-Scale Cubic Regularization by a Generalized Eigenvalue Problem},
  journal = {SIAM Journal on Optimization},
  volume  = {30},
  number  = {4},
  pages   = {3345-3358},
  year    = {2020}
}
